\documentclass[10pt]{amsart}
\usepackage{amsmath,amsfonts,amscd,amssymb,graphicx,mathrsfs,eufrak}
\DeclareSymbolFont{bbold}{U}{bbold}{m}{n}
\DeclareSymbolFontAlphabet{\mathbbold}{bbold}

\usepackage[dvips,all,arc,curve,color,frame]{xy}
\newxyColor{pink}{1.0 0.4 0.5}{rgb}{}
\usepackage[usenames]{color}
\newtheorem{defin}{Definition}[section]
\newtheorem{thm}[defin]{Theorem}
\newtheorem{prop}[defin]{Proposition}
\newtheorem{lemma}[defin]{Lemma}

\newtheorem{notation}[defin]{Notation}
\newtheorem{example}[defin]{Example}
\newtheorem{remark}[defin]{Remark}

\newcommand{\Mod}[1]{\mathfrak{Mod}#1}

\newcommand{\stq}[1]{\underline{\Mod{#1}}}

\newcommand{\Sstq}[1]{\textsf{St}(\stq{R})}

\newcommand{\C}[1]{\mathbb{C}^{#1}}

\newcommand{\s}[1]{\mathscr{#1}}
\renewcommand{\c}[1]{\mathcal{#1}}

\newcommand{\cl}[2]{c_{#1 #2}}
\newcommand{\an}[2]{a_{#1 #2}}
\newcommand{\CL}[2]{C_{#1 #2}}
\newcommand{\CLr}[2]{{{C_{#1 #2}}}}
\newcommand{\CLm}[2]{\hat{C}_{#1 #2}}
\newcommand{\AN}[2]{A_{#1 #2}}

\newcommand{\w}[1]{\widetilde{#1}}
\newcommand{\e}{\varepsilon}

\renewcommand{\t}[1]{\textnormal{#1}}
\renewcommand{\tt}[1]{\mathtt{#1}}

\def\GL{\mathop{\rm GL}\nolimits}
\def\SL{\mathop{\rm SL}\nolimits}
\def\CM{\mathop{\rm CM}\nolimits}
\def\Ext{\mathop{\rm Ext}\nolimits}
\def\Spec{\mathop{\rm Spec}\nolimits}
\def\proj{\mathop{\rm proj}\nolimits}
\def\rk{\mathop{\rm rank}\nolimits}
\begin{document}
\title{\textsc{Reconstruction Algebras of Type $D$ (I)}}
\author{Michael Wemyss}
\address{Michael Wemyss, The Maxwell Institute, School of Mathematics, James Clerk Maxwell Building, The King's Buildings, Mayfield Road, Edinburgh, EH9 3JZ, UK.}
\email{wemyss.m@googlemail.com}
\begin{abstract}
This is the second in a series of papers which give an explicit description of the reconstruction algebra as a quiver with relations; these algebras arise naturally as geometric generalizations of preprojective algebras of extended Dynkin diagrams.  This paper deals with dihedral groups $G=\mathbb{D}_{n,q}$ for which all special CM modules have rank one, and we show that all but four of the relations on such a reconstruction algebra are given simply as the relations arising from a reconstruction algebra of type $A$.  As a corollary, the reconstruction algebra reduces the problem of explicitly understanding the minimal resolution (=G-Hilb) to the same level of difficulty as the toric case. 
\end{abstract}
\maketitle
\parindent 20pt
\parskip 0pt

\tableofcontents

\section{Introduction}
It was discovered in \cite{Wemyss_reconstruct_A} that if $G$ is a finite small cyclic subgroup of $\GL(2,\C{})$ then the quiver of the endomorphism ring of the special CM $\C{}[x,y]^G$ modules (a ring built downstairs on the singularity) determines and is determined by the dual graph of  the minimal resolution $\w{X}\rightarrow\mathbb{C}^2/G$, labelled with self-intersection numbers.  Since the dual graph of such a group is always a Dynkin diagram of type $A$, we call the noncommutative ring in question the reconstruction algebra of type $A$.  The above is a correspondence purely on the level of the underlying quivers; it was further discovered that if we add in the extra information of the relations then in fact one can recover the whole space $\w{X}$ (not just the dual graph) as a certain GIT quotient, and also that the reconstruction algebra describes the derived category of $\w{X}$. 

Following the work of Bridgeland \cite{Bridgeland} and Van den Bergh \cite{VdB}, these ideas were pursued further in \cite{Wemyss_GL2} where the above statement on the level of quivers was proved for all complex rational surface singularities.  Two proofs of this fact were given, one non-explicit proof for the general case and one explicit proof which only covers the quotient case.  It is perhaps important to emphasize two points.  First, although it was shown that the number of relations for the quiver can also be obtained from the intersection theory, the relations themselves were not exhibited. Second, the non-explicit proof tells us nothing about the special CM modules (for example what they are) and very little about their structure.

The main purpose of this paper is to provide the relations in the case of certain dihedral groups $\mathbb{D}_{n,q}$ inside $\GL(2,\C{})$.  The companion paper \cite{Wemyss_reconstruct_D(ii)} deals with the remaining dihedral cases.  Although this may be technical the relations are important; we shall see that it is precisely the relations which allow us to compare the geometry of the minimal resolutions of different singularities and thus view the minimal resolutions as being very similar to spaces that we already understand.  We could not make such claims if we do not know the relations.  

It is worth emphasizing how small the reconstruction algebra is compared to the skew group ring $\C{}[x,y]\# \mathbb{D}_{n,q}$.  Although many statements can be made about the latter ring, the problem has always been extracting the geometric information that it encodes.  Here this is precisely because it contains far too much irrelevant information; indeed we believe it is entirely the wrong algebra to study.  Note that the reconstruction algebra and the skew group ring coincide if and only if the group is inside $\SL(2,\C{})$.  

The basic philosophy underpinning the geometry of the reconstruction algebra is that we should not view the minimal resolution of an affine rational surface as some sort of exotic beast, rather we should instead view it as being made from spaces that we are already familiar with.  The larger the fundamental cycle $Z_f$ of the minimal resolution the more the geometry should resemble the minimal resolutions arising from ADE quotients, since the reconstruction algebra quiver and relations look and behave more like a preprojective algebra (see \cite[\S 5]{Wemyss_reconstruct_D(ii)} for this case).  As $Z_f$ gets smaller and thus closer to being reduced (the case mainly considered in this paper), the more toric the geometry becomes since the reconstruction algebra quiver and relations begin to look and behave more like a reconstruction algebra of type $A$.  Via the work of Wunram it is $Z_{f}$ which dictates the rank of the special CM modules, and different ranks induce slightly different algebra structures because polynomials factor in different ways. 

The first step in this paper towards the goal of obtaining the relations is to fill in the gap in the classification of the specials CM modules left in \cite{Iyama_Wemyss_specials}.   Although the techniques in \emph{loc.\ cit} could plausibly be used for dihedral groups, the AR quiver splits into cases and the combinatorics are difficult. Both factors make it hard to write down a general proof.  The method used here is still quite combinatorially complex but it is at least possible to write down the proof.  It also has the added philosophical benefit of not requiring any knowledge of the McKay quiver, as instead (denoting $R:=\C{}[x,y]^{\mathbb{D}_{n,q}}$) we view the special CM $R$-modules as $R$-summands of the polynomial ring \cite{Herzog}, and by considering an appropriate quiver we argue by diagram chasing that the modules corresponding to the vertices must be minimally two-generated.  Furthermore we explicitly obtain their generators.  In fact this method gives all the rank one special CM modules for every $\mathbb{D}_{n,q}$, not just the case when $Z_f$ is reduced.

In the case when $Z_f$ is reduced, all special CM modules have rank one and so the special CM modules obtained above are them all.  Hence having explicitly obtained their generators, we then use this information to label the arrows in the known quiver of the reconstruction algebra in terms of polynomials involving $x$'s and $y$'s; our relations are then simply that `$x$ and $y$ commute'.  From this we write down some obvious relations and then an easy argument using the known number of relations tells us that these are them all.  Since there are choices for the generators of the special CM modules, in fact we obtain two different (but isomorphic) presentations of the endomorphism ring, corresponding to two different choices of such generators.

Once we have the quiver with relations we are able to describe the moduli spaces of representations using explicit techniques, justifying some of the philosophy above.  This is actually very easy since the dimension vector we use consists entirely of 1's and (almost all) the relations for the reconstruction algebras in this paper are just the relations from the reconstruction algebra of type $A$.  Consequently the explicit extraction of the geometry is really just the same as in \cite[\S 4]{Wemyss_reconstruct_A}. 

We now describe the contents and structure of this paper in more detail:  in Section 2 we define the groups $\mathbb{D}_{n,q}$ and recall some of the known properties of the singularities $\C{2}/\mathbb{D}_{n,q}$.  Furthermore we introduce and prove some combinatorics crucial to later arguments. The combinatorics are intricate, and can be skipped on first reading.  In Section 3 we exhibit the special rank one CM modules for every group $\mathbb{D}_{n,q}$ regardless of $Z_{f}$; when $Z_{f}$ is reduced these are all the special CM modules. Section 4 restricts to the case when $Z_f$ is reduced and in such a case we define $D_{n,q}$, the associated \emph{reconstruction algebra of type $D$}, and show that it is isomorphic to the endomorphism ring of the special CM modules.  In Section 5 we use this noncommutative algebra to exhibit explicitly the minimal resolution (which is equal to $\mathbb{D}_{n,q}$-Hilb via a result of Ishii \cite{Ishii}) in co-ordinates.  We produce an open affine cover in which every affine open is a smooth hypersurface in $\C{3}$, with equations given in terms of simple combinatorics.

Note that the discovery of the special CM modules and a similar open cover for dihedral groups has been discovered independently by Nolla de Celis \cite{Alvaro1},\cite{Alvaro2} by using the McKay quiver and combinatorics of $G$-Hilb.  However the benefits of using the reconstruction algebra over the McKay quiver is that it reduces the problem to the same level of difficulty as the toric case; thus from the viewpoint of the reconstruction algebra the geometry in this paper is not toric, but it may as well be.

Throughout, when working with quivers we shall write $ab$ to mean \textbf{$a$ followed by $b$.}  We work over the ground field $\mathbb{C}$ but any algebraically closed field of characteristic zero will suffice.

Various parts of this work were carried out when the author was a PhD student in Bristol, UK (funded by the EPSRC), when the author visited Nagoya University under the Cecil King Travel Scholarship from the London Mathematical Society, and also when the author was in receipt of a JSPS Postdoctoral Fellowship.  The author would like to thank the EPSRC, the Cecil King Foundation, the LMS and the JSPS, and also thank the University of Nagoya for kind hospitality.  Thanks are also due to Osamu Iyama, Martin Herschend and Alvaro Nolla de Celis for useful discussions, and an anonymous referee for some extremely helpful remarks and suggestions regarding the presentation of mathematics.

\section{Dihedral groups and combinatorics}
In this paper we follow the notation of Riemenschneider
\cite{Riemenschneider_invarianten}.
\begin{defin}
For $1<q<n$ with $(n,q)=1$ define the group $\mathbb{D}_{n,q}$ to
be
\[
\begin{array}{cc}
\mathbb{D}_{n,q}=\left\{ \begin{array}{cc} \langle \psi_{2q},
\tau, \varphi_{2(n-q)} \rangle& \mbox{if } n-q\equiv 1 \mbox{ mod
}2\\\langle \psi_{2q}, \tau\varphi_{4(n-q)} \rangle& \mbox{if }
n-q\equiv 0 \mbox{ mod }2
\end{array}\right.
\end{array}
\]
with the matrices
{\scriptsize{
\[
\begin{array}{ccc}
\psi_k=\left( \begin{array}{cc}\e_k & 0\\ 0& \e_k^{-1}
\end{array}\right) &\tau=\left( \begin{array}{cc}0 & \e_4\\ \e_4& 0
\end{array}\right)&\varphi_k=\left( \begin{array}{cc}\e_k & 0\\ 0& \e_k
\end{array}\right)
\end{array}
\]
}}where $\e_t$ is a primitive $t^{th}$ root of unity.
\end{defin}
The order of the group $\mathbb{D}_{n,q}$ is $4(n-q)q$.  The procedure to obtain the invariants
$\C{}[x,y]^{\mathbb{D}_{n,q}}$ is also well known: first develop
\[
\frac{n}{n-q}=[a_2,\hdots,a_{e-1} ]=a_2-\frac{1}{a_3 - \frac{1}{a_4 -
\frac{1}{(...)}}}
\]
as a Jung--Hirzebruch continued fraction expansion.  We fix this notation throughout the paper, noting the strange numbering.  Now define %\footnote{The last value is $c_{e}$ where $e=\sum_{i=1}^{N}(\alpha_1-2)+3$; if all $\alpha$'s are 2 (ie $SL(2,\C{})$) only get to 3.} 
series $c_j$, $d_j$ and $t_j$ for $2\leq j\leq e$ by
\[
\begin{array}{ll|ll}
c_2=1 & c_3=0 & c_4=1 & c_{j}=a_{j-1}c_{j-1}-c_{j-2} \mbox{ for
}5\leq j \leq e\\
d_2=0 & d_3=1 & d_4=a_3-1 & d_{j}=a_{j-1}d_{j-1}-d_{j-2} \mbox{
for
}5\leq j \leq e\\
t_2=a_2 & t_3=a_2-1 & t_4=a_3(a_2-1)-1 &
t_{j}=a_{j-1}t_{j-1}-t_{j-2} \mbox{ for
}5\leq j \leq e\\
\end{array}
\]
where the values to the right of the line exist only when $e>3$, i.e.\
when $n>q+1$, i.e. when the group $\mathbb{D}_{n,q}$ does not lie
inside $\SL(2,\C{})$. Also define the series $r_j$ for $2\leq j\leq e$
by
\[
r_j=(n-q)t_j-q(c_j+d_j) \mbox{ for } 2\leq j\leq e
\]
Note by definition that $r_j=a_{j-1}r_{j-1}-r_{j-2}$ for all
$5\leq j\leq e$.  Throughout this paper we set
\[
\begin{array}{ccc}
w_1:=xy & v_2:=x^{2q}+(-1)^{a_2}y^{2q} &
v_3:=x^{2q}+(-1)^{a_2-1}y^{2q}.
\end{array}
\]
\begin{thm}\cite[Satz 2]{Riemenschneider_invarianten}\label{Riemen_generates} The
polynomials $w_1^{2(n-q)}$ and
$w_1^{r_t}v_2^{c_t}v_3^{d_t}$ for $2\leq t\leq e$ generate the
ring $\C{}[x,y]^{\mathbb{D}_{n,q}}$.
\end{thm}
The main ingredient in the proof is Noether's bound on the degree
of the generators in characteristic zero; once this is used the
proof reduces to combinatorics.  In this paper we
shall need the following easy variant of the above: define
\begin{eqnarray*}
w_2&=&(x^q+y^q)(x^q+(-1)^{a_2}y^q)\\
w_3&=&(x^q-y^q)(x^q+(-1)^{a_2}y^q)
\end{eqnarray*}
\begin{lemma}\label{invariants}
The polynomials $w_1^{2(n-q)}$ and
$w_1^{r_t}w_2^{c_t}w_3^{d_t}$ for $2\leq t\leq e$ generate the
ring $\C{}[x,y]^{\mathbb{D}_{n,q}}$.%\footnote{I don't think we use the `cross-terms' version anywhere except in the examples section at the end: there we are ok by non-explicit GL since we know by that stage the arrows in the examples are all the irreducible maps between the specials.} \footnote{We do however need this proof to be independent of $\nu$, since we use it to get the 1d specials for all $\nu$} \footnote{Prove using the Hilbert polynomial - each of these guys is invariant, sitting in the degree where the generators are}
\end{lemma}
For $\mathbb{D}_{n,q}$, throughout this paper we fix the notation
\[
\frac{n}{q}=[\alpha_1,\hdots,\alpha_N]
\]
as the Hirzebruch--Jung continued fraction expansion of $\frac{n}{q}$.  By Riemenschneider duality
(see e.g. \cite[1.2]{Kidoh} or \cite[2.1]{Wemyss_reconstruct_A}) it is true that $e-2=1+\sum_{i=1}^{N}(\alpha_i-2)$.
\begin{defin}
Define $\nu<N$ to be the largest integer such that
$\alpha_1=\hdots=\alpha_\nu=2$, or 0 if no such integer exists.
\end{defin}
Now by \cite[2.11]{Brieskorn} the dual graph of the minimal resolution
of $\C{2}/\mathbb{D}_{n,q}$ is
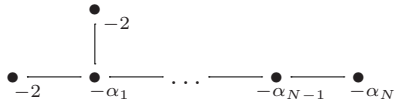
\[
\xymatrix@C=20pt@R=15pt{ &\bullet\ar@{-}[d]^<{-2}&&&\\
\bullet\ar@{-}[r]_<{-2} & \bullet\ar@{-}[r]_<{-\alpha_{1}}
&\hdots\ar@{-}[r] &\bullet\ar@{-}[r]_<{-\alpha_{N-1}} & \bullet
\ar@{}[r]_<{-\alpha_N}&}
\]
where the $\alpha$'s come from the Jung--Hirzebruch continued fraction expansion of $\frac{n}{q}$ above.  Notice that the only possible fundamental cycles $Z_{f}$ for dihedral groups
$\mathbb{D}_{n,q}$ are
\[
\begin{array}{c}\xymatrix@R=2pt@C=2pt{ &1&&&\\
1& 1&1&\hdots&1&1}\end{array} \quad\mbox{or}\quad
\begin{array}{c}\xymatrix@R=2pt@C=2pt{ &1&&&\\
1& 2&\hdots&2&1&\hdots&1}\end{array}
\]
where (since by definition $\nu<N$) the number of 2's in the right-hand picture is precisely $\nu$.  Thus $\nu$ records the number of 2's in $Z_f$.  
\begin{defin}\cite{Wunram_generalpaper}
Denote $R=\C{}[x,y]^{\mathbb{D}_{n,q}}$.  A CM $R$-module $M$ is called special if $(M\otimes \omega_R)/\t{tors} $ is CM, where $\omega_R$ is the canonical module of $R$.
\end{defin}
There are now many equivalent characterizations of the special CM $R$-modules, for example they are precisely the CM $R$-modules $X$ for which $\Ext_R^1(X,R)=0$ \cite[2.7]{Iyama_Wemyss_specials}.  The theory of special CM modules was first developed in the work of Wunram \cite{Wunram_generalpaper}, who proved the following results:

\begin{thm}\label{WunramMainResults}  Let $X=\Spec R$ be a complete--local rational normal surface singularity, and denote the minimal resolution by $Y\to \Spec R$.\\
(1) There is a 1--1 correspondence between the non-free indecomposable special CM $R$-modules and the exceptional curves in the minimal resolution $Y$.\\
(2) In the correspondence in (1), the rank of the indecomposable special CM $R$-module corresponding to a curve $E_i$ is equal to the coefficient of $E_i$ in the fundamental cycle $Z_f$.\\
(3) For the dihedral groups $\mathbb{D}_{n,q}$, there are precisely $(N+2-\nu)$ non-free rank one special CM modules, and $\nu$ rank two special indecomposable CM modules.
\end{thm}
\begin{proof}
(1) and (2). This is \cite[1.2]{Wunram_generalpaper}.\\
(3)  By the above discussion on $Z_f$ for dihedral groups, this follows directly from (2).
\end{proof}
Note that for dihedral groups $\mathbb{D}_{n,q}$, the rank two indecomposable special CM modules are known from \cite[6.2]{Iyama_Wemyss_specials}; in fact the classification of the special CM modules for all finite subgroups of $\GL(2,\C{})$ is complete with the exception of these $(N+2-\nu)$ non-free rank one special CM modules in the dihedral cases.

To be more precise, Wunram defined the specials using CM modules on the ring $\C{}[[x,y]]^G$ where such modules are of the form $(\rho\otimes_{\C{}}\C{}[[x,y]])^G$.  In this paper we shall mainly be working with the $\C{}[x,y]^G$-modules $(\rho\otimes_{\C{}}\C{}[x,y])^G$, i.e. we work in the non-complete case.  We are mainly interested in computing the endomorphism ring in the non-complete case, but later we shall reduce this problem to the complete case since the associated graded ring of $\t{End}_{\C{}[[x,y]]^G}(\bigoplus_{\rho\t{ special}}(\rho\otimes_{\C{}}\C{}[[x,y]])^G)$ is $\t{End}_{\C{}[x,y]^G}(\bigoplus_{\rho\t{ special}}(\rho\otimes_{\C{}}\C{}[x,y])^G)$.

To find these specials, and thus finish the classification, we need the combinatorial $i$-series (as in type $A$) together with some other combinatorial series:

\begin{defin}\label{combdata}
For any integers $1\leq m_1<m_2$ with $(m_1,m_2)=1$ we can associate to
the continued fraction expansion
$\frac{m_2}{m_1}=[\beta_1,\hdots,\beta_X]$ combinatorial series defined as follows:\\
1. The $i$-series, defined as
\[
\begin{array}{ccl}
i_0=m_2 & i_1=m_1 & i_{t}=\beta_{t-1}i_{t-1}-i_{t-2}\mbox{ for }2\leq
t\leq X+1.
\end{array}
\]
2. The $j$-series, defined as
\[
\begin{array}{ccl}
j_0=0 & j_1=1 & j_{t}=\beta_{t-1}j_{t-1}-j_{t-2}\mbox{ for }2\leq
t\leq X+1.
\end{array}
\]
3. The $l$-series, defined as
\[
l_j=2+\sum_{p=1}^{j}(\beta_p-2) \mbox{ for }1\leq j\leq X.
\]
4.  The $b$-series. Define $b_{0}:=1$, $b_{l_X-1}:=X$, and further for all  $1\leq t\leq l_X-2$ (if such $t$ exists), define $b_t$  to be the smallest integer $1\leq
b_t\leq X$ such that %\footnote{records the butt of $k_t$ for the reconstruction algebra of type $A$ associated to $[\beta_1,\hdots,\beta_X]$, where we label the $k$'s in the opposite direction to that of type $A$!}\footnote{Need $b_0=1$ (unlike $\tt{B}_0=+$ later in the relations) so that the $t=0$ case of Lemma~\ref{Relating_the_two_i_series} is ok.}
\[
t\leq \sum_{p=1}^{b_t}(\beta_p-2).
\]
\end{defin}
\begin{defin}
Given a continued fraction expansion $\frac{m_2}{m_1}=[\beta_1,\hdots,\beta_X]$ we call the continued fraction expansion of $\frac{m_2}{m_2-m_1}$ the dual continued fraction, and denote it by $[\beta_1,\hdots,\beta_X]^\vee$.
\end{defin}

Throughout this section we shall be using the above combinatorial series for many different continued fraction inputs, thus to avoid confusion we now fix some notation.
\begin{notation}
For $\mathbb{D}_{n,q}$, throughout this paper we fix $\frac{n}{n-q}=[a_2,\hdots,a_{e-1}]$ and $\frac{n}{q}=[\alpha_1,\hdots,\alpha_N]$.  We shall denote the combinatorial data in Definition~\ref{combdata} associated to the continued fraction expansion of $\frac{n}{q}$ exclusively by using the fonts and letters $(i,j,l,b)$.  For all other continued fraction inputs we shall denote the combinatorial data in Definition~\ref{combdata} by using different fonts and letters.
\end{notation}

We record some easy combinatorics.  %Throughout this section we shall often refer to Riemenschneider duality; for a reference for this fact see \cite[1.2]{Kidoh} or \cite[2.11]{Wemyss_reconstruct_A}.
\begin{lemma}\label{a2,q,r2,r3,r2}
For any $\mathbb{D}_{n,q}$ with any $\nu$,\\
$
\begin{array}{cl}
\t{(i)} & a_2=\nu+2.\\
\t{(ii)}& q=i_{\nu+1}+\nu(n-q).\\
\t{(iii)}& r_2=2(n-q)-i_{\nu+1}.\\
\t{(iv)} & r_3=(n-q)-i_{\nu+1}=i_{\nu}-2i_{\nu+1}.\\
\t{(v)} &r_2=2r_3+i_{\nu+1}.
\end{array}
$
\end{lemma}
\begin{proof}
(i) This is immediate by Riemenschneider duality. \\
%\[
%\frac{n}{n-q}=[\underbrace{2,\hdots,2}_{u_{\sigma_1}-v_{\sigma_1}},(\sigma_2-\sigma_1+2),\hdots].
%\]
%If $\nu=0$ then $u_{\sigma_1}>v_{\sigma_1}$ and so $a_2=2=\nu+2$.
%Hence assume $\nu>0$ in which case $u_{\sigma_1}=v_{\sigma_1}$
%with $1=\sigma_1<\sigma_2=\nu+1$ and so
%$a_2=\sigma_2-\sigma_1+2=\nu+2$.\\
(ii) Trivially this is true if $\nu=0$, so we can assume that
$\nu>0$. This being the case it is easy to see that
\[
i_t=tq-(t-1)n \mbox{ for all } 1\leq t\leq \nu+1 \tag{$\star$}
\]
since $\alpha_1=\hdots=\alpha_\nu=2$.  In particular
$i_{\nu+1}=(\nu+1)q-\nu n$ and so the result is trivial.\\
(iii) This follows from (i) and (ii) since
\[
r_2=a_2(n-q)-q(c_2+d_2)=(\nu+2)(n-q)-(i_{\nu+1}+\nu(n-q))(0+1)=2(n-q)-i_{\nu+1}.
\]
(iv) Immediate from (iii) since
\[
r_3=(a_2-1)(n-q)-q(c_3+d_3)=(a_2-1)(n-q)-q=r_2-(n-q).
\]
(v) Immediate from (iii) and (iv) above.
\end{proof}

\begin{lemma}\label{Relating_the_two_i_series}
Take some continued fraction expansion $\frac{m_2}{m_1}=[\beta_1,\hdots,\beta_X]$ and denote the combinatorial data from Definition~\ref{combdata} by $(I,J,L,B)$.  To the dual continued fraction $[\beta_1,\hdots,\beta_X]^\vee:=[\gamma_1,\hdots,\gamma_Y]$ denote the combinatorial data by $(\mathbbold{I},\mathbbold{J},\mathbbold{L},\mathbbold{B})$.  Then \\
$
\begin{array}{cl}
\t{(i)}& B_t=\mathbbold{L}_t-1 \mbox{ for all } 1\leq t\leq L_X-1.\\%\footnote{limit is 1 since $\mathbbold{L}_0$ is not defined}
\t{(ii)}& \mathbbold{I}_t=\mathbbold{I}_{t+1}+I_{B_t} \mbox{ for all } 0\leq t\leq {L_X-1}.\\
\t{(iii)}& J_{t+1}-J_{t}=\mathbbold{J}_{\mathbbold{B}_{t}} \mbox{ for all } 0\leq t\leq \mathbbold{L}_Y-1. \\
\t{(iv)}& \mathbbold{J}_{t+1}-\mathbbold{J}_{t}=J_{\mathbbold{L}_t-1} \mbox{ for all } 1\leq t\leq L_X-1.\\%\footnote{limit is 1 since use the limits in (i) to prove (iv)}
\t{(v)}& \mathbbold{J}_{\mathbbold{B}_{t}}=1+\sum_{p=1}^{\mathbbold{B}_{t}-1}J_{\mathbbold{L}_p-1} \mbox{ for all } 1\leq t\leq \mathbbold{L}_Y-1.  
\end{array}
$
\end{lemma}
\begin{proof}
(i) is an immediate consequence of Riemenschneider duality.\\
(ii) and (iii) are just a slight rephrasing of a result of Kidoh \cite[1.3]{Kidoh}.\\ %\footnote{for (i), (ii) and (iii) see my Combinatorics notes. Note that in the dictionary between Kidoh and me, her results don't give the outer limits in (ii) and (iii) - they are are true only by convention (i.e. by the definition of $B_0$ and $B_{L_X-1}$.)} 
(iv) By duality, swapping bold and non-bold in (iii) gives $\mathbbold{J}_{t+1}-\mathbbold{J}_{t}=J_{B_t}$ for all $0\leq t\leq L_X-1$.  The result then follows by (i).\\
(v) Follows immediately from (iv) since
\[
\mathbbold{J}_{\mathbbold{B}_t}-1=\mathbbold{J}_{\mathbbold{B}_t}-\mathbbold{J}_{1}=(\mathbbold{J}_{\mathbbold{B}_t}-\mathbbold{J}_{\mathbbold{B}_t-1})+\hdots+(\mathbbold{J}_{2}-\mathbbold{J}_{1})=\sum_{p=1}^{\mathbbold{B}_t-1}J_{\mathbbold{L}_p-1}.
\]
\end{proof}

The following is known, and can be found in \cite[p.214]{Riemenschneider_dihedral}.
\begin{lemma}\label{r_series_is_I_series}
For $\mathbb{D}_{n,q}$ as above with
$\frac{n}{n-q}=[a_2,a_3,\hdots,a_{e-1}]$, then the $r$-series is
simply the $i$-series for the data $[a_3+1,a_4,\hdots,a_{e-1}]$. More
precisely, denoting the $i$-series of $[a_3+1,a_4,\hdots,a_{e-1}]$ by
$\mathcal{I}_0,\mathcal{I}_1,\hdots$, we have $r_k=\mathcal{I}_{k-2}$ for all $2\leq k\leq e$.
\end{lemma}
\begin{proof}
By definition $\frac{n}{n-q}=[a_2,\hdots,a_{e-1}]=a_2-\frac{1}{[a_3,\hdots,a_{e-1}]}$ and so $\frac{n-q}{a_2(n-q)-n}=[a_3,\hdots,a_{e-1}]$.
But combining Lemma~\ref{a2,q,r2,r3,r2}(iii),(iv) we know that $r_3=r_2-(n-q)$ and further since $r_2=a_2(n-q)-q$ (by definition) we have
\[
\frac{r_2}{r_3}=\frac{a_2(n-q)-q}{a_2(n-q)-q-(n-q)}=\frac{n-q+a_2(n-q)-n}{a_2(n-q)-n}=[a_3+1,a_4,\hdots,a_{e-1}].
\] 
The result now follows since $r_4=(a_3+1)r_3-r_2$ %\footnote{easy manipulation: see Allinone.pdf p3}
and $r_t=a_{t-1}r_{t-1}-r_{t-2}$ for all $5\leq t\leq e$.
\end{proof}

Note in particular this means that $r_{e-1}=1$ and $r_e=0$.

\begin{lemma}\label{r_butt_lemma}
Consider $\mathbb{D}_{n,q}$.  Then for all $2\leq t\leq e-2$, %\footnote{shows that for $\nu=N-1$ that $r_t-\alpha_N-(t-1)$; see allinone.pdf p36 and 37.}\footnote{shows the ending of the $r$'s on vertex $N$ to be decreasing by 1 each time.}
\[
r_{t+1}=r_{t+2}+i_{b_t}.
\]
\end{lemma}
\begin{proof}
Denote the $i$-series associated to the following data as follows:
\[
\begin{array}{rcl}
\frac{n}{n-q}=[a_2,a_3,\hdots,a_{e-1}] &\t{by} & \iota_0,\iota_1,\hdots\\
{}[a_3,a_4,\hdots,a_{e-1}] &\t{by} & \tilde{\iota}_0,\tilde{\iota}_1,\hdots\\
{}[a_3+1,a_4,\hdots,a_{e-1}] & \t{by} & \mathcal{I}_0,\mathcal{I}_1,\hdots \t{(as in Lemma~\ref{r_series_is_I_series})}
\end{array}
\]
Since $\frac{\tilde{\iota}_0+\tilde{\iota}_1}{\tilde{\iota}_1}=[a_3+1,a_4,\hdots,a_{e-1}] $ it is clear that for $0\leq j\leq e$
\[
\mathcal{I}_j=\left\{  \begin{array}{cc}\tilde{\iota}_0+\tilde{\iota}_1 & \mbox{if
}j=0\\\tilde{\iota}_j &\mbox{if }j\geq 1
\end{array}\right. =\left\{  \begin{array}{cc}\iota_1+\iota_2 &
\mbox{if }j=0\\ \iota_{j+1} &\mbox{if }j\geq 1
\end{array}\right. .
\]
Now by Lemma~\ref{r_series_is_I_series} for all $2\leq k\leq e$ we
have $r_{k}=\mathcal{I}_{k-2}$. Thus %\footnote{Note to self: the $t=1$ case above (which we ignored) gives a new proof of Lemma~\ref{a2,q,r2,r3,r2}(v) since $r_2=I_0=\iota_1+\iota_2=2\iota_2+i_{b_1}=2I_1+i_{b_1}=2r_3+i_{\nu+1}.$\]} 
for $2\leq t\leq e-2$, by the above and Lemma~\ref{Relating_the_two_i_series} applied to $[\beta_1,\hdots,\beta_X]=[\alpha_1,\hdots,\alpha_N]$ we have
\[
r_{t+1}=\mathcal{I}_{t-1}=\iota_{t}=\iota_{t+1}+i_{b_{t}}=\mathcal{I}_{t}+i_{b_t}=r_{t+2}+i_{b_t}.
\]
\end{proof}

\begin{lemma}\label{r_is_difference_in_i_series}
Consider $\mathbb{D}_{n,q}$, then for all $\nu+1\leq t\leq N$
\[
r_{l_t}=i_t-i_{t+1}.
\]
\end{lemma}
\begin{proof}
Proceed by induction.  Consider first the base case $t=\nu+1$.  Notice it is always true (by definition) that $l_{\nu+1}=\alpha_{\nu+1}$.  Now to prove the base case requires two subcases:\\
(i) If $\alpha_{\nu+1}=3$ then by Lemma~\ref{a2,q,r2,r3,r2}(iv)
\[
r_{l_{\nu+1}}=r_{\alpha_{\nu+1}}=r_3=i_{\nu}-2i_{\nu+1}=i_{\nu+1}-i_{\nu+2}
\]
(ii) If $\alpha_{\nu+1}>3$ then Lemma~\ref{a2,q,r2,r3,r2}(iv) and
Lemma~\ref{r_butt_lemma} we have
\[
r_{l_{\nu+1}}=r_3-(l_{\nu+1}-3)i_{\nu+1}=(i_{\nu}-2i_{\nu+1})-(\alpha_{\nu+1}-3)i_{\nu+1}=i_{\nu+1}-i_{\nu+2}
\]
and so we are done.  This proves the base case $t=\nu+1$.  If $\nu=N-1$ we are done hence suppose $\nu<N-1$, let $t$ be such that $\nu+1<t\leq N$ and assume that the result is true for smaller $t$.  To prove the induction step again requires two subcases:\\
(i) If $\alpha_t=2$, then $l_t=l_{t-1}$ and so by inductive
hypothesis
\[
r_{l_t}=r_{l_{t-1}}=i_{t-1}-i_{t}=i_t-i_{t-1}.
\]
(ii) If $\alpha_t>2$ then by Lemma~\ref{r_butt_lemma} and
inductive hypothesis
\[
r_{l_t}=r_{l_{t-1}}-(l_t-l_{t-1})i_t=(i_{t-1}-i_{t})-(\alpha_t-2)i_t=i_t-i_{t+1}.
\]
\end{proof}

\begin{defin}%\footnote{Notice the $k-1$ on top of the sum!} 
Consider $\mathbb{D}_{n,q}$.  Define for $\nu+1\leq k\leq N+1$
\[
\Delta_k=1+\sum_{t=\nu+1}^{k-1} c_{l_t} \quad\mbox{and}\quad
\Gamma_k=\sum_{t=\nu+1}^{k-1} d_{l_t}
\]
where the convention is that for $k=\nu+1$ the sum is empty and so
equals zero.
\end{defin}

The following is the key technical result of this section.  Its main use can be seen visually in Example~\ref{444pattern}, and it forms the basis of the combinatorial arguments used in Proposition~\ref{loops_in_D1} (and hence the remainder of the paper).

\begin{lemma}\label{c_and_d_lemma} Consider $\mathbb{D}_{n,q}$ for any $\nu$.  Then for all $2\leq t\leq e-2$
\[
c_{t+2}=c_{t+1}+\Delta_{b_t}\quad\mbox{and}\quad d_{t+2}=d_{t+1}+\Gamma_{b_t}.
\]
%\t{(ii)} If $\nu=N-1$ then for all $3\leq t\leq e$
%\[
%c_{t}=t-3,\qquad d_{t}= 1
%\]
\end{lemma}
\begin{proof}
(i) The $c$ statement.  The trick is to interpret the $c$'s as the $j$-series associated to some continued fraction, then use the results of Lemma~\ref{Relating_the_two_i_series}.

We first prove that the lemma holds in the case $t=2$.  Notice that $\Delta_{b_2}=1$ since either $b_2=\nu+1$ and so the sum is empty (and so by convention zero), or $b_2>\nu+1$ in which case $l_{\nu+1}=\hdots=l_{b_2-1}=3$ and so the sum is $c_3+\hdots+c_3=0$.  Thus $c_4=c_3+\Delta_{b_2}$ follows since $c_4:=1$ and $c_3:=0$.  Hence the result is true for $t=2$ thus we may %\footnote{this is necessary to write the Riemenschneider duality out nicely.} 
and restrict our attention to the interval $3\leq t\leq e-2$.

Denote the $j$-series of $[a_4,\hdots,a_{e-1}]$ by $\tt{j}$.    It is clear from the definition of $c$ that
\[
c_t=\tt{j}_{t-3}\quad\mbox{ for all } 3\leq t\leq e.
\]
To $[a_4,\hdots,a_{e-1}]^\vee$ denote the $j$-series by $\mathbbold{j}$, the $b$-series by $\mathbbold{b}$ and the $l$-series by $\mathbbold{L}$.  By Lemma~\ref{Relating_the_two_i_series}(iii) applied to the data $[\beta_1,\hdots,\beta_X]=[a_4,\hdots,a_{e-1}]$ %\footnote{limits ok since $\mathbbold{L}_Y-1=e-4$, see footnote below}
\[
c_{t+2}-c_{t+1}=\tt{j}_{t-1}-\tt{j}_{t-2}=\mathbbold{j}_{\mathbbold{b}_{t-2}}\quad\mbox{ for all } 2\leq t\leq e-2.
\]
On the other hand
\[
\Delta_{b_t}=1+\sum_{p=\nu +1}^{b_t-1}c_{l_p}=1+\sum_{p=\nu+1}^{b_t-1}\tt{j}_{l_p-3},
\]
thus to prove the lemma we just need to show that %\footnote{ok as have already proved the $t=2$ case; note 3 important in the limits of $\mathbbold{b}$ below} 
\[
\mathbbold{j}_{\mathbbold{b}_{t-2}}=1+\sum_{p=\nu+1}^{b_t-1}\tt{j}_{l_p-3}\quad\mbox{ for all } 3\leq t\leq e-2.
\]
Now by Riemenschneider duality %\footnote{In the top necessarily $\nu<N-1$ by $e>4$ assumption, in the bottom $\nu=N-1$ is allowed.} \footnote{convention is that the left term always exists, whereas the ones after that may not}
\[
[a_4,\hdots,a_{e-1}]^\vee=\left\{  \begin{array}{cl} [\alpha_{b_2}-1,\alpha_{1+b_2},\hdots,\alpha_N]&\mbox{if }\alpha_{\nu+1}=3\\  
{[ \alpha_{b_2}-2,\alpha_{1+b_2},\hdots,\alpha_N]}  &\mbox{if }\alpha_{\nu+1}>3\end{array} 
\right.
\]
from which it is easy to see %\footnote{p4,5 proof of c} 
in either case that %\footnote{Note in both cases that $\mathbbold{L}_Y-1=1+\sum (\gamma-2)=e-4$ so this is the top $\mathbbold{b}$}
\[
\mathbbold{b}_0=1 \mbox{ whereas } \mathbbold{b}_{t}=b_{t+2}-(b_2-1) \mbox{ for all } 1\leq t\leq e-4
\]
and $\mathbbold{L}_s=l_{(b_2-1)+s}-2$ for all %\footnote{the extreme is the number of terms in the dual continued expansion, which is by definition the top for the $\mathbbold{L}$ series.} 
$1\leq s\leq N-(b_2-1)$.  Hence %\footnote{first equality certainly ok if $b_2=\nu+1$; else $l_{\nu+1}=\hdots=l_{b_2-1}=3$ and so ok since $\tt{j}_{3-3}=\tt{j}_0=0$}
\[
1+\sum_{p=\nu+1}^{b_t-1}\tt{j}_{l_p-3}=1+\sum_{p=b_2}^{b_t-1}\tt{j}_{l_p-3}=1+\tt{j}_{\mathbbold{L}_1-1}+\hdots+\tt{j}_{\mathbbold{L}_{(\mathbbold{b}_{t-2}-1)}-1}=1+\sum_{p=1}^{\mathbbold{b}_{t-2}-1}\tt{j}_{\mathbbold{L}_p-1}=\mathbbold{j}_{\mathbbold{b}_{t-2}}
\]
for all $3\leq t\leq e-2$ where the last equality holds by Lemma~\ref{Relating_the_two_i_series}(v) %\footnote{got the limits exactly right since as stated in previous footnote $\mathbbold{L}_Y-1=e-4$, so use limits in Lemma~\ref{Relating_the_two_i_series}(v) and add 2 to both sides due to the change in $t$}
applied to the data $[\beta_1,\hdots,\beta_X]=[a_4,\hdots,a_{e-1}]$.\\

(ii) The $d$ statement.  The trick is again to interpret the $d$-series as the $j$-series of something, but here it is a little bit more subtle.

Due to lack of suitable alternatives we recycle notation from the proof of (i): now denote the $j$-series of %\footnote{convention is that the left hand side always exists, whereas the terms to the left may not} 
$[a_{\alpha_{\nu+1}}-1,a_{\alpha_{\nu+1}+1},\hdots, a_{e-1}]$ by $\tt{j}$, %\footnote{Note that $\tt{j}$ here is defined with respect to a continued fraction involving $a$'s, just like before} 
and further for the dual continued fraction $[a_{\alpha_{\nu+1}}-1,a_{\alpha_{\nu+1}+1},\hdots, a_{e-1}]^\vee$ denote the $j$, $b$ and $l$ series by $\mathbbold{j}$, $\mathbbold{b}$ and $\mathbbold{L}$ respectively.

Now it is easy to see that $d_3=\hdots=d_{\alpha_{\nu+1}}=1$  and further $d_t=\tt{j}_{(t+1)-\alpha_{\nu+1}}$ for all $\alpha_{\nu+1}+1\leq t\leq e$. %\footnote{actually $l_{\nu+1}\leq t$ ok, but $d_{l_{\nu+1}}=1$ already covered.}   
Hence the result is certainly true (by the convention of the empty sum being zero) for the interval $2\leq t\leq \alpha_{\nu+1}-1$. Thus we are done in the case $\nu=N-1$ %\footnote{since $e-2=1+\sum_{p=1}^{N}(\alpha_p-2)=\alpha_{\nu+1}-1$ in this case} 
and also in the case $e=4$.  Thus we may assume that $\nu<N-1$ and that $e>4$, %\footnote{again I think this is need to write Riemenschneider duality nicely} 
and concentrate on the interval $(\alpha_{\nu+1}-2)+1\leq t\leq e-2=(\alpha_{\nu+1}-2)+(1+\sum_{p=\nu+2}^{N}(\alpha_p-2))$.

By Lemma~\ref{Relating_the_two_i_series}(iii) applied to the data $[\beta_1,\hdots,\beta_X]=[a_{\alpha_{\nu+1}}-1,a_{\alpha_{\nu+1}+1},\hdots, a_{e-1}]$
\[
d_{t+2}-d_{t+1}=\tt{j}_{(t+3)-\alpha_{\nu+1}}-\tt{j}_{(t+2)-\alpha_{\nu+1}}=\mathbbold{j}_{\mathbbold{b}_{(t+2)-\alpha_{\nu+1}}}
\]
for all $(\alpha_{\nu+1}-2)+1\leq t\leq (\alpha_{\nu+1}-2)+(1+\sum_{p=\nu+2}^{N}(\alpha_p-2))$.

On the other hand
\[
\Gamma_{b_t}=\sum_{p=\nu+1}^{b_t-1}d_{l_p}=\sum_{p=\nu+1}^{b_t-1}\tt{j}_{(l_p-\alpha_{\nu+1})+1}
\]
and so the result follows if we can show that
\[
\mathbbold{j}_{\mathbbold{b}_{(t+2)-\alpha_{\nu+1}}}=\sum_{p=\nu+1}^{b_t-1}\tt{j}_{(l_p-\alpha_{\nu+1})+1}
\]
for all $(\alpha_{\nu+1}-2)+1\leq t\leq (\alpha_{\nu+1}-2)+(1+\sum_{p=\nu+2}^{N}(\alpha_p-2))$, i.e.%\footnote{$t\mapsto t+(\alpha_{\nu+1}-2)$}
\begin{eqnarray}
\mathbbold{j}_{\mathbbold{b}_t}=\sum_{p=\nu+1}^{b_{(\alpha_{\nu+1}-2)+t}-1}\tt{j}_{(l_p-\alpha_{\nu+1})+1}\label{1}
\end{eqnarray}
for all  $1\leq t\leq 1+\sum_{p=\nu+2}^{N}(\alpha_p-2)$. %\footnote{just subtract $(\alpha_{\nu+1}-2)$ from the limits of the interval above!}  
Now by Riemenschneider duality %\footnote{note only exists if $\nu<N-1$; we've already reduced to this case.} 
\[
[a_{\alpha_{\nu+1}}-1,a_{\alpha_{\nu+1}+1},\hdots, a_{e-1}]^\vee=[\alpha_{\nu+2},\hdots,\alpha_N]
\]
and so $\mathbbold{L}_Y-1=1+\sum_{p=\nu+2}^{N}(\alpha_p-2)$.  Further it is easy to see that $\mathbbold{b}_0=1$, $\mathbbold{b}_t=b_{(\alpha_{\nu+1}-2)+t}-(\nu+1)$ for all $1\leq t\leq 1+\sum_{p=\nu+2}^{N}(\alpha_p-2)$ %\footnote{the top of $\mathbbold{b}$ series is just $\mathbbold{L}_Y-1=1+\sum_{p=\nu+2}^{N}(\alpha_p-2)=(l_N-1)-(\alpha_{\nu+1}-2)$ so substitute in then top $b$ we get is $b_{1+\sum_{p=\nu+1}^{N}(\alpha_p-2)}$ which is the top $b$} 
and $\mathbbold{L}_t=l_{(\nu+1)+t}-\alpha_{\nu+1}+2$ for all $1\leq t\leq N-\nu$. %\footnote{top extreme is just the number of terms in $[\alpha_{\nu+1}-1,\alpha_{\nu+2},\hdots, \alpha_N]$.} \footnote{recall $\alpha_{\nu+1}=l_{\nu+1}$ so really $\mathbbold{L}_t=l_{(\nu+1)+t}-l_{\nu+1}+2$} 
This implies that the sum in (\ref{1}) is simply%\footnote{$\mathbbold{b}_t+\nu=b_{(\alpha_{\nu+1}-2)+t}-1$}
\[
\sum_{p=\nu+1}^{\mathbbold{b}_t+\nu}\tt{j}_{(l_p-\alpha_{\nu+1})+1}=\tt{j}_1+\tt{j}_{(l_{\nu+2}-\alpha_{\nu+1})+1}+\hdots+\tt{j}_{(l_{\mathbbold{b}_t+\nu}-\alpha_{\nu+1})+1}=1+\tt{j}_{\mathbbold{L}_1-1}+\hdots+\tt{j}_{\mathbbold{L}_{\mathbbold{b}_t-1}-1}
\]
and so by Lemma~\ref{Relating_the_two_i_series}(v) applied to the data $[\beta_1,\hdots,\beta_X]=[a_{\alpha_{\nu+1}}-1,a_{\alpha_{\nu+1}+1},\hdots, a_{e-1}]$
\[
\sum_{p=\nu+1}^{b_{(\alpha_{\nu+1}-2)+t}-1}\tt{j}_{(l_p-\alpha_{\nu+1})+1}=1+\sum_{p=1}^{\mathbbold{b}_t-1}\tt{j}_{\mathbbold{L}_p-1}=\mathbbold{j}_{\mathbbold{b}_t}
\]
for all $1\leq t\leq \mathbbold{L}_Y-1=1+\sum_{p=\nu+2}^{N}(\alpha_p-2)$, as required.
\end{proof}
%Easy proof of \nu=N-1 case:
%
%Since $\nu=N-1$ we may write
%\[
%\frac{n}{q}=[\underbrace{2,\hdots,2}_{N-1},\alpha_N]
%\]
%By Riemenschneider duality
%\[
%\frac{n}{n-q}=[N+1,\underbrace{2,\hdots,2}_{\alpha_N-2}]
%\]
%and so $a_2=N+1$ whereas $a_3=\hdots=a_{e-1}=2$.  By definition of $c_t$ and $d_t$ the result is %now trivial.
%\end{proof}

\section{Computation of the rank one specials}

In this section we determine the rank one special CM modules for all $\mathbb{D}_{n,q}$, and obtain their generators.  Denoting $R:=\C{}[x,y]^{\mathbb{D}_{n,q}}$, we do this by first viewing CM $R$-modules as $R$-summands of the polynomial ring \cite{Herzog}.  We then factorize the generators of the invariant ring $R$ via a certain quiver, and argue by diagram chasing that the modules corresponding to the vertices are two-generated.  Thus in this section we view a certain quiver $\tt{D}_1$ (see Notation~\ref{D}) as a convenient method of factorizing certain polynomials; in later sections the quiver $\tt{D}_1$ forms part of the reconstruction algebra. 

We begin with the following simple lemma.
\begin{lemma}\label{maps1d}
Let $G$ be a small finite subgroup of $\GL(2,\C{})$.  For one-dimensional %\footnote{I don't think this is necessary, but to get the $G$-equiv of the last iso it is easiest to use 1d spaces and characters thereof; note that $G$ action on the Hom space is given on p85 Yoshino} 
representations $\rho,\sigma$, denote the corresponding rank one CM modules by $S_\rho$ and $S_\sigma$.  Then as $\C{}[x,y]^G$-modules
\[
\t{Hom}_{\C{}[x,y]^G}(S_{\rho},S_{\sigma})\cong
S_{\sigma\otimes\rho^\ast}.
\]
\end{lemma}
\begin{proof}
Since the group $G$ is small, it is well known that $\C{}[x,y]\# G\cong \t{End}_{\C{}[x,y]^G}(\C{}[x,y])$ and $\t{add}\C{}[x,y]=\CM\C{}[x,y]^G$.  Consequently 
\[
\begin{array}{rcl}
\proj \C{}[x,y]\# G & \approx & \CM \C{}[x,y]^G\\
M&\mapsto& M^G
\end{array}
\]
is an equivalence of categories (for details see for example \cite[10.9]{Yoshino1}) and thus 
\begin{eqnarray*}
\t{Hom}_{\C{}[x,y]^G}(S_{\rho},S_{\sigma})&\cong&\t{Hom}_{\C{}[x,y]\#G}(\C{}[x,y]\otimes_{\C{}}{\rho},\C{}[x,y]\otimes_{\C{}}{\sigma})\\
&\cong&\t{Hom}_{\C{}[x,y]}(\C{}[x,y]\otimes_{\C{}}{\rho},\C{}[x,y]\otimes_{\C{}}{\sigma})^G\\
&\cong&(\C{}[x,y]\otimes_{\C{}}(\sigma\otimes\rho^\ast))^G=S_{\sigma\otimes\rho^\ast}
\end{eqnarray*}
\end{proof}

\begin{defin}
For $1\leq t\leq n-q$ define $W_t$ to be the CM $\C{}[x,y]^{\mathbb{D}_{n,q}}$-module containing
$(xy)^t$.
\end{defin}
We should make two remarks.  First, the $W_t$ are well defined since $(xy)^t$ is a relative invariant for the one-dimensional representation%\footnote{Here, and for $W_+$,$W_-$, the convention is $(\C{}[x,y]\otimes \rho)^G$ where $G$ acts on the polynomial ring as inverses and on $\rho$ normally}
\[
\begin{array}{c|c}
n-q \mbox{ odd} & n-q \mbox{ even}\\ \hline \begin{array}{ccl}
\psi_{2q}&\mapsto& 1\\\tau & \mapsto &
(-1)^t\\\varphi_{2(n-q)}&\mapsto &\e_{n-q}^t
\end{array} & \begin{array}{ccl}
\psi_{2q}&\mapsto& 1 \\\tau\varphi_{4(n-q)}&\mapsto &
\e_{2(n-q)}^{(n-q)+t}
\end{array}
\end{array}
\]
Second, the assumption $t<n-q$ ensures that the $W_{t}$ are mutually non--isomorphic representations.  Now for any $\mathbb{D}_{n,q}$, by Lemma~\ref{a2,q,r2,r3,r2}(iv) it is true that $i_{\nu+1}<n-q$ and so we aim to prove that $W_{i_{\nu+1}},W_{i_{\nu+2}},\hdots,W_{i_{N}}$ are special CM modules.  Since they have rank one, we just need to show that they are two-generated:

\begin{lemma}\cite{Wunram_generalpaper}\label{2genlemma}
Suppose that $M$ is a rank one CM $\C{}[x,y]^G$-module which is minimally 2-generated.  Then $M$ is special.  
\end{lemma}
\begin{proof}
Denote the minimal resolution of $\C{2}/G$ by $\pi:\w{X}\rightarrow \C{2}/G$, and the irreducible exceptional curves by $\{ E_i\}_{i\in I}$.   For a sheaf $\s{F}$ on $\w{X}$, we denote $\s{F}^{\vee}$ to be the sheaf $\s{H}om_{\w{X}}(\s{F},\s{O}_{\w{X}})$ and we denote $\mathbf{T}(\s{F})$ to be the torsion subsheaf of $\s{F}$, i.e. the kernel of the natural map $\s{F}\to\s{F}^{\vee\vee}$. To ease notation, if $M$ is a CM $\C{}[x,y]^G$-module, we denote $\c{M}:=\pi^*M/\mathbf{T}(\pi^\ast M)$.

Now by \cite[2.1]{Wunram_generalpaper} (see also \cite[3.5]{IWtri}) $M$ is generated by $\rk M + Z_f \cdot c_1(\c{M})$ elements.   This then implies that $Z_f \cdot c_1(\c{M})=1$.  But the fundamental cycle $Z_f=\sum_{i\in I}a_iE_i$ where each $a_i\geq 1$, and also $c_1(\c{M})\cdot E_i\geq 0$ for all $i$.  Hence
\[
1=\sum_{i\in I}a_i (c_1(\c{M})\cdot E_i)
\]
forces $c_1(\c{M})\cdot E_i=\delta_{ij}$ for some $j\in I$ (where $\delta_{ij}$ is the Kronecker delta), and further $a_j=1$.  

But since $a_j=1$, by Wunram \cite[1.2]{Wunram_generalpaper} there exists a special CM module $M_j$, of rank one, such that $c_1(\c{M}_j)\cdot E_i=\delta_{ij}$.  Since line bundles on $\w{X}$ are uniquely determined by their first Chern classes (see e.g. \cite[3.4.3]{VdB}) we conclude that $\c{M}\cong\c{M}_j$, from which taking global sections yields $M\cong M_j$.
\end{proof}

The next lemma is trivial but is used extensively.
\begin{lemma}\label{polys_to_invariants}
Consider a rank one CM $\C{}[x,y]^G$-module $T$.  Let $f_1,f_2\in T$ be such that every element of $T$ may be written as $Af_1+Bf_2$ for some polynomials $A$ and $B$.  Then $T$ is generated as an $\C{}[x,y]^G$-module by $f_1$ and $f_2$.
\end{lemma}
\begin{proof}
Let $a\in T$ be written as $a=Af_1+Bf_2$.  We just need to show that we can replace $A$ and $B$ with polynomials inside $\C{}[x,y]^G$, then the result obviously follows.  Taking any element $g$ of the group $G$, since $f_1,f_2\in T$ we may act on $a$ and then cancel the relative invariant scalars, leaving
$
a=(g\cdot A)f_1+(g\cdot B)f_2
$.  
Thus summing over all group elements,
$
a=\frac{1}{|G|}(\sum_{g\in G} g\cdot A) f_1+\frac{1}{|G|}(\sum_{g\in G}g\cdot B) f_2.
$
\end{proof}
For technical reasons we split the proof that $W_{i_{\nu+1}},\hdots, W_{i_N}$ are special into two cases:

\medskip
\noindent
\textbf{Case 1: $0\leq \nu<N-1$}, so $\alpha_{\nu+1}\geq 3$.

\medskip\noindent  By definition and Lemma~\ref{maps1d}, certainly we have the following maps:
\[
\xymatrix@C=65pt{ W_{i_{\nu+1}}&W_{i_{\nu+2}}\ar@/_0.65pc/[l]|{(xy)^{i_{\nu+1}-i_{\nu+2}}}\ar@{}[r]|{\hdots} & W_{i_{N-1}} & W_{i_N}\ar@/_0.65pc/[l]|{(xy)^{i_{N-1}-i_N}}&R\ar@/_0.65pc/[l]|{xy} }
\]
Now by Lemma~\ref{a2,q,r2,r3,r2}(iii) $2(n-q)=r_2+i_{\nu+1}$, so
since $(xy)^{2(n-q)}$ is an invariant we also have a map
$(xy)^{r_2}:W_{i_{\nu+1}}\rightarrow R$.  Since
$(xy)^{r_2}w_2$ is an invariant this in turn means\\
(i) there is a map $w_2:R\rightarrow W_{i_{\nu+1}}$.\\
(ii) (since $2r_3+i_{\nu+1}=r_2$) there is a map $g:=(xy)^{2r_3}w_2:W_{i_{\nu+1}}\rightarrow R$.\\
But also $(xy)^{r_{l_t}}w_2^{c_{l_t}}w_3^{d_{l_t}}$ is invariant,
so by Lemma~\ref{r_is_difference_in_i_series} we have, for each
$\nu+1\leq t\leq N-1$, a map
\[
w_2^{c_{l_t}}w_3^{d_{l_t}}:W_{i_t}\rightarrow W_{i_{t+1}}
\]
and also a map $w_2^{c_{l_N}}w_3^{d_{l_N}}:W_{i_N}\rightarrow R$.
Thus we have justified all the maps in the following picture:
%\[
%\xymatrix@C=80pt@R=50pt{
% W_{i_{\nu+1}}\ar@/_1pc/[d]|(0.6){g}\ar@<-0.5ex>@/_1.5pc/[d]|{(xy)^{r_2}\quad}\ar[r]|{w_2^{c_{l_{\nu+1}}}w_3^{d_{l_{\nu+1}}}}
%&W_{i_{\nu+2}}\ar@/_1pc/[l]|{(xy)^{r_{l_{\nu+1}}}}\ar@{}[r]|{\hdots} & W_{i_{N-1}}\ar[r]|{w_2^{c_{l_{N-1}}}w_3^{d_{l_{N-1}}}} & W_{i_{N}}\ar@/_1pc/[l]|{(xy)^{r_{l_{N-1}}}}\ar@<0.25ex>[1,-3]|{w_2^{c_{l_{N}}}w_3^{d_{l_{N}}}} \\
%R\ar[u]|{w_2}\ar@<-0.4ex>@/_1pc/[-1,3]|(0.45){(xy)^{r_{l_N}}}&&&}
%\]
\[
\xy0;/r.32pc/:
\POS(0,0)*+{R}="s",(75,15)*+{W_{i_{N}}}="4",(50,15)*+{W_{i_{N-1}}}="3",(25,15)*+{W_{i_{\nu+2}}}="2",(0,15)*+{W_{i_{\nu+1}}}="1"
\ar"s";"1"|{{}_{w_2}}
\ar"1";"2"|{{}_{w_2^{c_{l_{\nu+1}}}w_3^{d_{l_{\nu+1}}}}}
\ar@{.}"2";"3"
\ar@<0.5ex>"3";"4"|{{}_{w_2^{c_{l_{N-1}}}w_3^{d_{l_{N-1}}}}}
\ar@<0.5ex>"4";"s"|{{}_{w_2^{c_{l_{N}}}w_3^{d_{l_{N}}}}}
\ar@<-0.1ex>@/_1pc/"1";"s"|{{}_{g}}
\ar@<-0.2ex>@/_1.75pc/"1";"s"|(0.4){{}_{(xy)^{r_2}}\quad}
\ar@<-0.6ex>@/_0.8pc/"s";"4"|(0.4){\quad{}_{(xy)^{r_{l_N}}}}
\ar@<-0.6ex>@/_1pc/"4";"3"|{{}_{(xy)^{r_{l_{N-1}}}}}
\ar@/_1pc/"2";"1"|{{}_{(xy)^{r_{l_{\nu+1}}}}}
\endxy
\]
In general there will be more maps.  If $\alpha_{\nu+1}>3$ then
for every $2\leq t\leq \alpha_{\nu+1}-2$ add an extra map
$W_{i_{\nu+1}}\rightarrow R$ labelled %\footnote{don't put in terms of $l_{nu}=l_{(\nu+1)-1}$ like below since $l_\nu$ not defined for $\mu=0$.  Its true that $l_{\nu}=2$ for $\nu>0$ so we could define $l_0=2$ but then in all cases $l_{\nu}=2$ so we may as well substitute this in} 
$(xy)^{r_{t+2}}w_2^{c_{t+1}}w_3^{d_{t+1}}$.  For all such $t$ it is true that $b_{t}=\nu+1$ and so these maps go where they should since $(xy)^{r_{t+1}}w_2^{c_{t+1}}w_3^{d_{t+1}}$
is an invariant and $r_{t+1}=r_{t+2}+i_{b_t}$ by
Lemma~\ref{r_butt_lemma}.

Now if $s$ is such that
$\nu+1<s\leq N$ with $\alpha_s>2$, then for every $1\leq t\leq
\alpha_s-2$ add an extra map $W_{i_{s}}\rightarrow R$ labelled $(xy)^{r_{t+l_{s-1}}}w_2^{c_{t-1+l_{s-1}}}w_3^{d_{t-1+l_{s-1}}}$.  For all such $t$ it is true that $b_{t-2+l_{s-1}}=s$ and so by Lemma~\ref{r_butt_lemma} $r_{t-1+l_{s-1}}=r_{t+l_{s-1}}+i_{b_{t-2+l_{s-1}}}=r_{t+l_{s-1}}+i_s$.  Thus these maps also go where they should since
$(xy)^{r_{t-1+l_{s-1}}}w_2^{c_{t-1+l_{s-1}}}w_3^{d_{t-1+l_{s-1}}}$
is an invariant.

\begin{notation}\label{D}
We denote by $\tt{D}_1$ the above quiver with all the extra arrows (if these extra arrows exist).  We denote by $\tt{D}_2$ the quiver obtained from $\tt{D}_1$ by making the substitution $w_2\mapsto v_2$ and $w_3\mapsto v_3$ wherever $w_2$ and $w_3$ appear in the labels of the arrows in $\tt{D}_1$. 
\end{notation}

To ease notation in the proof of Proposition~\ref{loops_in_D1} below we denote the arrows in $\tt{D}_1$ by
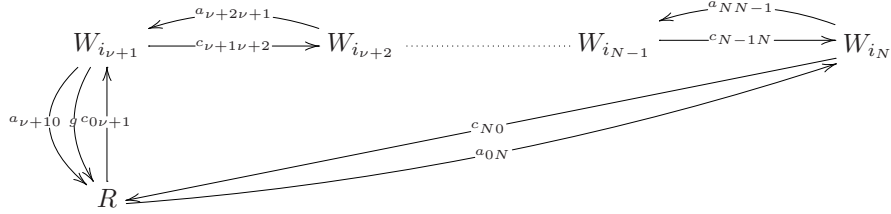
\[
\xy0;/r.32pc/:
\POS(0,0)*+{R}="s",(75,15)*+{W_{i_{N}}}="4",(50,15)*+{W_{i_{N-1}}}="3",(25,15)*+{W_{i_{\nu+2}}}="2",(0,15)*+{W_{i_{\nu+1}}}="1"
\ar"s";"1"|{{}_{\cl{0}{\nu+1}}}
\ar"1";"2"|{{}_{\cl{\nu+1}{\nu+2}}}
\ar@{.}"2";"3"
\ar@<0.5ex>"3";"4"|{{}_{\cl{N-1}{N}}}
\ar@<0.5ex>"4";"s"|{{}_{\cl{N}{0}}}
\ar@<-0.1ex>@/_1pc/"1";"s"|{{}_{g}}
\ar@<-0.2ex>@/_1.75pc/"1";"s"|{{}_{\an{\nu+1}{0}}\quad}
\ar@<-0.6ex>@/_0.8pc/"s";"4"|{{}_{\an{0}{N}}}
\ar@<-0.6ex>@/_1pc/"4";"3"|{{}_{\an{N}{N-1}}}
\ar@/_1pc/"2";"1"|{{}_{\an{\nu+2}{\nu+1}}}
\endxy
\]
Note that $\cl{t}{t+1}$ should read `clockwise from $t$ to $t+1$' and similarly $\an{t}{t-1}$ should read `anticlockwise from $t$ to $t-1$'.  Further we denote the extra arrows (if they exist) by $k_2,k_3,\hdots$ labelled from left to right.

\begin{example}\label{444pattern}
\t{Consider the group $\mathbb{D}_{56,15}$ of order 2460.  Now using 
\[
\frac{56}{56-15}=[2 , 2 ,  3,   2,   3,   2,   2 ]
\]
we can calculate
\parindent 20pt
\parskip 0pt
\[
\begin{array}{c|cccccccc}
& 2& 3&4&5&6&7&8&9\\ \hline
r&67 & 26  &11&  7 &  3 &  2 &  1   &0 \\
c&1  & 0  & 1  & 3  & 5  & 12 & 19 & 26 \\
d&0  & 1 &  1   &2   &3   &7   &11&  15
\end{array}
\]
and so in this example $\nu=0$, and the quiver $\tt{D}_1$ is 
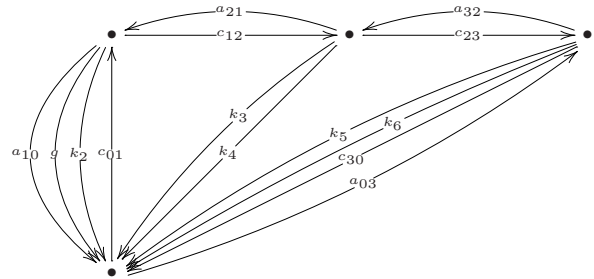
\[
\def\objectstyle{\scriptstyle}
\begin{array}{c}
\xy0;/r.5pc/:
\POS(0,0)*+{\bullet}="s",(30,15)*+{\bullet}="3",(15,15)*+{\bullet}="2",(0,15)*+{\bullet}="1"
\ar"s";"1"|{{}_{\cl{0}{1}}}
\ar"1";"2"|{{}_{\cl{1}{2}}}
\ar"2";"3"|{{}_{\cl{2}{3}}}
\ar@<0.5ex>"3";"s"|{{}_{\cl{3}{0}}}
\ar@<-0.1ex>@/_1.75pc/"1";"s"|{{}_{g}}
\ar@<-0.2ex>@/_2.5pc/"1";"s"|{{}_{\an{1}{0}}\,\,}
\ar@<-0.6ex>@/_0.8pc/"s";"3"|{{}_{\an{0}{3}}}
\ar@/_0.75pc/"3";"2"|{{}_{\an{3}{2}}}
\ar@/_0.75pc/"2";"1"|{{}_{\an{2}{1}}}
\ar@/_1pc/|{{}_{k_2}}"1";"s"
\ar@/_0.6pc/@<-0.25ex>|(0.4){{}_{k_3}}"2";"s"
\ar@<-0.2ex>|{{}_{k_4}}"2";"s"
\ar@<0.3ex>@/_0.8pc/|(0.5){{}_{k_5}}"3";"s"
\ar@<0.4ex>@/_0.4pc/|(0.4){{}_{k_6}}"3";"s"
\endxy
\end{array}
\]
where
\[
\begin{array}{ccc}\begin{array}{l}
k_2= (xy)^{11}w_3\\
k_3= (xy)^{7}w_2^{}w_3^{}\\
k_4= (xy)^{3}w_2^{3}w_3^{2}\\
k_5= (xy)^{2}w_2^{5}w_3^{3}\\
k_6= (xy)w_2^{12}w_3^{7}\\
\end{array}& 
\begin{array}{l}
\cl{0}{1}=w_2\\
\cl{1}{2}= w_2w_3\\
\cl{2}{3}= w_2^{5}w_3^{3}\\
\cl{3}{0}= w_2^{19}w_3^{11}
\end{array}&
\begin{array}{l}
\an{0}{3}=xy\\
\an{3}{2}= (xy)^3\\
\an{2}{1}= (xy)^{11}\\
\an{1}{0}= (xy)^{67}\\
\\
g=(xy)^{52}w_2
\end{array}\end{array}
\]
%\[
%\def\objectstyle{\scriptscriptstyle}
%\begin{array}{cc}
%\begin{array}{c}
%\xy0;/r.3pc/:
%\POS(0,0)*+{R}="s",(30,15)*+{W_1}="3",(15,15)*+{W_4}="2",(0,15)*+{W_{15}}="1"
%\ar"s";"1"|{{}_{w_2}}
%\ar"1";"2"|{{}_{w_2w_3}}
%\ar"2";"3"|{{}_{w_2^5w_3^3}}
%\ar@<0.5ex>"3";"s"|{{}_{\quad w_2^{19}w_3^{11}}}
%\ar@<-0.1ex>@/_1.25pc/"1";"s"|{{}_{g}}
%\ar@<-0.2ex>@/_1.75pc/"1";"s"|(0.4){{}_{(xy)^{67}}\quad}
%\ar@<-0.6ex>@/_0.8pc/"s";"3"|{{}_{xy}}
%\ar@/_0.75pc/"3";"2"|{{}_{(xy)^{3}}}
%\ar@/_0.75pc/"2";"1"|{{}_{(xy)^{11}}}
%\ar@/_0.75pc/|(0.6){{}_{k_2}}"1";"s"
%\ar@/_0.4pc/@<-0.25ex>|(0.4){{}_{k_3}}"2";"s"
%\ar@<-0.2ex>|{{}_{k_4}}"2";"s"
%\ar@<0.3ex>@/_0.6pc/|(0.5){{}_{k_5}}"3";"s"
%\ar@<0.4ex>@/_0.3pc/|(0.4){{}_{k_6}}"3";"s"
%\endxy
%\end{array}
%&
%{\scriptsize{\begin{array}{rcl}
%g&=&(xy)^{52}w_2\\
%k_2&=& (xy)^{11}w_3\\
%k_3&=& (xy)^{7}w_2^{}w_3^{}\\
%k_4&=& (xy)^{3}w_2^{3}w_3^{2}\\
%k_5&=& (xy)^{2}w_2^{5}w_3^{3}\\
%k_6&=& (xy)w_2^{12}w_3^{7}\\
%\end{array}}}
%\end{array}
%\]
Below are the generators of the ring of invariants, and how to view them as paths in $\tt{D}_1$:}
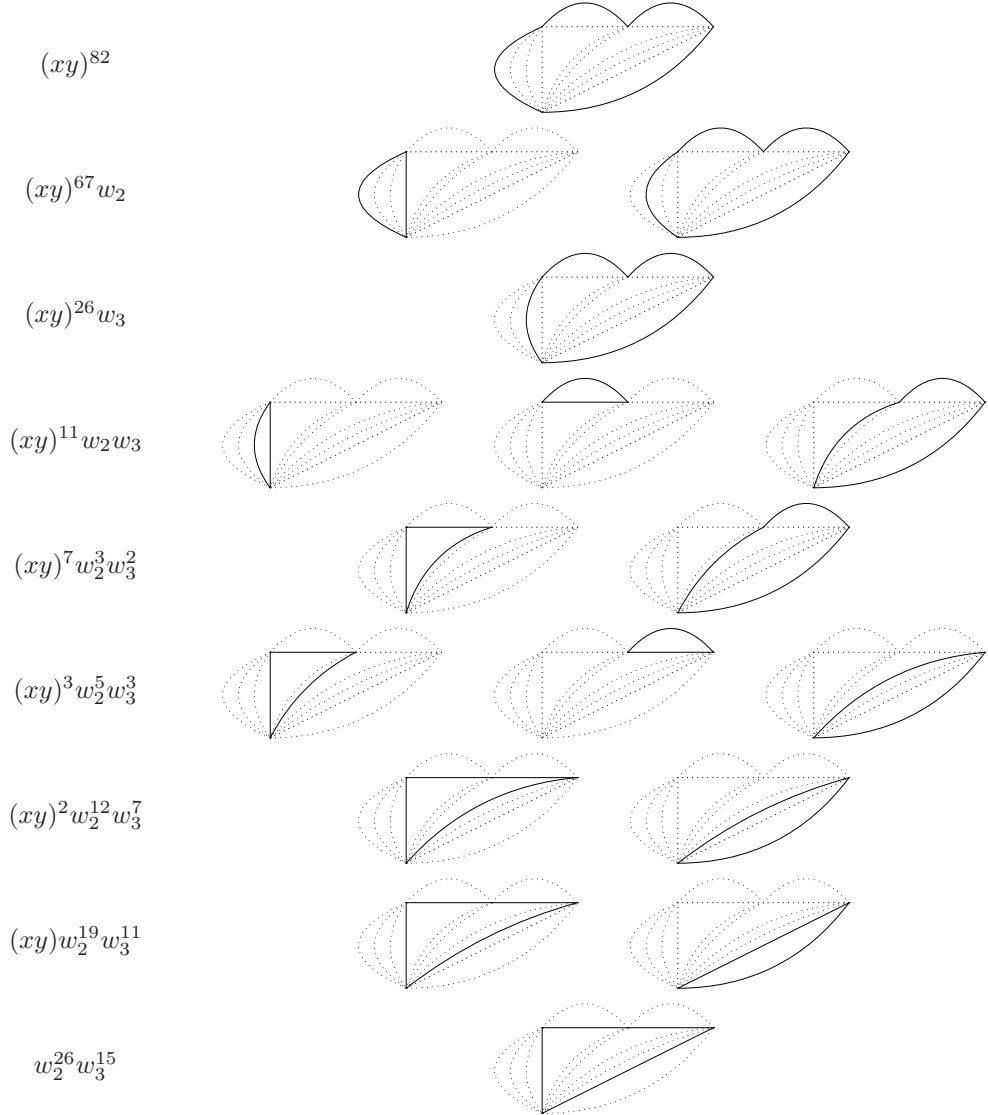
\[
\begin{array}{ccc}
(xy)^{82}&&\begin{array}{c}
\def\objectstyle{\scriptscriptstyle}
\xy0;/r.18pc/:
\POS(0,0)*{}="s",(30,15)*{}="3",(15,15)*{}="2",(0,15)*{}="1"
\ar@{.}"s";"1"
\ar@{.}"1";"2"
\ar@{.}"2";"3"
\ar@{.}"3";"s"

\ar@{.}@/_1pc/"1";"s"
\ar@{-}@/_1.5pc/"1";"s"
\ar@{-}@/_0.75pc/"s";"3"
\ar@{-}@/_0.75pc/"3";"2"
\ar@{-}@/_0.75pc/"2";"1"

\ar@{.}@/_0.5pc/"1";"s"
\ar@{.}@/_0.5pc/"2";"s"
\ar@{.}@/_0.3pc/"2";"s"
\ar@{.}@/_0.6pc/"3";"s"
\ar@{.}@/_0.3pc/"3";"s"
\endxy
\end{array}\\
(xy)^{67}w_2&&\begin{array}{cc} \begin{array}{c}
\def\objectstyle{\scriptscriptstyle}
\xy0;/r.18pc/:
\POS(0,0)*{}="s",(30,15)*{}="3",(15,15)*{}="2",(0,15)*{}="1"
\ar@{-}"s";"1"
\ar@{.}"1";"2"
\ar@{.}"2";"3"
\ar@{.}"3";"s"

\ar@{.}@/_1pc/"1";"s"
\ar@{-}@/_1.5pc/"1";"s"
\ar@{.}@/_0.75pc/"s";"3"
\ar@{.}@/_0.75pc/"3";"2"
\ar@{.}@/_0.75pc/"2";"1"

\ar@{.}@/_0.5pc/"1";"s"
\ar@{.}@/_0.5pc/"2";"s"
\ar@{.}@/_0.3pc/"2";"s"
\ar@{.}@/_0.6pc/"3";"s"
\ar@{.}@/_0.3pc/"3";"s"
\endxy
\end{array}&\begin{array}{c}
\def\objectstyle{\scriptscriptstyle}
\xy0;/r.18pc/:
\POS(0,0)*{}="s",(30,15)*{}="3",(15,15)*{}="2",(0,15)*{}="1"
\ar@{.}"s";"1"
\ar@{.}"1";"2"
\ar@{.}"2";"3"
\ar@{.}"3";"s"

\ar@{-}@/_1pc/"1";"s"
\ar@{.}@/_1.5pc/"1";"s"
\ar@{-}@/_0.75pc/"s";"3"
\ar@{-}@/_0.75pc/"3";"2"
\ar@{-}@/_0.75pc/"2";"1"

\ar@{.}@/_0.5pc/"1";"s"
\ar@{.}@/_0.5pc/"2";"s"
\ar@{.}@/_0.3pc/"2";"s"
\ar@{.}@/_0.6pc/"3";"s"
\ar@{.}@/_0.3pc/"3";"s"
\endxy
\end{array} \end{array}\\
(xy)^{26}w_3 && \begin{array}{c}
\def\objectstyle{\scriptscriptstyle}
\xy0;/r.18pc/:
\POS(0,0)*{}="s",(30,15)*{}="3",(15,15)*{}="2",(0,15)*{}="1"
\ar@{.}"s";"1"
\ar@{.}"1";"2"
\ar@{.}"2";"3"
\ar@{.}"3";"s"

\ar@{.}@/_1pc/"1";"s"
\ar@{.}@/_1.5pc/"1";"s"
\ar@{-}@/_0.75pc/"s";"3"
\ar@{-}@/_0.75pc/"3";"2"
\ar@{-}@/_0.75pc/"2";"1"

\ar@{-}@/_0.5pc/"1";"s"
\ar@{.}@/_0.5pc/"2";"s"
\ar@{.}@/_0.3pc/"2";"s"
\ar@{.}@/_0.6pc/"3";"s"
\ar@{.}@/_0.3pc/"3";"s"
\endxy
\end{array}\\
 (xy)^{11}w_2w_3 &&\begin{array}{ccc} \begin{array}{c}
\def\objectstyle{\scriptscriptstyle}
\xy0;/r.18pc/:
\POS(0,0)*{}="s",(30,15)*{}="3",(15,15)*{}="2",(0,15)*{}="1"
\ar@{-}"s";"1"
\ar@{.}"1";"2"
\ar@{.}"2";"3"
\ar@{.}"3";"s"

\ar@{.}@/_1pc/"1";"s"
\ar@{.}@/_1.5pc/"1";"s"
\ar@{.}@/_0.75pc/"s";"3"
\ar@{.}@/_0.75pc/"3";"2"
\ar@{.}@/_0.75pc/"2";"1"

\ar@{-}@/_0.5pc/"1";"s"
\ar@{.}@/_0.5pc/"2";"s"
\ar@{.}@/_0.3pc/"2";"s"
\ar@{.}@/_0.6pc/"3";"s"
\ar@{.}@/_0.3pc/"3";"s"
\endxy
\end{array}&
\begin{array}{c}
\def\objectstyle{\scriptscriptstyle}
\xy0;/r.18pc/:
\POS(0,0)*{}="s",(30,15)*{}="3",(15,15)*{}="2",(0,15)*{}="1"
\ar@{.}"s";"1"
\ar@{-}"1";"2"
\ar@{.}"2";"3"
\ar@{.}"3";"s"

\ar@{.}@/_1pc/"1";"s"
\ar@{.}@/_1.5pc/"1";"s"
\ar@{.}@/_0.75pc/"s";"3"
\ar@{.}@/_0.75pc/"3";"2"
\ar@{-}@/_0.75pc/"2";"1"

\ar@{.}@/_0.5pc/"1";"s"
\ar@{.}@/_0.5pc/"2";"s"
\ar@{.}@/_0.3pc/"2";"s"
\ar@{.}@/_0.6pc/"3";"s"
\ar@{.}@/_0.3pc/"3";"s"
\endxy
\end{array}&
\begin{array}{c}
\def\objectstyle{\scriptscriptstyle}
\xy0;/r.18pc/:
\POS(0,0)*{}="s",(30,15)*{}="3",(15,15)*{}="2",(0,15)*{}="1"
\ar@{.}"s";"1"
\ar@{.}"1";"2"
\ar@{.}"2";"3"
\ar@{.}"3";"s"

\ar@{.}@/_1pc/"1";"s"
\ar@{.}@/_1.5pc/"1";"s"
\ar@{-}@/_0.75pc/"s";"3"
\ar@{-}@/_0.75pc/"3";"2"
\ar@{.}@/_0.75pc/"2";"1"

\ar@{.}@/_0.5pc/"1";"s"
\ar@{-}@/_0.5pc/"2";"s"
\ar@{.}@/_0.3pc/"2";"s"
\ar@{.}@/_0.6pc/"3";"s"
\ar@{.}@/_0.3pc/"3";"s"
\endxy
\end{array} \end{array}\\
(xy)^{7}w_2^3w_3^2 &&\begin{array}{cc} \begin{array}{c}
\def\objectstyle{\scriptscriptstyle}
\xy0;/r.18pc/:
\POS(0,0)*{}="s",(30,15)*{}="3",(15,15)*{}="2",(0,15)*{}="1"
\ar@{-}"s";"1"
\ar@{-}"1";"2"
\ar@{.}"2";"3"
\ar@{.}"3";"s"

\ar@{.}@/_1pc/"1";"s"
\ar@{.}@/_1.5pc/"1";"s"
\ar@{.}@/_0.75pc/"s";"3"
\ar@{.}@/_0.75pc/"3";"2"
\ar@{.}@/_0.75pc/"2";"1"

\ar@{.}@/_0.5pc/"1";"s"
\ar@{-}@/_0.5pc/"2";"s"
\ar@{.}@/_0.3pc/"2";"s"
\ar@{.}@/_0.6pc/"3";"s"
\ar@{.}@/_0.3pc/"3";"s"
\endxy
\end{array}&
\begin{array}{c}
\def\objectstyle{\scriptscriptstyle}
\xy0;/r.18pc/:
\POS(0,0)*{}="s",(30,15)*{}="3",(15,15)*{}="2",(0,15)*{}="1"
\ar@{.}"s";"1"
\ar@{.}"1";"2"
\ar@{.}"2";"3"
\ar@{.}"3";"s"

\ar@{.}@/_1pc/"1";"s"
\ar@{.}@/_1.5pc/"1";"s"
\ar@{-}@/_0.75pc/"s";"3"
\ar@{-}@/_0.75pc/"3";"2"
\ar@{.}@/_0.75pc/"2";"1"

\ar@{.}@/_0.5pc/"1";"s"
\ar@{.}@/_0.5pc/"2";"s"
\ar@{-}@/_0.3pc/"2";"s"
\ar@{.}@/_0.6pc/"3";"s"
\ar@{.}@/_0.3pc/"3";"s"
\endxy
\end{array}  \end{array}\\
(xy)^{3}w_2^5w_3^3  &&\begin{array}{ccc} \begin{array}{c}
\def\objectstyle{\scriptscriptstyle}
\xy0;/r.18pc/:
\POS(0,0)*{}="s",(30,15)*{}="3",(15,15)*{}="2",(0,15)*{}="1"
\ar@{-}"s";"1"
\ar@{-}"1";"2"
\ar@{.}"2";"3"
\ar@{.}"3";"s"

\ar@{.}@/_1pc/"1";"s"
\ar@{.}@/_1.5pc/"1";"s"
\ar@{.}@/_0.75pc/"s";"3"
\ar@{.}@/_0.75pc/"3";"2"
\ar@{.}@/_0.75pc/"2";"1"

\ar@{.}@/_0.5pc/"1";"s"
\ar@{.}@/_0.5pc/"2";"s"
\ar@{-}@/_0.3pc/"2";"s"
\ar@{.}@/_0.6pc/"3";"s"
\ar@{.}@/_0.3pc/"3";"s"
\endxy
\end{array}&
\begin{array}{c}
\def\objectstyle{\scriptscriptstyle}
\xy0;/r.18pc/:
\POS(0,0)*{}="s",(30,15)*{}="3",(15,15)*{}="2",(0,15)*{}="1"
\ar@{.}"s";"1"
\ar@{.}"1";"2"
\ar@{-}"2";"3"
\ar@{.}"3";"s"

\ar@{.}@/_1pc/"1";"s"
\ar@{.}@/_1.5pc/"1";"s"
\ar@{.}@/_0.75pc/"s";"3"
\ar@{-}@/_0.75pc/"3";"2"
\ar@{.}@/_0.75pc/"2";"1"

\ar@{.}@/_0.5pc/"1";"s"
\ar@{.}@/_0.5pc/"2";"s"
\ar@{.}@/_0.3pc/"2";"s"
\ar@{.}@/_0.6pc/"3";"s"
\ar@{.}@/_0.3pc/"3";"s"
\endxy
\end{array} &
\begin{array}{c}
\def\objectstyle{\scriptscriptstyle}
\xy0;/r.18pc/:
\POS(0,0)*{}="s",(30,15)*{}="3",(15,15)*{}="2",(0,15)*{}="1"
\ar@{.}"s";"1"
\ar@{.}"1";"2"
\ar@{.}"2";"3"
\ar@{.}"3";"s"

\ar@{.}@/_1pc/"1";"s"
\ar@{.}@/_1.5pc/"1";"s"
\ar@{-}@/_0.75pc/"s";"3"
\ar@{.}@/_0.75pc/"3";"2"
\ar@{.}@/_0.75pc/"2";"1"

\ar@{.}@/_0.5pc/"1";"s"
\ar@{.}@/_0.5pc/"2";"s"
\ar@{.}@/_0.3pc/"2";"s"
\ar@{-}@/_0.6pc/"3";"s"
\ar@{.}@/_0.3pc/"3";"s"
\endxy
\end{array} \end{array}\\
%\end{array}
%\]
%\[
%\begin{array}{ccc}
(xy)^2w_2^{12}w_3^7 &&\begin{array}{cc}  \begin{array}{c}
\def\objectstyle{\scriptscriptstyle}
\xy0;/r.18pc/:
\POS(0,0)*{}="s",(30,15)*{}="3",(15,15)*{}="2",(0,15)*{}="1"
\ar@{-}"s";"1"
\ar@{-}"1";"2"
\ar@{-}"2";"3"
\ar@{.}"3";"s"

\ar@{.}@/_1pc/"1";"s"
\ar@{.}@/_1.5pc/"1";"s"
\ar@{.}@/_0.75pc/"s";"3"
\ar@{.}@/_0.75pc/"3";"2"
\ar@{.}@/_0.75pc/"2";"1"

\ar@{.}@/_0.5pc/"1";"s"
\ar@{.}@/_0.5pc/"2";"s"
\ar@{.}@/_0.3pc/"2";"s"
\ar@{-}@/_0.6pc/"3";"s"
\ar@{.}@/_0.3pc/"3";"s"
\endxy
\end{array}&
\begin{array}{c}
\def\objectstyle{\scriptscriptstyle}
\xy0;/r.18pc/:
\POS(0,0)*{}="s",(30,15)*{}="3",(15,15)*{}="2",(0,15)*{}="1"
\ar@{.}"s";"1"
\ar@{.}"1";"2"
\ar@{.}"2";"3"
\ar@{.}"3";"s"

\ar@{.}@/_1pc/"1";"s"
\ar@{.}@/_1.5pc/"1";"s"
\ar@{-}@/_0.75pc/"s";"3"
\ar@{.}@/_0.75pc/"3";"2"
\ar@{.}@/_0.75pc/"2";"1"

\ar@{.}@/_0.5pc/"1";"s"
\ar@{.}@/_0.5pc/"2";"s"
\ar@{.}@/_0.3pc/"2";"s"
\ar@{.}@/_0.6pc/"3";"s"
\ar@{-}@/_0.3pc/"3";"s"
\endxy
\end{array} \end{array} \\
(xy)w_2^{19}w_3^{11} &&\begin{array}{cc} \begin{array}{c}
\def\objectstyle{\scriptscriptstyle}
\xy0;/r.18pc/:
\POS(0,0)*{}="s",(30,15)*{}="3",(15,15)*{}="2",(0,15)*{}="1"
\ar@{-}"s";"1"
\ar@{-}"1";"2"
\ar@{-}"2";"3"
\ar@{.}"3";"s"

\ar@{.}@/_1pc/"1";"s"
\ar@{.}@/_1.5pc/"1";"s"
\ar@{.}@/_0.75pc/"s";"3"
\ar@{.}@/_0.75pc/"3";"2"
\ar@{.}@/_0.75pc/"2";"1"

\ar@{.}@/_0.5pc/"1";"s"
\ar@{.}@/_0.5pc/"2";"s"
\ar@{.}@/_0.3pc/"2";"s"
\ar@{.}@/_0.6pc/"3";"s"
\ar@{-}@/_0.3pc/"3";"s"
\endxy
\end{array}&
\begin{array}{c}
\def\objectstyle{\scriptscriptstyle}
\xy0;/r.18pc/:
\POS(0,0)*{}="s",(30,15)*{}="3",(15,15)*{}="2",(0,15)*{}="1"
\ar@{.}"s";"1"
\ar@{.}"1";"2"
\ar@{.}"2";"3"
\ar@{-}"3";"s"

\ar@{.}@/_1pc/"1";"s"
\ar@{.}@/_1.5pc/"1";"s"
\ar@{-}@/_0.75pc/"s";"3"
\ar@{.}@/_0.75pc/"3";"2"
\ar@{.}@/_0.75pc/"2";"1"

\ar@{.}@/_0.5pc/"1";"s"
\ar@{.}@/_0.5pc/"2";"s"
\ar@{.}@/_0.3pc/"2";"s"
\ar@{.}@/_0.6pc/"3";"s"
\ar@{.}@/_0.3pc/"3";"s"
\endxy
\end{array}  \end{array}\\
w_2^{26}w_3^{15}&&\begin{array}{c}
\def\objectstyle{\scriptscriptstyle}
\xy0;/r.18pc/:
\POS(0,0)*{}="s",(30,15)*{}="3",(15,15)*{}="2",(0,15)*{}="1"
\ar@{-}"s";"1"
\ar@{-}"1";"2"
\ar@{-}"2";"3"
\ar@{-}"3";"s"

\ar@{.}@/_1pc/"1";"s"
\ar@{.}@/_1.5pc/"1";"s"
\ar@{.}@/_0.75pc/"s";"3"
\ar@{.}@/_0.75pc/"3";"2"
\ar@{.}@/_0.75pc/"2";"1"

\ar@{.}@/_0.5pc/"1";"s"
\ar@{.}@/_0.5pc/"2";"s"
\ar@{.}@/_0.3pc/"2";"s"
\ar@{.}@/_0.6pc/"3";"s"
\ar@{.}@/_0.3pc/"3";"s"
\endxy
\end{array}
\end{array}
\]
\end{example}

The next proposition follows the above example closely; in fact the cycle pattern is the same as in type $A$, but the proof here is a little more complicated than in type $A$ since we need to rely heavily on the combinatorics from Section 2.  In what follows it is convenient to write $\AN{i}{j}$ for the composition of
anticlockwise paths $a$ from vertex $i$ to vertex $j$, and
similarly $\CL{i}{j}$ as the composition of clockwise paths, where
by $\CL{i}{i}$ (resp. $\AN{i}{i}$) we mean not the empty path at
vertex $i$ but the path from $i$ to $i$ round each of the
clockwise (resp. anticlockwise) arrows precisely once.  So, in Example~\ref{444pattern} above, $(xy)^{82}=A_{00}$ and $w_2^{26}w_3^{15}=C_{00}$.  We also refer to the vertex $W_{i_t}$ as vertex $t$, and the vertex $R$ as $\star$.

\begin{prop}\label{loops_in_D1}
At every vertex in the quiver $\tt{D}_1$, all invariant polynomials exist as a sum of compositions of cycles at that vertex.  Furthermore, if we remove any one arrow this is no longer true.  Both statements also hold for $\tt{D}_2$.
\end{prop}
\begin{proof}
We first prove the statements for $\tt{D}_1$.  To simplify the exposition we consider cycles up to rotation --- for example by `$\CL{0}{0}$ at vertex $t$' we mean $\CL{t}{t}$ (i.e. we rotate suitably so that the cycle starts and finishes at the vertex we want).  It is clear that we can see $(xy)^{2(n-q)}$ at every vertex as $\AN{0}{0}$.  By Lemma~\ref{invariants} it suffices to justify why at every vertex there is also the cycle $(xy)^{r_t}w_2^{c_t}w_3^{d_t}$ for all $2\leq t\leq e$. 

First, consider $t=2$, i.e. $(xy)^{r_2}w_2$.  At the vertex $\star$ and the vertex ${\nu+1}$ we can see this as $\cl{0}{\nu+1}\an{\nu+1}{0}$.  At all other vertices we can see $(xy)^{r_2}w_2$ as $\AN{0}{\nu+1}g$  since by Lemma~\ref{a2,q,r2,r3,r2}(v) $r_2=2r_3+i_{\nu+1}$.  

Now consider $t=3$, i.e. $(xy)^{r_3}w_2^{c_3}w_3^{d_3}$.  At vertices $b_2,\hdots,N$ and $\star$ we can view this invariant as $\AN{0}{b_2}k_2$ since $\AN{0}{b_2}k_2=(xy)^{i_{b_2}+r_4}w_2^{c_3}w_3^{d_3}=(xy)^{r_3}w_2^{c_3}w_3^{d_3}$ by Lemma~\ref{r_butt_lemma}.  Thus if $b_2=\nu+1$ these are all the vertices, else $b_2>\nu+1$.  In this case we can view  $(xy)^{r_3}w_2^{c_3}w_3^{d_3}$ at vertices $\nu+1\leq s<b_2$ as $\cl{s}{s+1}\an{s+1}{s}$ since for these $s$ by definition $l_s=3$. 

From here the pattern in the invariants is exactly the same way as in type $A$ --- consider now $(xy)^{r_t}w_2^{c_t}w_3^{d_t}$ with $4\leq t\leq e-1$.  At the vertices $b_{t-1},\hdots, N$ and $\star$ can see this invariant as $\AN{0}{b_{t-1}}k_{t-1}$ since $\AN{0}{b_{t-1}}k_{t-1}=(xy)^{i_{b_{t-1}}+r_{t+1}}w_2^{c_{t}}w_3^{d_{t}}=(xy)^{r_{t}}w_2^{c_{t}}w_3^{d_{t}}$ by Lemma~\ref{r_butt_lemma}.  At vertices $\nu+1,\hdots,b_{t-2}$ and $\star$ the invariant may be viewed as $\CL{0}{b_{t-2}}k_{t-2}$ since $\CL{0}{b_{t-2}}k_{t-2}=(xy)^{r_t}w_2^{c_{t-1}+\Delta_{b_{t-2}}}w_3^{d_{t-1}+\Delta_{b_{t-2}}}=(xy)^{r_t}w_2^{c_t}w_3^{d_t}$ by Lemma~\ref{c_and_d_lemma}.  Thus if $b_{t-2}=b_{t-1}$ these are all vertices, else $b_{t-2}<b_{t-1}$.  In this case we can view  the invariant at vertices $b_{t-2}\leq s<b_{t-1}$ as $\cl{s}{s+1}\an{s+1}{s}$ since for these $s$ we have $l_s={t}$. 

The last invariant $(xy)^{r_e}w_2^{c_e}w_3^{d_e}=w_2^{c_e}w_3^{d_e}$ can be seen at all vertices as $\CL{0}{0}$.  This is because $l_N=e-1$ and so the $\cl{N}{0}=w_2^{c_{e-1}}w_3^{d_{e-1}}$.  Consequently $\CL{0}{0}=w_2^{c_{e-1}+\Delta_{N}}w_3^{d_{e-1}+\Gamma_{N}}$ which is $w_2^{c_e}w_3^{d_e}$ by Lemma~\ref{c_and_d_lemma} since $b_{e-2}=b_{l_N-1}$ is defined to be $N$.  

Since by Lemma~\ref{invariants}  the collection $(xy)^{2(n-q)}$, $(xy)^{r_t}w_2^{c_t}w_3^{d_t}$ for all $2\leq t\leq e$ generate $\C{}[x,y]^{\mathbb{D}_{n,q}}$ it follows that we can see, at each vertex, all invariants.  If we remove any one arrow from $\tt{D}_1$ then by inspection of the above proof this is no longer true.  To obtain the final statement regarding $\tt{D}_2$, note that in the above proof if we make the substitution $w_2\mapsto v_2$ and $w_3\mapsto v_3$ throughout and appeal to Theorem~\ref{Riemen_generates} instead of Lemma~\ref{invariants}, nothing in the proof is affected.
\end{proof}

The above may be a simple observation, but it allows us to diagram chase over $\tt{D}_1$ to deduce that $W_{i_{\nu+1}},\hdots,W_{i_N}$ are 2-generated and thus (by Lemma~\ref{2genlemma}) are special:

\begin{thm}\label{Wi's_are_special_most}
For $\mathbb{D}_{n,q}$ with $0\leq \nu<N-1$ and all $1\leq t\leq
\nu+1$, $W_{i_t}$ is generated as a $\C{}[x,y]^{\mathbb{D}_{n,q}}$
module by $(xy)^{i_{t}}$ and $w_2^{\Delta_t}w_3^{\Gamma_t}$, and
so is special.  Alternatively we can take as generators $(xy)^{i_{t}}$ and $v_2^{\Delta_t}v_3^{\Gamma_t}$.
\end{thm}
\begin{proof}
We restrict ourselves to proving the generators $(xy)^{i_{t}}$ and $w_2^{\Delta_t}w_3^{\Gamma_t}$ by using $\tt{D}_1$; the other generators follow immediately from the proof below by making the substitutions $w_2\mapsto v_2$ and $w_3\mapsto v_3$ throughout and working with $\tt{D}_2$ instead.  

We first verify the case $W_{i_N}$, then proceed by induction on decreasing $t$.  To prove the $W_{i_N}$ case, let $f\in W_{i_N}$, and split into 2 subcases:\\
\emph{Subcase 1:} $\alpha_N=2$.   Consider $w_2^{c_{l_N}}w_3^{d_{l_N}}f$.  It is an invariant, so we can view it as a sum of cycles at the vertex $W_{i_N}$.  These must all leave the vertex, so since there are only 2 arrows out we can write
\[
w_2^{c_{l_N}}w_3^{d_{l_N}}f=w_2^{c_{l_N}}w_3^{d_{l_N}}p_{0,N}+(xy)^{i_{N-1}-i_N}p_{N-1,N}
\]
where $p_{0,N}$ is a sum of paths from $R$ to $W_{i_N}$ and $p_{N-1,N}$ is a sum of paths from $W_{i_{N-1}}$ to $W_{i_N}$.  Note that $i_{N-1}-i_N=1$ since $\alpha_N=2$.  Viewing everything as polynomials $w_2^{c_{l_N}}w_3^{d_{l_N}}$ must divide $p_{N-1,N}$ and so after cancelling the $w_2^{c_{l_N}}w_3^{d_{l_N}}$ term we may write
\[
f=p_{0,N}+(xy)A
\]
for some polynomial $A$.  By inspection of the quiver $\tt{D}_1$ there are only two paths from $R$ to $W_{i_N}$ that don't involve cycles, so after moving all cycles to the end of the path (which we can do since there are all invariants at all vertices) we may write $p_{0,N}$ as
\[
p_{0,N}=(xy)^{r_{l_N}}B_1+(w_2^{\Delta_N}w_3^{\Gamma_N})B_2
\]
for some polynomials $B_1$ and $B_2$.  Thus since $r_{l_N}=1$ we see that
\[
f=(xy)(B_1+A)+(w_2^{\Delta_N}w_3^{\Gamma_N})B_2
\]
thus by Lemma~\ref{polys_to_invariants} it follows that $W_{i_N}$ is generated by $xy=(xy)^{i_N}$ and $w_2^{\Delta_N}w_3^{\Gamma_N}$, and so is special.\\
\emph{Subcase 2:} $\alpha_N>2$.  In this subcase there is an extra arrow from $W_{i_N}$ to $R$ labelled by $(xy)^{r_{1+l_{N-1}}}w_2^{c_{l_{N-1}}}w_3^{d_{l_{N-1}}}$.  Further by Lemma~\ref{r_is_difference_in_i_series} and Lemma~\ref{r_butt_lemma}
\[
 i_{N-1}-i_N=r_{l_{N-1}}=r_{1+l_{N-1}}+i_{b_{-1+l_{N-1}}}
\]
thus since $b_{-1+l_{N-1}}=N$ and $i_N=1$ it follows that $r_{1+l_{N-1}}=(i_{N-1}-i_N)-1$.  Consequently the polynomial $(xy)^{(i_{N-1}-i_N)-1}w_2^{c_{l_{N-1}}}w_3^{d_{l_{N-1}}}f$ is an invariant and so we can view it as a sum of cycles at the vertex $W_{i_N}$.  They all must leave the vertex, so we may write
\begin{multline*}
(xy)^{(i_{N-1}-i_N)-1}w_2^{c_{l_{N-1}}}w_3^{d_{l_{N-1}}}f=\\ (xy)^{i_{N-1}-i_N}A_1+(xy)^{(i_{N-1}-i_N)-1}w_2^{c_{l_{N-1}}}w_3^{d_{l_{N-1}}}p_{0,N}+w_2^{\Delta_t}w_3^{\Gamma_t}B_1
\end{multline*}
where $p_{0,N}$ is a sum of paths from $R$ to $W_{i_N}$,  $A_1$ is a sum of paths from $W_{i_{N-1}}$ to $W_{i_N}$,  and $B_1$ is some polynomial; we can do this since all the other arrows leaving $W_{i_N}$ are divisible by $w_2^{\Delta_t}w_3^{\Gamma_t}$ by Lemma~\ref{c_and_d_lemma}. %\footnote{its really just saying the reconstruction relations}  
Viewing everything as polynomials, $(xy)^{(i_{N-1}-i_N)-1}$ must divide $B_1$ and further  $w_2^{c_{l_{N-1}}}w_3^{d_{l_{N-1}}}$ must divide $A_1$, thus after cancelling these terms we see
\[
f=(xy)A_2+p_{0,N}+w_2^{\Delta_t}w_3^{\Gamma_t}B_2
\]
for some polynomials $A_2$ and $B_2$.  By inspection of the quiver $\tt{D}_1$ there are only two paths from $R$ to $W_{i_N}$ that don't involve cycles (one is $xy=(xy)^{r_{l_N}}$, the other is $w_2^{\Delta_N}w_3^{\Gamma_N}$), so after moving all cycles to the end of the path (which we can do since there are all invariants at all vertices) we may write $p_{0,N}$ as
\[
p_{0,N}=(xy)^{r_{l_N}}C_1+(w_2^{\Delta_N}w_3^{\Gamma_N})C_2
\]
for some polynomials $C_1$ and $C_2$.  Thus since $r_{l_N}=1$ we see that
\[
f=(xy)(C_1+A_2)+(w_2^{\Delta_N}w_3^{\Gamma_N})(B_2+C_2)
\]
and so by Lemma~\ref{polys_to_invariants} it follows that $W_{i_N}$ is generated by $xy=(xy)^{i_N}$ and $w_2^{\Delta_N}w_3^{\Gamma_N}$, hence is special.\\
For the induction step suppose we are considering $W_{i_t}$ with $\nu+1\leq t<N$ and we have established the result for $W_{i_{t+1}}$.  Let $f\in W_{i_t}$ and consider $w_2^{c_{l_t}}w_3^{d_{l_t}}f\in W_{i_{t+1}}$. %\footnote{going this way avoids `edge issues'}  
By inductive hypothesis we may write
\[
w_2^{c_{l_t}}w_3^{d_{l_t}}f=(xy)^{i_{t+1}}A_1 +w_2^{\Delta_{t+1}}w_3^{\Gamma_{t+1}}B_1
\]
for some invariant polynomials $A_1,B_1$.  Viewing everything as polynomials we see that $w_2^{c_{l_t}}w_3^{d_{l_t}}$ divides $A_1$ and so after cancelling this factor
\[
f=(xy)^{i_{t+1}}A_2+w_2^{\Delta_{t}}w_3^{\Gamma_{t}}B_1
\]
for some polynomial $A_2$.  Since $B_1$ is invariant $w_2^{\Delta_{t}}w_3^{\Gamma_{t}}B_1\in W_{i_t}$, hence since $f\in W_{i_t}$ we get $(xy)^{i_{t+1}}A_2\in W_{i_t}$.  This is turn implies that $A_2w_2^{c_{l_t}}w_3^{d_{l_t}}$ is an invariant, and so we can view it as a sum of cycles at vertex $W_{i_{t+1}}$.  From here we split into 2 subcases:\\ %\footnote{The following depends on $t+1$, thus largely avoiding side issues.  However can't go past $\nu+1$ and get other specials since to get to this stage we use $fw_2^{c_{l_t}}w_3^{d_{l_t}}\in W_{i_{t+1}}$ which needs the assumption greater than $\nu+1$ so that we can use previous combinatorics} 
\emph{Subcase 1:} $\alpha_{t+1}=2$.  Then there are only two arrows out of $W_{i_{t+1}}$, thus
\[
A_2w_2^{c_{l_t}}w_3^{d_{l_t}}=(xy)^{i_t-i_{t+1}}p_{t,t+1}+w_2^{c_{l_{t+1}}}w_3^{d_{l_{t+1}}}p_{t+2,t+1}
\]
where $p_{t,t+1}$ is a sum of paths from $W_{i_t}$ to $W_{i_{t+1}}$ and $p_{t+2,t+1}$ is a sum of paths from $W_{i_{t+2}}$ to $W_{i_{t+1}}$.
But since $\alpha_{t+1}=2$ it follows that $c_{l_{t+1}}=c_{l_t}$ and $d_{l_{t+1}}=d_{l_t}$, thus viewing everything as polynomials $w_2^{c_{l_t}}w_3^{d_{l_t}}$ divides $p_{t,t+1}$ and so after cancelling this factor
\[
A_2=(xy)^{i_t-i_{t+1}}D_1+p_{t+2,t+1}
\]
for some polynomial $D_1$.  Now any path from $W_{i_{t+2}}$ to $W_{i_{t+1}}$ must either factor through the map $(xy)^{i_{t+1}-i_{t+2}}$ or go via $R$ and end through the composition of maps $w_2^{\Delta_{t+1}}w_3^{\Gamma_{t+1}}$, thus we may write $p_{t+2,t+1}$ as
\[
p_{t+2,t+1}=(xy)^{i_{t+1}-i_{t+2}}E_1+w_2^{\Delta_{t+1}}w_3^{\Gamma_{t+1}}E_2
\]
for some polynomials $E_1$ and $E_2$.  But $i_t-i_{t+1}=i_{t+1}-i_{t+2}$ since $\alpha_{t+1}=2$, hence
\[
A_2=(xy)^{i_t-i_{t+1}}(D_1+E_1)+w_2^{\Delta_{t+1}}w_3^{\Gamma_{t+1}}E_2.
\]
Consequently
\begin{eqnarray*}
f&=&(xy)^{i_{t+1}}((xy)^{i_t-i_{t+1}}(D_1+E_1)+w_2^{\Delta_{t+1}}w_3^{\Gamma_{t+1}}E_2)+w_2^{\Delta_{t}}w_3^{\Gamma_{t}}B_1\\
&=& (xy)^{i_t}F_1+w_2^{\Delta_{t}}w_3^{\Gamma_{t}}F_2
\end{eqnarray*}
for some polynomials $F_1$ and $F_2$.  Hence by Lemma~\ref{polys_to_invariants} it follows that $W_{i_t}$ is generated by $(xy)^{i_t}$ and $w_2^{\Delta_t}w_3^{\Gamma_t}$, thus is special.\\
\emph{Subcase 2:} $\alpha_{t+1} >2$.  Here we may write
\[
A_2w_2^{c_{l_t}}w_3^{d_{l_t}}=(xy)^{i_t-i_{t+1}}p_{t,t+1}+(xy)^{i_t-2i_{t+1}}w_2^{c_{l_t}}w_3^{d_{l_t}}p_{0,t+1} +w_2^{\Delta_{t+1}}w_3^{\Gamma_{t+1}}G_1
\]
where $p_{t,t+1}$ is a sum of paths from $W_{i_t}$ to $W_{i_{t+1}}$, $p_{0,t+1}$ is a sum of paths from $R$ to $W_{i_{t+1}}$ and $G_1$ is some polynomial; we can do this since by Lemma~\ref{c_and_d_lemma} %\footnote{again its really just saying the reconstruction relations} 
all other paths out of $W_{i_{t+1}}$ when viewed as polynomials are divisible by $w_2^{\Delta_{t+1}}w_3^{\Gamma_{t+1}}$.  We are also using the fact that $r_{1+l_{t}}=i_{t}-2i_{t+1}$, which is true %\footnote{very similar proof to Case 2 in proof of $W_{i_N}$} 
by Lemma~\ref{r_is_difference_in_i_series} and Lemma~\ref{r_butt_lemma}.  By inspection of $\tt{D}_1$ we may write
\[
p_{0,t+1}=(xy)^{i_{t+1}}H_1+w_2^{\Delta_{t+1}}w_3^{\Gamma_{t+1}}H_2
\]
for some polynomials $H_1, H_2$ and so
\begin{eqnarray*}
A_2w_2^{c_{l_t}}w_3^{d_{l_t}}&=&(xy)^{i_t-i_{t+1}}(p_{t,t+1}+w_2^{c_{l_t}}w_3^{d_{l_t}}H_1)+ w_2^{\Delta_{t+1}}w_3^{\Gamma_{t+1}}(G_1+(xy)^{i_t-2i_{t+1}}H_2)\\
&=& (xy)^{i_t-i_{t+1}}J_1+ w_2^{\Delta_{t+1}}w_3^{\Gamma_{t+1}}K_1
\end{eqnarray*}
say.  Thus $w_2^{c_{l_t}}w_3^{d_{l_t}}$ must divide $J_1$ and so after cancelling this factor we get
\[
A_2=(xy)^{i_t-i_{t+1}}J_2+ w_2^{\Delta_{t}}w_3^{\Gamma_{t}}K_1
\]
for some polynomial $J_2$.  Consequently
\begin{eqnarray*}
f&=&(xy)^{i_{t+1}}((xy)^{i_t-i_{t+1}}J_2+ w_2^{\Delta_{t}}w_3^{\Gamma_{t}}K_1)+w_2^{\Delta_{t}}w_3^{\Gamma_{t}}B_1\\
&=& (xy)^{i_t}J_2+w_2^{\Delta_{t}}w_3^{\Gamma_{t}}(B_1+(xy)^{i_{t+1}}K_1)
\end{eqnarray*}
and so by Lemma~\ref{polys_to_invariants} it follows that $W_{i_t}$ is generated by $(xy)^{i_t}$ and $w_2^{\Delta_t}w_3^{\Gamma_t}$, thus is special.

This establishes the induction step, so the result now follows.
\end{proof}

\noindent
\textbf{Case 2: $\nu=N-1$}.
\begin{thm}\label{Wi's_are_special_nu_top}
For $\mathbb{D}_{n,q}$ with $\nu=N-1$, $W_{i_N}=W_{1}$ is generated as a $\C{}[x,y]^{\mathbb{D}_{n,q}}$ module by $xy$ and $w_2$, thus is special.  Alternatively we can take as generators $xy$ and $v_2$.
\end{thm}
\begin{proof}
We prove the statement regarding $xy$ and $w_2$; the other statement follows immediately since either $\nu=N-1$ is odd and so $v_2=w_2$, or its even in which case (since $q=i_{\nu+1}+\nu(n-q)=1+\nu(n-q)$ by Lemma~\ref{a2,q,r2,r3,r2}(ii)) $v_2=w_2-2(xy)((xy)^{2(n-q)})^{\frac{\nu}{2}}$ with $((xy)^{2(n-q)})^{\frac{\nu}{2}}$ an invariant. 

If $\alpha_N=2$ then the group is in $\SL(2,\C{})$ and so the result is well known, %\footnote{can't do SL case using this technique here as $w_3$ is invariant; need GL so that all invariants are either divisible by $xy$ or $w_2$; this is why we get all invariants}
hence we can assume that $\alpha_N>2$.  It is clear that we have the following maps from $R$ to $W_{i_N}$, and from $W_{i_N}$ to $R$:
\[
\xymatrix@R=80pt@C=100pt{R\ar@/^1.5pc/[r]|{xy}\ar@/^1pc/|{w_2}[r] &
W_{i_N}\ar@/^1pc/|{(xy)^{r_2}}[l]
\ar[l]|{(xy)^{-1+r_2}w_2^{c_2}w_3^{d_2}}}
\]
Note $(xy)^{-1+r_2}w_2^{c_2}w_3^{d_2}=(xy)^{-1+r_2}w_2$ and also
by Lemma~\ref{a2,q,r2,r3,r2} $-1+r_2=2r_3-2$ which is greater than zero since $\alpha_N>2$ forces $r_3>1$.  Now if  $\alpha_N\geq 3$, for every $t$ with $4\leq t\leq \alpha_N+1$ add an extra arrow %\footnote{so there are $\alpha_N-2$ extra arrows, since we include the endpoints} 
from $W_{1}$ to $R$ labelled
\[
(xy)^{r_t}w_2^{-1+c_t}w_3^{d_t}=(xy)^{\alpha_N+1-t}w_2^{t-3}w_3
\]
where the equality holds since %\footnote{eg $r_3=\alpha_N-2$, which is true} 
$r_t=\alpha_N+1-t$ by Lemma~\ref{r_butt_lemma} whilst $c_t=t-3$ and $d_t=1$ by Lemma~\ref{c_and_d_lemma}.  The extra maps %\footnote{Its curious: for the $\alpha_{N}=3$ case the maps back are $(xy)^3$, $(xy)^2w_2$ and $w_3$.  Without $(xy)^2w_2$ we still have all invariants at all vertices, but $(xy)^2w_2$ doesn't come for free.  However trying to prove that $W_3$ is 2-generated fails since the $(xy)^2w_2$ pops up at the end; see sheet for working} 
go where they should since $(xy)^{r_t}w_2^{c_t}w_3^{d_t}$ is an invariant for all $4\leq t\leq e=\alpha_N+1$ and we have a map $w_2$ from $R$ to $W_{i_N}$.  It is easy to see that at both vertices we have all invariants.  Let $f\in W_{i_N}$ and consider the invariant %\footnote{$t=\alpha_N+1$ extra arrow}
$w_2^{\alpha_N-2}w_3f$.   Now the two original arrows from $W_{i_N}$ to $R$ both have factor $xy$ and further all the extra arrows from $W_{i_N}$ to $R$ have factor $xy$, except $w_2^{\alpha_N-2}w_3$. Consequently viewing $w_2^{\alpha_N-2}w_3f$ as cycles at the vertex $W_{i_N}$ we can write
\[
w_2^{\alpha_N-2}w_3f=w_2^{\alpha_N-2}w_3A+(xy)B
\]
where $A$ is a sum of paths from $R$ to $W_{i_N}$ and $B$ is some polynomial.  But there are only two arrows from $R$ to $W_{i_N}$ (namely $xy$ and $w_2$) and so writing $A$ in terms of them
\[
w_2^{\alpha_N-2}w_3f=w_2^{\alpha_N-1}w_3A_1+(xy)(B+w_2^{\alpha_N-2}w_3A_2)
\]
for some polynomials $A_1,A_2$.  Thus $w_2^{\alpha_N-2}w_3$ must divide $B+w_2^{\alpha_N-2}w_3A_2$ and so after cancelling these terms we get
\[
f=w_2A_1+(xy)B_2
\]
for some polynomial $B_2$.  Hence by Lemma~\ref{polys_to_invariants} it follows that $W_1=W_{i_N}$ is generated by $xy$ and $w_2$, and so is special.
\end{proof}

We now search for the remaining rank one specials.

\begin{defin}
Define $W_t$ to be the CM $\C{}[x,y]^{\mathbb{D}_{n,q}}$-module containing $x^q+y^q$, and $W_{-}$
to be the CM module containing $x^q-y^q$.
\end{defin}
These are well defined since $x^q\pm y^q$ is a relative invariant
for the one-dimensional representations
\[
\begin{array}{c|c}
n-q \mbox{ odd} & n-q \mbox{ even}\\ \hline \begin{array}{ccl}
\psi_{2q}&\mapsto& \mp 1\\\tau & \mapsto &
\e_4^q\\\varphi_{2(n-q)}&\mapsto &\e_{2(n-q)}^q
\end{array} & \begin{array}{ccl}
\psi_{2q}&\mapsto& \mp 1 \\\tau\varphi_{4(n-q)}&\mapsto &
\e_4^q\e_{4(n-q)}^{q}
\end{array}
\end{array}
\]
Note also that $W_+$ and $W_-$ correspond to non-isomorphic representations, and that they are also distinct from the $W_t$ defined earlier.
\begin{lemma}\label{W_plus_minus_special}
For any $\mathbb{D}_{n,q}$ with any $\nu$, $W_{+}$ is generated by the two elements $(xy)^{n-q}(x^q-y^q), x^q+y^q$ whilst $W_{-}$ is generated by the two elements $(xy)^{n-q}(x^q+y^q), x^q-y^q$.  Hence both are special.
\end{lemma}
\begin{proof}
Let $f\in W_{+}$.  We first claim that $w_2=(x^q+y^q)(x^q+(-1)^{\nu}y^q)\in W_{i_{\nu+1}}$.  If $\nu=N-1$, then this follows since $w_2$ generates $W_{i_{\nu+1}}=W_{i_N}=W_1$ by Theorem~\ref{Wi's_are_special_nu_top}; if $0\leq \nu<N-1$ it follows by inspection of $\tt{D}_1$.  It follows that $(x^q+(-1)^\nu y^q)f\in W_{i_{\nu+1}}$.  But by combining Theorem~\ref{Wi's_are_special_most} and Theorem~\ref{Wi's_are_special_nu_top}, for any $\nu$ we know that $W_{i_{\nu+1}}$ is generated by $(xy)^{i_{\nu+1}}$ and $w_2^{\Delta_{\nu+1}}w_2^{\Gamma_{\nu+1}}=w_2=(x^q+y^q)(x^q+(-1)^\nu y^q)$ and so we may write
\[
(x^q+(-1)^\nu y^q)f=(xy)^{i_{\nu+1}}C_1+(x^q+y^q)(x^q+(-1)^\nu y^q)C_2
\]
for some invariant polynomials $C_1,C_2$.  This means $x^q+(-1)^\nu y^q$ must divide $C_1$ and so by inspection of the list of generators of the invariant ring (Lemma~\ref{invariants}) we may write%\footnote{minus sign correct as its $r_3$ and thus $w_3$}
\[
C_1=(xy)^{r_3}(x^q-y^q)(x^q+(-1)^{\nu}y^q)D_1+(x^q+ y^q)(x^q+(-1)^{\nu}y^q)D_2
\]
for some polynomials $D_1,D_2$,  thus
\[
f=(xy)^{r_3+i_{\nu+1}}(x^q-y^q)D_1+(x^q+y^q)(C_2+(xy)^{i_{\nu+1}}D_2).
\]
Now by Lemma~\ref{a2,q,r2,r3,r2}, $r_3+i_{\nu+1}=n-q$ and so
\[
f=(xy)^{n-q}(x^q-y^q)D_1+(x^q+y^q)(C_2+(xy)^{i_{\nu+1}}D_2),
\]
hence by Lemma~\ref{polys_to_invariants} it follows that $W_+$ is generated by $(xy)^{n-q}(x^q-y^q)$ and $x^q+y^q$, thus is special.  The argument for $W_-$ is symmetrical.
\end{proof}

Summarizing what we have proved:
\begin{thm}\label{specials_generators}
For any $\mathbb{D}_{n,q}$, the following CM modules are special and
further they are 2-generated as
$\C{}[x,y]^{\mathbb{D}_{n,q}}$-modules by the following elements:
\[
\begin{array}{l|cc|cc}
W_{+} & x^q+y^q &(xy)^{n-q}(x^q-y^q)&x^q+y^q &(xy)^{n-q}(x^q-y^q) \\
W_{-} & x^q-y^q & (xy)^{n-q}(x^q+y^q)& x^q-y^q & (xy)^{n-q}(x^q+y^q) \\
W_{i_{\nu+1}} &(xy)^{i_{\nu+1}} & w_2=w_2^{\Delta_{\nu+1}}w_3^{\Gamma_{\nu+1}}&(xy)^{i_{\nu+1}} & v_2=v_2^{\Delta_{\nu+1}}v_3^{\Gamma_{\nu+1}}\\
W_{i_{\nu+2}} &(xy)^{i_{\nu+2}} & w_2^{\Delta_{\nu+2}}w_3^{\Gamma_{\nu+2}}&(xy)^{i_{\nu+2}} & v_2^{\Delta_{\nu+2}}v_3^{\Gamma_{\nu+2}}\\
&\vdots&&&\\
W_{i_N} & (xy)^{i_N} & w_2^{\Delta_{N}}w_3^{\Gamma_{N}}& (xy)^{i_N} & v_2^{\Delta_{N}}v_3^{\Gamma_{N}}
\end{array}
\]
where the left column is one such choice of generators, and the right-hand column is another choice.
Further there are no other non-free rank one specials.
\end{thm}
\begin{proof}
Combine Lemma~\ref{W_plus_minus_special},
Theorem~\ref{Wi's_are_special_most} and
Theorem~\ref{Wi's_are_special_nu_top}.  Since they are distinct and we have $N+2-\nu$ of them, by Theorem~\ref{WunramMainResults}(3) these are them all. 
\end{proof}

We can go further and assign to each of the above specials the corresponding vertex in the minimal resolution.  The $\nu>0$ version of the following lemma can be found in \cite{Wemyss_reconstruct_D(ii)}.  
\begin{lemma}\label{specials2vertices}
Consider $\mathbb{D}_{n,q}$ with $\nu=0$.  Then the special CM modules above correspond to the dual graph of the minimal resolution in the following way
\[
\begin{array}{cc}
\begin{array}{c}
\xymatrix@C=20pt@R=15pt{ &\bullet\ar@{-}[d]^<{-2}&&&\\
\bullet\ar@{-}[r]_<{-2} & \bullet\ar@{-}[r]_<{-\alpha_{1}}
&\hdots\ar@{-}[r] &\bullet\ar@{-}[r]_<{-\alpha_{N-1}} & \bullet
\ar@{}[r]_<{-\alpha_N}&}
\end{array} &
\begin{array}{c}
\xymatrix@C=15pt@R=10pt{ &W_-\ar@{-}[d]&&&\\
W_+\ar@{-}[r] & W_{i_{1}}\ar@{-}[r]
&\hdots\ar@{-}[r] &W_{i_{N-1}}\ar@{-}[r] & W_{i_N}}
\end{array}
\end{array}
\]
\end{lemma}
\begin{proof}
The assumption $\nu=0$ translates into the condition $\alpha_1\geq 3$.  Here the fundamental cycle $Z_f$ is reduced from which, denoting the exceptional curves by $\{  E_i\}_{i\in I}$, we can easily calculate
\[
\begin{array}{cc}
\begin{array}{c}
\xymatrix@C=20pt@R=15pt{ &\bullet\ar@{-}[d]^<{1}&&&\\
\bullet\ar@{-}[r]_(0){1} & \bullet\ar@{-}[r]_(0){-3+\alpha_{1}}\ar@{-}[r]
&\bullet\ar@{-}[r]_(0){-2+\alpha_{2}}&\hdots\ar@{-}[r] &\bullet\ar@{-}[r]_(0){-2+\alpha_{N-1}} & \bullet
\ar@{}[r]_(0.2){-1+\alpha_N}&}
\end{array}&
\begin{array}{c}\xymatrix@C=20pt@R=15pt{ &\bullet\ar@{-}[d]^<{1}&&&\\
\bullet\ar@{-}[r]_<{1} & \bullet\ar@{-}[r]_<{-1}\ar@{-}[r]
&\bullet\ar@{-}[r]_<{0}&\hdots\ar@{-}[r] &\bullet\ar@{-}[r]_<{0} & \bullet
\ar@{}[r]_<{1}&}
\end{array}\\
(-Z_f\cdot E_i)_{i\in I}&((Z_K-Z_f)\cdot E_i)_{i\in I}
\end{array}
\]
%\begin{eqnarray*}
%(-Z_f\cdot E_i)_{i\in I}&=&\begin{array}{c}
%\xymatrix@C=20pt@R=15pt{ &\bullet\ar@{-}[d]^<{1}&&&\\
%\bullet\ar@{-}[r]_(0){1} & \bullet\ar@{-}[r]_(0){-3+\alpha_{1}}\ar@{-}[r]
%&\bullet\ar@{-}[r]_(0){-2+\alpha_{2}}&\hdots\ar@{-}[r] &\bullet\ar@{-}[r]_(0){-2+\alpha_{N-1}} & \bullet
%\ar@{}[r]_(0.2){-1+\alpha_N}&}
%\end{array}\\
%((Z_K-Z_f)\cdot E_i)_{i\in I}&=&\begin{array}{c}\xymatrix@C=20pt@R=15pt{ &\bullet\ar@{-}[d]^<{1}&&&\\
%\bullet\ar@{-}[r]_<{1} & \bullet\ar@{-}[r]_<{-1}\ar@{-}[r]
%&\bullet\ar@{-}[r]_<{0}&\hdots\ar@{-}[r] &\bullet\ar@{-}[r]_<{0} & \bullet
%\ar@{}[r]_<{1}&}
%\end{array}
%\end{eqnarray*}
where the canonical cycle $Z_K$ is the rational cycle defined by the condition $Z_K\cdot E_i=-K_{\w{X}}\cdot E_i$ for all $i\in I$.  Now denoting $R=\C{}[x,y]^{\mathbb{D}_{n,q}}$ if we consider the quiver of $\t{End}_{R}(R\oplus W_+\oplus W_-\oplus_{t=1}^{N}W_{i_t})$ we must be able to see the generators of the specials as compositions of irreducible maps out of $R$.  But by \cite[Corollary 3.1]{Wemyss_GL2} using the above intersection theory calculation it follows that we must see the generators of the specials using only compositions (containing no cycles) of the maps
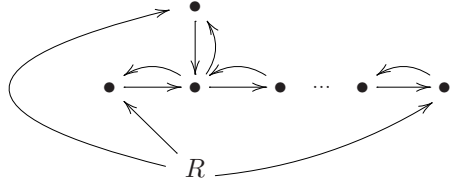
\[
\xymatrix@C=20pt@R=20pt{ &\bullet\ar[d]&&&\\
\bullet\ar[r]&
\bullet\ar@/_0.65pc/[u]\ar@/_0.65pc/[l]\ar[r]
&\bullet\ar@/_0.65pc/[l]\ar@{}[r]\ar@{}[r]|{\hdots} & \bullet\ar[r] & \bullet\ar@/_0.65pc/[l] \\
&R\ar@<1ex>@/^5.5pc/[-2,0]\ar[-1,-1]\ar@<-0.6ex>@/_0.55pc/[-1,3]&&&}
\]
Inspecting the list of generators of the specials, it is clear that $xy\in W_{i_N}$ cannot factor through any of the other specials, thus we must have this as an arrow in the quiver and so consequently $W_{i_N}$ must correspond to one of the vertices above to which $R$ connects.  The same analysis holds for $(+):=x^q+y^q\in W_+$ and $(-):=x^q-y^q\in W_-$ and so these too must correspond to vertices to which $R$ connects.  Now both $(+)^2$ and $(-)^2$ belong to $W_{i_1}$, and it is clear that they factor as 
$R\xrightarrow{(+)}W_+\xrightarrow{(+)}W_{i_1}$ and $R\xrightarrow{(-)}W_-\xrightarrow{(-)}W_{i_1}$ respectively.  By inspection of the generators of the specials, both the maps $W_+\xrightarrow{(+)}W_{i_1}$ and $W_-\xrightarrow{(-)}W_{i_1}$ are irreducible, hence $W_{i_1}$ must be a common neighbour to both $W_+$ and $W_-$.  This forces the positions of $W_+$, $W_-$ and $W_{i_1}$  in the dual graph as in the statement, and also forces the position of $W_{i_N}$ since it must occupy the final vertex which is connected to $R$.  Now the polynomial $(xy)^{i_1}$ factors as
\[
\xymatrix@C=65pt{ W_{i_{1}}&W_{i_{2}}\ar@/_0.65pc/[l]|{(xy)^{i_{1}-i_{2}}}\ar@{}[r]|{\hdots} & W_{i_{N-1}} & W_{i_N}\ar@/_0.65pc/[l]|{(xy)^{i_{N-1}-i_N}}&R,\ar@/_0.65pc/[l]|{xy} }
\]
forcing the remaining positions.
\end{proof}

\section{The reconstruction algebra}
In this section we define the reconstruction algebra $D_{n,q}$
with parameter $\nu=0$; when $\nu>0$ a corresponding definition
can be found in \cite{Wemyss_reconstruct_D(ii)}.   In fact we give two different presentations of this algebra, and prove both are isomorphic to the endomorphism ring of the sum of the special CM modules.

Consider, for $N\in\mathbb{N}$ with $N\geq 2$ and for positive integers $\alpha_1\geq 3$ and $\alpha_i\geq 2$ for all $2\leq i\leq N$, the labelled Dynkin diagram of type $D$:
\[
\xymatrix@C=20pt@R=15pt{ &\bullet\ar@{-}[d]^<{-2}&&&\\
\bullet\ar@{-}[r]_<{-2} & \bullet\ar@{-}[r]_<{-\alpha_{1}}
&\hdots\ar@{-}[r] &\bullet\ar@{-}[r]_<{-\alpha_{N-1}} & \bullet
\ar@{}[r]_<{-\alpha_N}&}
\]
We call the left-hand vertex the $+$ vertex, the top vertex
the $-$ vertex and the vertex corresponding to $\alpha_i$ the
$i^{th}$ vertex.  To the labelled Dynkin diagram above we add an extended vertex
(called $\star$) and `double-up' as follows:
\[%@C=50pt@R=50pt is best!
\xymatrix@C=40pt@R=40pt{ &\bullet\ar@<-2ex>@/_7.25pc/[2,0]|(0.4){\an{-}{0}}\ar[d]|{\cl{-}{1}\quad}&&&\\
\bullet\ar[r]|{\cl{+}{1}}\ar@/_0.65pc/[1,1]|{\an{+}{0}}&
\bullet\ar@/_0.65pc/[u]|{\quad\an{1}{-}}\ar@/_0.65pc/[l]|{\an{1}{+}}\ar[r]|{\cl{1}{2}}
&\bullet\ar@/_0.65pc/[l]|{\an{2}{1}}\ar@{}[r]|{\hdots} & \bullet\ar[r]|{\cl{N-1}{N}} & \bullet\ar@/_0.65pc/[l]|{\an{N}{N-1}}\ar@<0.25ex>[1,-3]|{\cl{N}{0}} \\
&\star\ar@<1ex>@/^6.25pc/[-2,0]|(0.5){\cl{0}{-}}\ar[-1,-1]|(0.55){\cl{0}{+}}\ar@<-0.6ex>@/_0.55pc/[-1,3]|(0.45){\an{0}{N}}&&&}
\]
Now if $\sum_{i=1}^{N}(\alpha_i-2)\geq 2$, we add extra arrows to the above picture in the following way:
\begin{itemize}
\item If $\alpha_1>3$, then add an extra $\alpha_1-3$ arrows from the
$1^{st}$ vertex to $\star$.
\item If $\alpha_i>2$ with $i\geq 2$, then add an extra $\alpha_i-2$ arrows from the $i^{th}$ vertex to $\star$.
\end{itemize}
Label the new arrows %\footnote{since counting the endpoints there are $1-2+\sum_{i=\nu+1}^{N}(\alpha_i-2)=-1+\sum_{i=\nu+1}^{N}(\alpha_i-2)=(\alpha_{\nu+1}-3)+\sum_{i=\nu+2}^{N}(\alpha_i-2)$ of them} 
(if they exist) by $k_2,\hdots,k_{\sum_{i=1}^{N}
(\alpha_i-2)}$ starting from the $1^{st}$ vertex and working to
the right.  Name this new quiver $Q$.
\begin{example}\label{D185quiver} 
\t{Consider $\mathbb{D}_{18,5}$ then $\frac{18}{5}=[4,3,2]$ and so
the quiver $Q$ is}
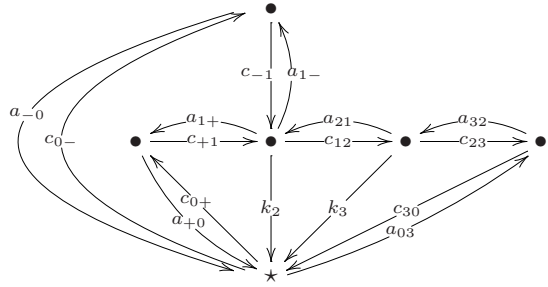
\[%@C=50pt@R=50pt is best!
\xymatrix@C=40pt@R=40pt{ &\bullet\ar@<-2ex>@/_7.25pc/[2,0]|(0.4){\an{-}{0}}\ar[d]|{\cl{-}{1}\quad}&&&\\
\bullet\ar[r]|{\cl{+}{1}}\ar@/_0.65pc/[1,1]|{\an{+}{0}}&
\bullet\ar@/_0.65pc/[u]|{\quad\an{1}{-}}\ar@/_0.65pc/[l]|{\an{1}{+}}\ar[d]|{k_2}\ar[r]|{\cl{1}{2}}
&\bullet\ar[1,-1]|{k_3}\ar@/_0.65pc/[l]|{\an{2}{1}}\ar[r]|{\cl{2}{3}} & \bullet\ar@/_0.65pc/[l]|{\an{3}{2}}\ar@<0.25ex>[1,-2]|{\cl{3}{0}} \\
&\star\ar@<1ex>@/^6.25pc/[-2,0]|(0.5){\cl{0}{-}}\ar[-1,-1]|(0.55){\cl{0}{+}}\ar@<-0.4ex>@/_0.45pc/[-1,2]|(0.45){\an{0}{3}}&&&}
\]
\end{example}
\begin{example}\label{D5211quiver}
\t{Consider $\mathbb{D}_{52,11}$ then $\frac{52}{11}=[5,4,3]$ and so
the quiver $Q$ is}
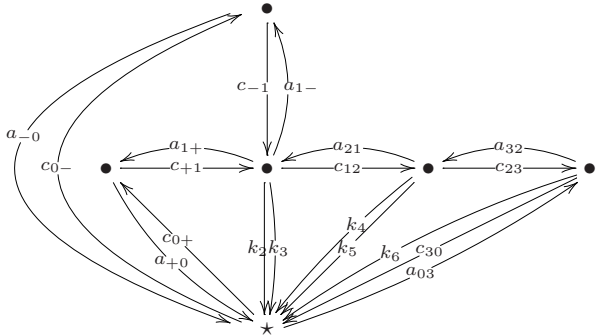
\[
\xymatrix@C=50pt@R=50pt{ &\bullet\ar@<-2ex>@/_7.25pc/[2,0]|(0.4){\an{-}{0}}\ar[d]|{\cl{-}{1}\quad}&&&\\
\bullet\ar[r]|{\cl{+}{1}}\ar@/_0.65pc/[1,1]|{\an{+}{0}}&
\bullet\ar@/_0.65pc/[u]|{\quad\an{1}{-}}\ar@/_0.65pc/[l]|{\an{1}{+}}\ar|{k_2\quad\,\,}@/^0.25pc/[d]|{\,\,k_3}\ar@<-0.25ex>[d]\ar[r]|{\cl{1}{2}}
&\bullet\ar@/_0.25pc/@<-0.25ex>[1,-1]|(0.4){k_4}\ar@<0ex>[1,-1]|{k_5}\ar@/_0.65pc/[l]|{\an{2}{1}}\ar[r]|{\cl{2}{3}} & \bullet\ar@<0.25ex>@/_0.5pc/[1,-2]|(0.6){k_6}\ar@/_0.65pc/[l]|{\an{3}{2}}\ar@<0.25ex>[1,-2]|{\cl{3}{0}} \\
&\star\ar@<1ex>@/^6.25pc/[-2,0]|(0.5){\cl{0}{-}}\ar[-1,-1]|(0.55){\cl{0}{+}}\ar@<-0.4ex>@/_0.45pc/[-1,2]|(0.45){\an{0}{3}}&&&}
\]
\end{example}
Now for every $\mathbb{D}_{n,q}$ with $\nu=0$ denote $k_1:=\an{+}{0}$ and $k_{1+\sum_{i=1}^{N}(\alpha_i-2)}:=\cl{N}{0}$.
\begin{defin}
For all %\footnote{extremes need the specially defined $k$} 
$1\leq r\leq 1+\sum_{i=1}^{N}(\alpha_i-2)$ define $\tt{B}_{r}$ to be the number (or $+$) of the vertex associated to the tail of the arrow $k_r$.
\end{defin}
Notice %\footnote{OK as the old $b_2$ recorded $\nu+1$ if $\alpha_{\nu+1}>3$ or else the next biggest $\alpha>2$.}\footnote{for the $1+\sum_{i=1}^{N}(\alpha_i-2)$ extreme need to use the special definitions of $k_{1+\sum_{i=1}^{N}(\alpha_i-2)}$ and $b_{1+\sum_{i=1}^{N}(\alpha_i-2)}=N$. } 
for all $2\leq r\leq 1+\sum_{i=1}^{N}(\alpha_i-2)$ it is true that $\tt{B}_r=b_r$ where $b_r$ is the $b$-series of $\frac{n}{q}$ defined in Section 2. However %\footnote{we also had defined $b_0$ in a special way before; here we don't use it}
$\tt{B}_1\neq b_1$ since $b_1=\nu+1=1$ whilst $\tt{B}_1=+$ by the definition of $k_1$.

Now define $u_+=1$ and further for $1\leq i\leq N$ denote%\footnote{top is $e-2$ so that $u_N=e-2=1+\sum_{p=1}^{N}(\alpha_p-2)$, hence top is $k_{e-2}:=\cl{0}{N}$.  This way the Step $N$ relations don't need a special ending}
\[
\begin{array}{c}
u_i:=\t{max}\{ j: 2\leq j\leq e-2\mbox{ with } b_{j}=i \}\\ 
v_i:=\t{min}\{ j: 2\leq j\leq e-2\mbox{ with } b_{j}=i \}
\end{array}
\]
if such things exist (i.e. vertex $i$ has an extra arrow leaving it).  Also define $W_1:=+$ and for every $2\leq i\leq N$ consider the set
\[
\c{S}_i=\{\mbox{vertex }j: 1\leq j< i\mbox{ and }j\mbox{ has an extra arrow leaving it} \}.
\]
For $2\leq i\leq N$ define
\[
W_i=\left\{ \begin{array}{cl} + & \mbox{if }\c{S}_i\mbox{ is empty}\\ \mbox{the maximal number in }\c{S}_i & \mbox{else}
\end{array}\right. 
\]
and so $W_i$ is defined for all $1\leq i\leq N$.  The idea behind it is that $W_i$ records the closest vertex to the left of vertex $i$ which has a $k$ leaving it; since we have defined $k_1:=\an{+}{0}$ this is always possible to find.  Now define, for all $1\leq i\leq N$, $V_{i}=u_{W_i}$.  Thus $V_i$ records the number of the largest $k$ to the left of the vertex $i$, where since $k_1:=\an{+}{0}$ and $u_+=1$ it  always exists.

We now define the reconstruction algebra as a presentation by generators and relations.  In what follows, we give two presentations --- one which we call the `symmetric presentation', and one which we call the `moduli presentation'.
  
\begin{defin}\label{reconD}
For $\frac{n}{q}=[\alpha_1,\hdots,\alpha_N]$ with $\nu=0$ (i.e.
$\alpha_1\geq 3$) define $D_{n,q}$, the reconstruction algebra of
type $D$, to be the path algebra of the quiver $Q$ defined above
subject to the relations
\[
\begin{array}{ccc}
1.&&\cl{0}{+}\cl{+}{1}-\cl{0}{-}\cl{-}{1}=4\AN{0}{1}\\
2.&&\cl{0}{+}\an{+}{0}=\cl{0}{-}\an{-}{0}\\
3.&&\an{-}{0}\cl{0}{-}=\cl{-}{1}\an{1}{-}\\
4.&&\an{1}{+}\cl{+}{1}=\an{1}{-}\cl{-}{1}
\end{array}
\]
together with the relations defined algorithmically as:%\footnote{the definition of $u_N$ is such that $k_{u_N}=\cl{0}{N}$}
\[
\begin{array}{ccl}
\mbox{Step 0:} & & \an{+}{0}\cl{0}{+}=\cl{+}{1}\an{1}{+}\\
&&\\
  \mbox{Step 1:} & \mbox{If }\alpha_1=3  & \cl{1}{2}\an{2}{1}=\an{1}{+}\cl{+}{1}\\
 & \mbox{If }\alpha_1> 3 & k_2\AN{0}{1}=\an{1}{+}\cl{+}{1}, \AN{0}{1}k_2=\cl{0}{+}\an{+}{0}\\
 &&k_t\CLr{0}{1}=k_{t+1}\AN{0}{1}, \CLr{0}{1}k_t=\AN{0}{1}k_{t+1}\,\,\forall\,\,2\leq t< u_1\\  & &
k_{u_1}\CLr{0}{1}=\cl{1}{2}\an{2}{1}\\  &&\vdots 
\end{array}
\]
\[
\begin{array}{ccl}
  \mbox{Step i:} & \mbox{If }\alpha_i=2  & \cl{i}{i+1}\an{i+1}{i}=\an{i}{i-1}\cl{i-1}{i}\\
 & \mbox{If }\alpha_i>2 & k_{v_i}\AN{0}{i}=\an{i}{i-1}\cl{i-1}{i}, \AN{0}{i}k_{v_i}=\CLr{0}{\tt{B}_{V_i}}k_{V_i}  \\
 & & k_{t}\CLr{0}{i}=k_{t+1}\AN{0}{i}, \CLr{0}{i}k_t=\AN{0}{i}k_{t+1}\,\,\forall\,\,v_i\leq t< u_i\\
 & & k_{u_i}\CLr{0}{i}=\cl{i}{i+1}\an{i+1}{i}\\  &&\vdots \\
  \mbox{Step N:} & \mbox{If }\alpha_N=2  & \cl{N}{0}\an{0}{N}=\an{N}{N-1}\cl{N-1}{N}, \CLr{0}{\tt{B}_{V_N}}k_{V_N}=\an{0}{N}\cl{N}{0} \\
 & \mbox{If }\alpha_N>2 & k_{v_N}\an{0}{N}=\an{N}{N-1}\cl{N-1}{N},\an{0}{N}k_{v_N}=\CLr{0}{\tt{B}_{V_N}}k_{V_N}  \\
 & & k_{t}\CLr{0}{N}=k_{t+1}\an{0}{N}, \CLr{0}{N}k_{t}=\an{0}{N}k_{t+1}\,\,\forall\,\,v_N\leq t<
 u_N\\
 %& & k_{u_N}\CLr{0}{N}=\cl{N}{0}\an{0}{N}, \an{0}{N}\cl{N}{0}=\CLr{0}{N}k_{u_N}
\end{array}
\]
where $\AN{0}{t}:=\an{0}{N}\hdots\an{t+1}{t}$ for every $1\leq t\leq N$.  The $\CLr{}{}$'s are defined as follows: we define $\CLr{0}{+}:=\cl{0}{+}$, and $\CLr{0}{t}:=\CLr{0}{1}\cl{1}{2}\hdots\cl{t-1}{t}$ for all $1\leq t\leq N$.  The only thing that remains to be defined is $\CLr{0}{1}$, and it is this which changes according to the presentation, namely 
\[
\CLr{0}{1}:=\left\{ \begin{array}{cl}\cl{0}{+}\cl{+}{1}&\mbox{in moduli presentation}\\ \frac{1}{2}(\cl{0}{+}\cl{+}{1}+\cl{0}{-}\cl{-}{1})&\mbox{in symmetric presentation}\end{array}\right. .
\]
\end{defin}

\begin{remark}
\t{We should explain why we give two presentations.  The symmetric case is pleasing since it treats the  two $(-2)$ horns equally, so that the algebra produced is independent of how we view the dual graph (see Lemma~\ref{duality}).  On the other hand the moduli presentation treats one of the $(-2)$ horns (namely $+$) to be `better' as relations go through that vertex and not the $-$ vertex.  The moduli presentation makes the explicit geometry easier to write down in Section 5, and is also satisfactory from the viewpoint of Remark~\ref{recon_A_part} below.  We show in Theorem~\ref{iso_of_rings} that the two presentations yield isomorphic algebras, but note that the explicit isomorphism is difficult to write down.  For the moment denote $D_{n,q}$ for the moduli presentation and $D_{n,q}^\prime$ for the symmetric presentation, as a priori they may be different.}
\end{remark}
\begin{remark}\label{recon_A_part}
\t{In the moduli presentation (i.e. $\CLr{0}{1}=\cl{0}{+}\cl{+}{0}$) the algorithmic relations are \emph{precisely} the same as those for the reconstruction algebra of type $A$ associated to the data
\[
\xymatrix@C=40pt{ \bullet\ar@{-}[r]^<{-2} &
\bullet\ar@{-}[r]^<{-(\alpha_{1}-1)} &\bullet\ar@{-}[r]^<{-\alpha_{2}} &\hdots\ar@{-}[r]
&\bullet\ar@{-}[r]^<{-\alpha_{N-1}} & \bullet \ar@{}[r]^<{-\alpha_N}&}
\]
Consequently the moduli presentation of the reconstruction algebra of type $D$ for $\nu=0$ is simply a reconstruction algebra of type $A$, together with an extra piece stuck on to compensate for the dihedral horns.}
\end{remark}
\begin{remark}
\t{Since $\CLr{0}{+}:=\cl{0}{+}$, the two presentations are exactly the same if and only if $\alpha_1=3$ and $\alpha_2=\hdots=\alpha_N=2$.  This corresponds to the `base case' of \cite[Lemma 3.3]{Wemyss_GL2}.  Note also that %\footnote{Note also the benefit of defining $k_1:=\an{+}{0}$ - the relations do not fork!} 
the use of $\tt{B}$ (instead of the $b$-series $b$) in the above definition is not a typo since $W_N=+$ is certainly possible (e.g. in the family $\mathbb{D}_{2s+1,s}$) in which case $V_N=u_+=1$; consequently $\tt{B}_{V_N}=\tt{B}_1=+$, which is different to $b_1=1$.  Thus using the $b$-series $b$  can take us to the wrong vertex.}
\end{remark}
\begin{remark}
\t{Double care must be taken to get the algorithmic relations here from the ones in \cite{Wemyss_reconstruct_A}.  First, \emph{loc. cit} labels the extra arrows $k$ in the other direction, and second it also begins labelling with a $k_1$ instead of $k_2$ (which we use here).  The fact that the direction of the labelling of the $k$'s has changed is awkward but it doesn't matter due to the duality for reconstruction algebras of type $A$ (see \cite[2.10]{Wemyss_reconstruct_A}).}
\end{remark}
\begin{example}\label{D185relations}
\t{For the group $\mathbb{D}_{18,5}$, the symmetric presentation of the reconstruction algebra is the quiver in Example~\ref{D185quiver} subject to the relations}
\[
\begin{array}{c}
\cl{0}{+}\cl{+}{1}-\cl{0}{-}\cl{-}{1}=4(\an{0}{3}\an{3}{2}\an{2}{1})\\
\cl{0}{+}\an{+}{0}=\cl{0}{-}\an{-}{0}\\
\an{-}{0}\cl{0}{-}=\cl{-}{1}\an{1}{-}\\
\an{1}{+}\cl{+}{1}=\an{1}{-}\cl{-}{1}\\
\an{+}{0}\cl{0}{+}=\cl{+}{1}\an{1}{+}
\end{array}
\]
\[
\begin{array}{c}
\begin{array}{rl}
\an{1}{+}\cl{+}{1}=k_2(\an{0}{3}\an{3}{2}\an{2}{1})& (\an{0}{3}\an{3}{2}\an{2}{1})k_2=\cl{0}{+}\an{+}{0}
\end{array}\\
k_2(\frac{1}{2}(\cl{0}{+}\cl{+}{1}+\cl{0}{-}\cl{-}{1}))=\cl{1}{2}\an{2}{1}\\
\begin{array}{rl}
\an{2}{1}\cl{1}{2}=k_3(\an{0}{3}\an{3}{2})& (\an{0}{3}\an{3}{2})k_3=(\frac{1}{2}(\cl{0}{+}\cl{+}{1}+\cl{0}{-}\cl{-}{1}))k_2
\end{array}\\
k_3(\frac{1}{2}(\cl{0}{+}\cl{+}{1}+\cl{0}{-}\cl{-}{1})\cl{1}{2})=\cl{2}{3}\an{3}{2}\\
\begin{array}{rl}
\an{3}{2}\cl{2}{3}=\cl{3}{0}(\an{0}{3})& (\an{0}{3})\cl{3}{0}=(\frac{1}{2}(\cl{0}{+}\cl{+}{1}+\cl{0}{-}\cl{-}{1})\cl{1}{2})k_3
\end{array}
\end{array}
\]
\end{example}
\begin{example}\label{D5211relations}
\t{For the group $\mathbb{D}_{52,11}$, the moduli presentation of the reconstruction algebra is the quiver in Example~\ref{D5211quiver} subject to the relations}
\[
\begin{array}{c}
\cl{0}{+}\cl{+}{1}-\cl{0}{-}\cl{-}{1}=4\an{0}{3}\an{3}{2}\an{2}{1}\\
\cl{0}{+}\an{+}{0}=\cl{0}{-}\an{-}{0}\\
\an{-}{0}\cl{0}{-}=\cl{-}{1}\an{1}{-}\\
\an{1}{+}\cl{+}{1}=\an{1}{-}\cl{-}{1}\\
\an{+}{0}\cl{0}{+}=\cl{+}{1}\an{1}{+}\\
\begin{array}{rl}
\an{1}{+}\cl{+}{1}=k_2(\an{0}{3}\an{3}{2}\an{2}{1})& (\an{0}{3}\an{3}{2}\an{2}{1})k_2=\cl{0}{+}\an{+}{0}\\
k_2(\cl{0}{+}\cl{+}{1})=k_3(\an{0}{3}\an{3}{2}\an{2}{1})& (\an{0}{3}\an{3}{2}\an{2}{1})k_3=(\cl{0}{+}\cl{+}{1})k_2
\end{array}\\
k_3(\cl{0}{+}\cl{+}{1})=\cl{1}{2}\an{2}{1}\\
\begin{array}{rl}
\an{2}{1}\cl{1}{2}=k_4(\an{0}{3}\an{3}{2})& (\an{0}{3}\an{3}{2})k_4=(\cl{0}{+}\cl{+}{1})k_3\\
k_4(\cl{0}{+}\cl{+}{1}\cl{1}{2})=k_5(\an{0}{3}\an{3}{2})& (\an{0}{3}\an{3}{2})k_5=(\cl{0}{+}\cl{+}{1}\cl{1}{2})k_4
\end{array}\\
k_5(\cl{0}{+}\cl{+}{1}\cl{1}{2})=\cl{2}{3}\an{3}{2}\\
\begin{array}{rl}
\an{3}{2}\cl{2}{3}=k_6(\an{0}{3})& (\an{0}{3})k_6=(\cl{0}{+}\cl{+}{1}\cl{1}{2})k_5\\
k_6(\cl{0}{+}\cl{+}{1}\cl{1}{2}\cl{2}{3})=\cl{3}{0}(\an{0}{3})& (\an{0}{3})\cl{3}{0}=(\cl{0}{+}\cl{+}{1}\cl{1}{2}\cl{2}{3})k_6
\end{array}
\end{array}
\]
\end{example}

The following is the main theorem of this paper.
\begin{thm}\label{iso_of_rings}
For a group $\mathbb{D}_{n,q}$ with parameter $\nu=0$ (i.e. reduced fundamental cycle), denote
$R=\C{}[x,y]^{\mathbb{D}_{n,q}}$ and let $T_{n,q}=R\oplus
W_{+}\oplus W_{-}\oplus_{t=1}^{N}W_{i_t}$ be the sum of the
special CM modules.  Then
\[
D_{n,q}\cong \t{End}_R(T_{n,q})\cong D_{n,q}^\prime.
\]
\end{thm}
\begin{proof}
We prove both statements at the same time, by making different choices for the generators of the specials.   Using the intersection theory in the proof of Lemma~\ref{specials2vertices} it follows immediately from \cite[Corollary 3.1]{Wemyss_GL2} that the quiver of the endomorphism ring of the specials is precisely that of the quiver $Q$ defined above.  We first find representatives for the known number of arrows:

As before denote $x^q+y^q=(+)$ and $x^q-y^q=(-)$.  We must reach the generators of the specials as paths out of $R$ (i.e. $\star$).  We know from the proof of Lemma~\ref{specials2vertices} we may choose $\cl{0}{+}=(+)$, $\cl{+}{1}=(+)$, $\cl{0}{-}=(-)$, $\cl{-}{1}=(-)$, $\an{0}{N}=xy$ and $\an{t+1}{t}=(xy)^{r_{l_t}}=(xy)^{i_t-i_{t+1}}$ (for all $1\leq t<N$) as representatives.  Now the generator $(xy)^{n-q}(+)$ of $W_-$ must be reached through $W_{i_1}$.  Since $(xy)^{n-q}(+)=(\an{0}{N}\hdots\an{2}{1})(xy)^{n-2q}(+)$ we may choose $\an{1}{-}=(xy)^{n-2q}(+)=(xy)^{r_3}(+)$.   By symmetry we may choose $\an{1}{+}=(xy)^{r_3}(-)$.

Now consider the generator $w_2^{\Delta_2}w_3^{\Gamma_2}$ of $W_{i_2}$.  We already have the generator $w_2=\cl{0}{+}\cl{+}{1}$ from $R$ to $W_{i_1}$, so it is clear that we may choose $\cl{1}{2}=w_2^{c_{l_2}}w_3^{d_{l_2}}$.   If we consider the generator $v_2$ of $W_{i_1}$ instead (which we have as $\frac{1}{2}(\cl{0}{+}\cl{+}{1}+\cl{0}{-}\cl{-}{1})$) and want to obtain the generator $v_2^{\Delta_2}v_3^{\Gamma_3}$ of $W_{i_2}$, instead choose $\cl{1}{2}=v_2^{c_{l_2}}v_3^{d_{l_2}}$. Continuing like this we can choose $\cl{t}{t+1}=w_2^{c_{l_t}}w_3^{d_{l_t}}$ for all $1\leq t< N$.

We now claim that we may choose $\cl{N}{0}=w_2^{c_{l_N}}w_3^{d_{l_N}}=w_2^{c_{e-1}}w_3^{d_{e-1}}$. To see this, first note that $w_2^{c_{e-1}}w_3^{d_{e-1}}$ doesn't factor through $\an{N}{N-1}$ (the only possible map to the non-trivial specials).  Second, note that $w_2^{c_{e-1}}w_3^{d_{e-1}}$ can't factor as some map from $W_{i_N}$ to $R$ multiplied by a non-scalar invariant %\footnote{this covers $W_{i_N}\rightarrow R\rightarrow R$ and $W_{i_N}\rightarrow W_{i_N}\rightarrow R$} 
else the invariant generator $\an{0}{N}\cl{N}{0}=(xy)^{r_{e-1}}w_2^{c_{e-1}}w_3^{d_{e-1}}$ factors into two non-scalar invariants, contradicting the embedding dimension.   A similar argument shows that we may choose $\an{+}{0}=(xy)^{r_3}(-)$ and $\an{-}{0}=(xy)^{r_3}(+)$. %\footnote{or use non-explicit pre-projective relation at -2 curve}  

Hence we have labelled all arrows in $Q$ by polynomials, except $k_2,\hdots,k_{e-3}$ (if the $k$ arrows exist).  How to do this is obvious by the quiver $\tt{D}_1$ in Section 3:  if the $k$'s exist we can choose $k_t$ labelled by $k_t=(xy)^{r_{t+2}}w_2^{c_{t+1}}w_3^{d_{t+1}}$.  The argument that these choices don't factor through other specials via maps of strictly positive lower degrees (and so can be chosen as representatives) is similar to the above --- for example if $\alpha_1\geq 4$ consider $(xy)^{r_{4}}w_2^{c_{3}}w_3^{d_{3}}$.  First, it does not factor through maps we have already chosen, since if it does we may write $(xy)^{r_{4}}w_2^{c_{3}}w_3^{d_{3}}=\an{1}{-}f+\an{1}{+}g+\cl{1}{2}h=(xy)^{r_3}F+w_2^{c_{l_1}}w_3^{d_{l_1}}h$.  But by looking at $xy$ powers, we know $(xy)^{r_4}$ divides $h$ and so after cancelling factors we may write $w_2^{c_3}w_3^{d_3}=(xy)^{r_3-r_4}F+w_2^{c_{l_1}}w_3^{d_{l_1}}h_1$.  After cancelling $w_2^{c_3}w_3^{d_3}$ (which $F$ must be divisible by) $1=w_2^{c_{l_1}-c_3}w_3^{d_{l_1}-d_3}h_1+(xy)^{r_3-r_4}F^\prime$ which is impossible since the right-hand side cannot have degree zero terms.  Second, $(xy)^{r_{4}}w_2^{c_{3}}w_3^{d_{3}}$ does not factor as a map $W_{i_1}\rightarrow R$ followed by a non-scalar invariant since again this would contradict the embedding dimension. Hence we may choose $k_2=(xy)^{r_{4}}w_2^{c_{3}}w_3^{d_{3}}$ in this case. %\footnote{the next $k$ is harder, since the vector space is 'the maps which don't factor through another special through maps of strictly positive smaller degree'; the linear span of $k_2$ isn't just $\{ \lambda k_2 \}$ but also has all other factorizations added in} 
Continue like this: if $\alpha_1\geq 5$ we want to choose $k_3=(xy)^{r_{5}}w_2^{c_{4}}w_3^{d_{4}}$.  If it factors through maps we have already chosen then we may write $k_3=\an{1}{-}f+\an{1}{+}g+k_2j+\cl{1}{2}h=(xy)^{r_4}F+\cl{1}{2}h$, so just  repeating the argument above shows that this is impossible.% \footnote{Contradict the embedding dimension again too.  Note that this argument works all the way to and including vertex $N$ since the choice of $\cl{N}{0}$ has been made already.}\footnote{In general to get the big $(xy)$ factor out to compose by needs not just the $a$ arrow out, but also }

By the above, we have justified that we may choose the following as representatives of all the irreducible maps between the special CM modules%\footnote{Note this gives simples view too should we need it}
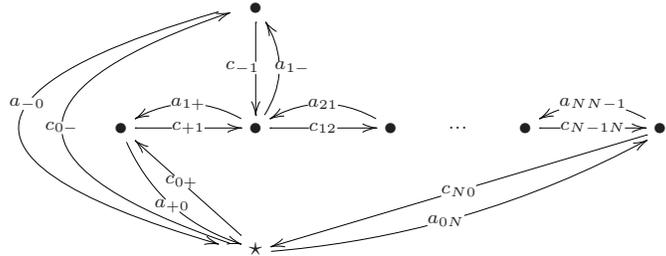
\[
\begin{array}{c}
\begin{array}{c}
\xymatrix@C=40pt@R=35pt{ &\bullet\ar@<-2ex>@/_6.75pc/[2,0]|(0.4){\an{-}{0}}\ar[d]|{\cl{-}{1}\quad}&&&\\
\bullet\ar[r]|{\cl{+}{1}}\ar@/_0.75pc/[1,1]|{\an{+}{0}}&
\bullet\ar@/_0.75pc/[u]|{\quad\an{1}{-}}\ar@/_0.75pc/[l]|{\an{1}{+}}\ar[r]|{\cl{1}{2}}
&\bullet\ar@/_0.75pc/[l]|{\an{2}{1}}\ar@{}[r]|{\hdots} & \bullet\ar[r]|{\cl{N-1}{N}} & \bullet\ar@/_0.75pc/[l]|{\an{N}{N-1}}\ar@<0.25ex>[1,-3]|{\cl{N}{0}} \\
&\star\ar@<1ex>@/^5.75pc/[-2,0]|(0.5){\cl{0}{-}}\ar[-1,-1]|(0.55){\cl{0}{+}}\ar@<-0.4ex>@/_0.75pc/[-1,3]|(0.45){\an{0}{N}}&&&}
\end{array}
\\ \\ {\small{
\begin{array}{l}
\begin{array}{lcl}
\cl{0}{+}=\cl{+}{1}=x^q+y^q &&\an{1}{+}=\an{+}{0}=(xy)^{r_3}(x^q-y^q)\\
\cl{0}{-}=\cl{-}{1}=x^q-y^q && \an{1}{-}=\an{-}{0}=(xy)^{r_3}(x^q+y^q)\\
\cl{N}{0}=w_2^{c_{l_N}}w_3^{d_{l_N}}&&
\an{0}{N}=(xy)^{r_{l_N}}=(xy)^{i_N-i_{N+1}}=xy
\end{array}\\
\\
\left.\begin{array}{c}
\cl{t}{t+1}=w_2^{c_{l_t}}w_3^{d_{l_t}}\\
\qquad\an{t+1}{t}=(xy)^{r_{l_t}}=(xy)^{i_t-i_{t+1}}
\end{array}\right\}\mbox{for all }1\leq t< N
\end{array}}}
\end{array}
\]
and further in the above picture we also have, if $\sum_{i=1}^{N}(\alpha_i-2)\geq 2$, for every $2\leq t\leq \sum_{i=1}^{N}(\alpha_i-2)$ the extra arrows $k_t$ labelled by $k_t=(xy)^{r_{t+2}}w_2^{c_{t+1}}w_3^{d_{t+1}}$. %\footnote{Notice (since $r_e=0$) the pattern that $\cl{N}{0}=k_{\sum(\alpha_i-2)+1}=k_{e-2}$, if $k_{\sum(\alpha_i-2)+1}$ was defined as in the line above.  The top defined value of $r$ series is $r_e$, where $e=v\sum(\alpha_i-2)+3=(\sum(\alpha_i-2)+1)+2$, thus $\cl{N}{0}$ has $r_e=0$.}
The symmetric presentation choices are identical, except everywhere we replace $w_2$ by $v_2$ and $w_3$ by $v_3$.

Denote the relations in Definition~\ref{reconD} by $\c{S}^\prime$.  For the relations part of the proof, below we are really working in the completed case (so we can use \cite[Corollary 3.1]{Wemyss_GL2} and \cite[3.4]{BIRS}) and we prove that the completion of the endomorphism ring of the specials is given as the completion of $\C{}Q$ (denoted $\C{}\hat{Q}$) modulo the closure of the ideal $\langle \c{S}^\prime\rangle$ (denoted $\overline{\langle \c{S}^\prime\rangle}$).  The non-completed version of the theorem then follows by simply taking the associated graded ring of both sides of the isomorphism.

Now denote the kernel %\footnote{the kernel is a closed ideal by \cite[3.1]{BIRS}}
of the surjection $\C{}\hat{Q}\rightarrow \t{End}_{\C{}[[x,y]]^G}(T_{n,q}):=\Lambda$ by $I$, denote the radical of $\C{}\hat{Q}$ by $J$ and further for $t\in\{ \star,+,-,1,\hdots,N \}$ denote by $S_t$ the simple corresponding to the vertex $t$ of $Q$.  In Lemma~\ref{LIrelations} below we show that the elements of $\c{S}^\prime$ are linearly independent in $I/ (IJ+JI)$.  Thus we may extend $\c{S}^\prime$ to a basis $\c{S}$ of $I/ (IJ+JI)$.  Now by \cite[3.4(a)]{BIRS} $I=\overline{\langle \c{S}\rangle}$ and further by \cite[3.4(b)]{BIRS}
\[
\t{dim}_{\C{}}\t{Ext}^2_\Lambda(S_{a},S_{b})=\# (e_{a}\C{}\hat{Q}e_{b})\cap \c{S}
\]
for all $a,b\in\{ \star,+,-,1,\hdots,N \}$.  But on the other hand using the intersection theory in the proof of Lemma~\ref{specials2vertices} it follows immediately from \cite[Corollary 3.1]{Wemyss_GL2} that %\footnote{$\star$ OK: see p3 that $-Z_f\cdot Z_f-1=(-Z_f\cdot Z_f+1)-2=e-2=\sum_{i=1}^{N}(\alpha_1-2)+1$ where $e$ is embedding.}
\[
\begin{array}{lcl}
\t{dim}_{\C{}}\t{Ext}^2_\Lambda(S_+,S_+)&=&1\\
\t{dim}_{\C{}}\t{Ext}^2_\Lambda(S_-,S_-)&=&1\\
\t{dim}_{\C{}}\t{Ext}^2_\Lambda(S_i,S_i)&=&\alpha_i-1\\
\t{dim}_{\C{}}\t{Ext}^2_\Lambda(S_\star,S_\star)&=&1+\sum_{p=1}^{N}(\alpha_p-2)\\
\t{dim}_{\C{}}\t{Ext}^2_\Lambda(S_\star,S_1)&=&1
\end{array}
\]
for all $1\leq i\leq N$, and further all other $\t{Ext}^2$'s between the simples are zero.   By inspection of both the above information and the relations $\c{S}^\prime$ we notice that 
\[
\t{dim}_{\C{}}\t{Ext}^2_\Lambda(S_{a},S_{b})=\# (e_{a}\C{}\hat{Q}e_{b})\cap \c{S}^\prime
\]
for all $a,b\in\{ \star,+,-,1,\hdots,N \}$.  Hence 
\[
\# (e_{a}\C{}\hat{Q}e_{b})\cap \c{S}=\# (e_{a}\C{}\hat{Q}e_{b})\cap \c{S}^\prime
\]
for all $a,b\in\{ \star,+,-,1,\hdots,N \}$, proving that the number of elements in $\c{S}$ and $\c{S}^\prime$ are the same.  Hence $\c{S}^\prime=\c{S}$ and so $I=\overline{\langle \c{S}^\prime\rangle}$, as required.
\end{proof}

\begin{lemma}\label{LIrelations}
With notation from the above proof, the members of $\c{S}^\prime$ are linearly independent in $I/(IJ+JI)$.
\end{lemma}
\begin{proof}
First, it is easy to verify that all the members in $\c{S}^\prime$ are satisfied by the chosen representatives of the arrows and so belong to $I$.  To see this, note that the first four relations follow immediately by inspection (independent of presentation), as does the Step 0 relation.  For the moduli presentation the remaining algorithmic relations are simply the pattern between the invariants and cycles in $\tt{D}_1$ from Proposition~\ref{loops_in_D1}, which is the same as the pattern in type $A$.  The symmetric pattern is just a small modification of this, namely the pattern between the invariants and cycles in the $\tt{D}_2$ version of Proposition~\ref{loops_in_D1}. %\footnote{For $\nu=0$ it is true that $w_3=v_3$. Now the first $k_2$ is given by $k_2=(xy)^{r_4}w_2^{c_3}w_3^{d_3}=(xy)^{r_4}w_3=(xy)^4v_3$ and so the $\an{1}{+}\cl{+}{1}=\hdots=\AN{0}{1}k_2$ holds for both presentations.  It is the other $k_2$ relation (involving the $c$'s) which changes depending on whether $v$ or $w$ choice is made.}  
Thus in either presentation $\c{S}^\prime\subseteq I$.  

We grade the arrows in $Q$ by the degree of the polynomial representing that arrow above.  In what follows we say that a word $w$ in the path algebra $\C{}\hat{Q}$ satisfies condition (A) if
\[
\begin{array}{cl}
\t{(i)}&\t{It does not contain some proper subword which is a cycle.}\\
\t{(ii)}&\t{It does not contain some proper subword which is a path from $\star$ to $1$.}
\end{array}
\]
As stated in the proof of the above theorem, we know from the intersection theory and \cite[Corollary 3.1]{Wemyss_GL2} that the ideal $I$ is generated by one relation from $\star$ to $1$, whereas all other generators are cycles.  Consequently if a word $w$ satisfies (A) then $w\notin IJ+JI$.  It is also clear that $\cl{0}{+}\cl{+}{1}-\cl{0}{-}\cl{-}{1}-4\AN{0}{1}\notin IJ+JI$. %\footnote{by the $J$ part, we need something longer than a path $\star$ to $1$.}

Now since all members of $\c{S}^\prime$ are either cycles at some vertex or paths from $\star$ to $1$,   to prove that the members of $\c{S}^\prime$ are linearly independent in $I/(IJ+JI)$ we just need to show that 
\begin{itemize}
\item[1.] the elements of $\c{S}^\prime$ which are paths from $\star$ to $1$ are linearly independent in $e_\star(I/(IJ+JI))e_1$.
\item[2.] for all $t\in\{ \star,+,-,1,\hdots,N \}$, the elements of $\c{S}^\prime$ that are cycles at $t$ are linearly independent in $e_t(I/(IJ+JI))e_t$.
\end{itemize}

The first condition is easy, since the only relation in $\c{S}^\prime$ from $\star$ to $1$ is $\cl{0}{+}\cl{+}{1}-\cl{0}{-}\cl{-}{1}=4\AN{0}{1}$, and we have already noted that it does not belong to $IJ+JI$, thus it is non-zero and so linearly independent in $e_\star(I/(IJ+JI))e_1$.  Note also that there are only three paths of  minimal grade from $\star$ to $1$, namely $\{\cl{0}{+}\cl{+}{1},\cl{0}{-}\cl{-}{1},\AN{0}{1}\}$, and by inspection of the polynomials they represent we do not have any other relation from $\star$ to $1$ of this grade.

For the second condition, we must check $t$ case by case: \\
\emph{Case $t=+$.}  Here the only relation in $\c{S}^\prime$ from $+$ to $+$ is $\an{+}{0}\cl{0}{+}=\cl{+}{1}\an{1}{+}$ (the Step 0 relation) so it is linearly independent provided it is non-zero, i.e.  $\an{+}{0}\cl{0}{+}-\cl{+}{1}\an{1}{+}\notin IJ+JI$.  But since $\an{+}{0}\cl{0}{+}$ satisfies condition (A) we know that $\an{+}{0}\cl{0}{+}\notin IJ+JI$.  Thus if $\an{+}{0}\cl{0}{+}-\cl{+}{1}\an{1}{+}\in IJ+JI$ we may write
\[
\an{+}{0}\cl{0}{+}=\cl{+}{1}\an{1}{+}+u
\]
in the free algebra $\C{}\hat{Q}$ for some $u\in IJ+JI$.  But the term $\an{+}{0}\cl{0}{+}$ does not appear in the right-hand side, a contradiction.\\
\emph{Case $t=-$.} The only relation in $\c{S}^\prime$ from $-$ to $-$ is $\an{-}{0}\cl{0}{-}=\cl{-}{1}\an{1}{-}$ (the third relation).  A symmetrical argument to the above shows that it is linearly independent.\\
\emph{Case $t=1$.} If $\alpha_1=3$ then the only members of $\c{S}^\prime$ from $1$ to $1$ are $\cl{1}{2}\an{2}{1}=\an{1}{+}\cl{+}{1}$ (the Step 1 relation) and $\an{1}{+}\cl{+}{1}=\an{1}{-}\cl{-}{1}$ (the fourth relation).  Suppose that 
\[
\lambda_1(\cl{1}{2}\an{2}{1}-\an{1}{+}\cl{+}{1})+\lambda_2(\an{1}{+}\cl{+}{1}-\an{1}{-}\cl{-}{1})=0
\]
in $e_1(I/(IJ+JI))e_1$.  Then we may write
\[
\lambda_1\cl{1}{2}\an{2}{1}+\lambda_2\an{1}{+}\cl{+}{1}=\lambda_1\an{1}{+}\cl{+}{1}+\lambda_2\an{1}{-}\cl{-}{1}+u
\]
in the free algebra $\C{}\hat{Q}$ for some $u\in IJ+JI$.  But $\cl{1}{2}\an{2}{1}$ satisfies (A) so doesn't appear anywhere in the right-hand side, forcing $\lambda_1=0$.  But now since $\lambda_1=0$ and further $\an{1}{+}\cl{+}{1}$ satisfies (A), it cannot appear on the right hand side, forcing $\lambda_2=0$.    This proves the assertion when $\alpha_1=3$ and hence we may assume that $\alpha_1>3$, in which case the only relations in $\c{S}^\prime$ from $1$ to $1$ are
\[
\begin{array}{l}
\an{1}{-}\cl{-}{1}=\an{1}{+}\cl{+}{1}\\
 k_2\AN{0}{1}=\an{1}{+}\cl{+}{1}\\
 k_t\CLr{0}{1}=k_{t+1}\AN{0}{1}\,\,\forall\,\,2\leq t< u_1\\  
k_{u_1}\CLr{0}{1}=\cl{1}{2}\an{2}{1}.
\end{array}
\]
(i.e. the fourth relation and some of the Step $1$ relations). Now suppose that
\begin{multline*}
\lambda_1(\an{1}{-}\cl{-}{1}-\an{1}{+}\cl{+}{1})+\lambda_2(\an{1}{+}\cl{+}{1}-k_2\AN{0}{1})+\sum_{p=2}^{u_1-1}\lambda_{p+1}(k_{p+1}\AN{0}{1}-k_p\CLr{0}{1})\\+\lambda_{\alpha_1-1}(\cl{1}{2}\an{2}{1}-k_{u_1}\CLr{0}{1})=0
\end{multline*}
in $e_1(I/(IJ+JI))e_1$ then
\begin{multline*}
\lambda_1\an{1}{-}\cl{-}{1}+\lambda_2\an{1}{+}\cl{+}{1}+\sum_{p=2}^{u_1-1}\lambda_{p+1}k_{p+1}\AN{0}{1}+\lambda_{\alpha_1-1}\cl{1}{2}\an{2}{1}\\=\lambda_1\an{1}{+}\cl{+}{1}+\lambda_2k_2\AN{0}{1}+\sum_{p=2}^{u_1-1}\lambda_{p+1}k_p\CLr{0}{1}+\lambda_{\alpha_1-1}k_{u_1}\CLr{0}{1}+u
\end{multline*}
in the free algebra $\C{}\hat{Q}$ for some $u\in IJ+JI$. Now $\an{1}{-}\cl{-}{1}$ satisfies (A) and so cannot appear on the right-hand side; consequently $\lambda_1=0$.  This combined with the fact that $\an{1}{+}\cl{+}{1}$ satisfies (A) implies that $\lambda_2=0$. Similarly $\cl{1}{2}\an{2}{1}$ satisfies (A) and so can't appear on the right hand side, so $\lambda_{\alpha_1-1}=0$.  This leaves
\begin{eqnarray}
\sum_{p=2}^{u_1-1}\lambda_{p+1}k_{p+1}\AN{0}{1}=\sum_{p=2}^{u_1-1}\lambda_{p+1}k_p\CLr{0}{1}+u \label{B}
\end{eqnarray}
in the free algebra $\C{}\hat{Q}$ and so
\[
\lambda_{u_1}k_{u_1}\AN{0}{1}\equiv \t{terms starting with k of strictly smaller index}
\]
mod $IJ+JI$.  But since $k_{u_1}\AN{0}{1}$ does not have any subwords which are cycles, the only way we can change it mod $IJ+JI$ is to bracket as $k_{u_1}(\AN{0}{1})$ and use the relation in $I$ from $\star$ to $1$.  Doing this we get $k_{u_1}\AN{0}{1}\equiv \frac{1}{4}k_{u_1}(\cl{0}{+}\cl{+}{1}-\cl{0}{-}\cl{-}{1})$.  This still does not start with a $k$ of strictly lower index, so we must again use relations in $IJ+JI$ to change the terms.  But again the words contain no subwords which are cycles, which means we must either change $\cl{0}{+}\cl{+}{1}$ or $\cl{0}{-}\cl{-}{1}$.  But there is only one relation between $\cl{0}{+}\cl{+}{1}$, $\cl{0}{-}\cl{-}{1}$and $\AN{0}{1}$ so no matter what we do we arrive back as
\[
k_{u_1}\AN{0}{1}\equiv \frac{1}{4}k_{u_1}(\cl{0}{+}\cl{+}{1}-\cl{0}{-}\cl{-}{1})\equiv k_{u_1}\AN{0}{1}
\]
mod $IJ+JI$.  Thus mod $IJ+JI$ it is impossible to transform $k_{u_1}\AN{0}{1}$ into an expression involving $k$ terms with strictly lower index, and so we must have $\lambda_{u_1}=0$.  With this in mind we may re-arrange (\ref{B}) to get
\[
\lambda_{u_1-1}k_{u_1-1}\AN{0}{1}\equiv \t{terms starting with k of strictly smaller index}
\]
mod $IJ+JI$ and so repeating the above argument gives $\lambda_{u_1-1}=0$.  Continuing in this fashion we deduce all $\lambda$'s are zero, as required.\\
\emph{Case $t$ for $1<t\leq N$.}  %\footnote{Here the case $t=N$ is no different from others, except we need to be careful with $u_N$.  Note also that although there is a difference in the algorithmic relations between Step $i$ and Step $N$, this is only at vertex $\star$.}  
If $\alpha_t=2$ then the only relation in $\c{S}^\prime$ from $t$ to $t$ is $\an{t}{t-1}\cl{t-1}{t}=\cl{t}{t+1}\an{t+1}{t}$ (the Step $t$ relation).  This is linearly independent in $e_t(I/(IJ+JI))e_t$ using the same argument as in the case $t=+$.  Hence we may assume $\alpha_t>2$ in which case the only relations in $\c{S}^\prime$ from $t$ to $t$ are
\[
\begin{array}{l}
k_{v_t}\AN{0}{t}=\an{t}{t-1}\cl{t-1}{t}\\
k_{p}\CLr{0}{t}=k_{p+1}\AN{0}{t}\,\,\forall\,\,v_t\leq p< E\\
k_{E}\CLr{0}{t}=\cl{t}{t+1}\an{t+1}{t}
\end{array}
\]
(i.e. some of the Step $t$ relations) where $E:=\left\{\begin{array}{cl}u_t& \mbox{if }t<N\\ u_N-1&\mbox{if }t=N \end{array}\right. $ and also we mean $\cl{N}{N+1}=\cl{N}{0}$ and $\an{N+1}{N}=\an{0}{N}$ if at any place the subscripts become too large.  Now if%\footnote{Note that the sum is empty if $\alpha_t=3$, and also notice that $u_t-v_t+1=\alpha-2$ so have $\alpha_t-1$ relations}
\[
\lambda_1(k_{v_t}\AN{0}{t}-\an{t}{t-1}\cl{t-1}{t})+\sum_{p=0}^{E-v_t-1}\lambda_{p+2}(k_{p+v_t}\CLr{0}{t}-k_{p+v_t+1}\AN{0}{t})+\lambda_{\alpha_t-1}(k_{E}\CLr{0}{t}-\cl{t}{t+1}\an{t+1}{t})=0
\]
in $e_t(I/(IJ+JI))e_t$ then%\footnote{Notice we get round the problem in Case $t=1$ by just appealing to condition (A) now.}
\begin{multline*}
\lambda_1\an{t}{t-1}\cl{t-1}{t}+\sum_{p=0}^{E-v_t-1}\lambda_{p+2}k_{p+v_t+1}\AN{0}{t}+\lambda_{\alpha_t-1}\cl{t}{t+1}\an{t+1}{t}=\\ \lambda_1k_{v_t}\AN{0}{t}+\sum_{p=0}^{E-v_t-1}\lambda_{p+2}k_{p+v_t}\CLr{0}{t}+\lambda_{\alpha_t-1}k_{E}\CLr{0}{t}+u
\end{multline*}
in the free algebra $\C{}\hat{Q}$  for some $u\in IJ+JI$.  But now all terms on the left hand side satisfy (A) and so none of them can appear on the right hand side, forcing $\lambda_1=\hdots=\lambda_{\alpha_t-1}=0$ as required.\\
\emph{Case $t=\star$}.  The only relations left in $\c{S}^\prime$ are those from $\star$ to $\star$, which is the second relation together with the remaining relations from the algorithm.  These are precisely %\footnote{Note that there are $1+(\alpha_1-3)+(\alpha_2-2)+\hdots+(\alpha_N-2+1)=1+\sum_{p=1}^{N}(\alpha_p-2)$ of them!}
\[
\begin{array}{cc}
&\cl{0}{+}\an{+}{0}=\cl{0}{-}\an{-}{0}\\
  \mbox{If }\alpha_1>3& \left\{ \begin{array}{c} \AN{0}{1}k_2=\cl{0}{+}\an{+}{0}\\
  \CLr{0}{1}k_t=\AN{0}{1}k_{t+1}\,\,\forall\,\,2\leq t< u_1\end{array}\right.\\
 &\vdots\\
 \mbox{If }\alpha_i>2 &  \left\{\begin{array}{c}\AN{0}{i}k_{v_i}=\CLr{0}{\tt{B}_{V_i}}k_{V_i}  \\
  \CLr{0}{i}k_t=\AN{0}{i}k_{t+1}\,\,\forall\,\,v_i\leq t< u_i\end{array}\right.\\
 &\vdots  \\
  \mbox{If }\alpha_N=2 &  \an{0}{N}\cl{N}{0}=\CLr{0}{\tt{B}_{V_N}}k_{V_N} \\
  \mbox{If }\alpha_N>2 &\left\{\begin{array}{c} \an{0}{N}k_{v_N}=\CLr{0}{\tt{B}_{V_N}}k_{V_N}  \\
 \an{0}{N}k_{t+1}=\CLr{0}{N}k_{t}\,\,\forall\,\,v_N\leq t<u_N\\
 %\an{0}{N}\cl{N}{0}=\CLr{0}{N}k_{u_N}
 \end{array}\right.
\end{array}
\]
Hence if $\alpha_1=3$, $\alpha_2=\hdots=\alpha_N=2$ then there are no extra arrows $k$ and so the only relations at $\star$ are $\cl{0}{+}\an{+}{0}=\cl{0}{-}\an{-}{0}$ and $\an{0}{N}\cl{N}{0}=\cl{0}{+}\an{+}{0}$.  These are linearly independent by using the same argument as in the case $t=1$ with $\alpha_1=3$.  Thus we may assume that some extra $k$ arrows exist, in which case the proof is similar to the case $t=1$ with $\alpha_1>3$.%\footnote{Line up with $\lambda_1$,etc.  Inspecting $\cl{0}{-}\an{-}{0}$ gives $\lambda_1=0$ which \emph{then} forces $\lambda_2=0$.  Always get $\lambda_{last}=0$ by inspecting $\an{0}{N}\cl{N}{0}$.  Then go through the remainder - loads will have property (A) (eg $\CL{0}{1}k_t=k_{t+1}\AN{0}{3}$, the one on the right has property (A) so use it) and so get rid off those; the remainder that are left have the form $k_t\CL{0}{1}k_t=\AN{0}{1}k_{t+1}$ (note the 1 on both sides) which we get rid off by the same argument as $t=1$ with $\alpha_1>3$.}
\end{proof}

The input to define the reconstruction algebra is a certain labelled Dynkin diagram of type $D$, where the two `horns' are both $(-2)$ curves.  These are geometrically indistinguishable in the sense that their positions in the dual graph can be swapped and the input for the reconstruction algebra does not change.  Thus the reconstruction algebra should be invariant under this change of labels, which leads us to the following.

\begin{lemma}\label{duality}
For $\lambda\in\C{*}$ denote $D_{n,q}^{\prime}(\lambda)$ to be the algebra obtained from the symmetric presentation $D_{n,q}^{\prime}$ by replacing the number $4$ in the first relation by $\lambda$.  Similarly define $D_{n,q}(\lambda)$.  Then for all $\lambda\in\C{*}$
\[
D_{n,q}^{\prime}(\lambda)\cong D_{n,q}^{\prime}\cong D_{n,q}\cong D_{n,q}(\lambda).
\]
In particular the algebra obtained from $D_{n,q}^{\prime}$ by everywhere swapping $+$ and $-$ is isomorphic to $D_{n,q}^{\prime}$.
\end{lemma}
\begin{proof}
We prove that $D_{n,q}^{\prime}(\lambda)\cong D_{n,q}^{\prime}$; the proof that $D_{n,q}(\lambda)\cong D_{n,q}$ is identical.  All we must do is change the choices of the labels of the arrows made for $D_{n,q}^{\prime}$ in the proof of Theorem~\ref{iso_of_rings}, such that the relation $(1^\prime.)\,\,\cl{0}{+}\cl{+}{1}-\cl{0}{-}\cl{-}{1}=\lambda\AN{0}{1}$ together with all the other original relations (except relation 1) hold.  Since the only changes we shall make to the choices of arrows in the proof of Theorem~\ref{iso_of_rings} is by multiplying them by a non-zero scalar (for an arrow $p$, denote by $\kappa_p$ this non-zero scalar), the argument of Theorem~\ref{iso_of_rings} again goes through to show that $D_{n,q}^{\prime}(\lambda)\cong \t{End}_R(T_{n,q})$ and so in particular $D_{n,q}^{\prime}(\lambda)\cong D_{n,q}^\prime$.  

Choose $\kappa_{\cl{0}{+}}=\kappa_{\an{+}{0}}=\kappa_{\cl{0}{-}}=\kappa_{\an{-}{0}}=\kappa_{\cl{+}{1}}=\kappa_{\an{1}{+}}=\kappa_{\cl{-}{1}}=\kappa_{\an{1}{-}}=1$, $\kappa_{\an{0}{N}}=\frac{\lambda}{4}$ and $\kappa_{\an{N}{N-1}}=\hdots=\kappa_{\an{2}{1}}=1$.  Then certainly relations $1^\prime$, $2$, $3$ and $4$ hold, as does the Step 0 relation.  What remains is to choose $\kappa_{\cl{1}{2}},\hdots,\kappa_{\cl{N-1}{N}},\kappa_{\cl{N}{0}}$ (and the $\kappa_{k}$ for the $k$ arrows if they exist) and to verify the remaining relations.   

Consider the $b$-series of $\frac{n}{q}$.  If $b_2>1$ then choosing $\kappa_{\cl{1}{2}}=\hdots=\kappa_{\cl{b_2-1}{b_2}}=1$ it is clear that the step 1 to the step $b_2-1$ relations hold.  Thus in all cases we can consider the $b_2$ relations, knowing the previous step relations hold.

Now if $\alpha_{b_2}=2$ %\footnote{can only happen if $\alpha_1=3$ with $\alpha_2=\hdots=\alpha_N=2$!}
then $b_2=N$ and there are no $k$ arrows, so by choosing $\kappa_{\cl{N}{0}}=\frac{4}{\lambda}$ it is clear that the Step $N$ relations hold so we are done.  Hence we can assume that $\alpha_{b_2}>2$, in which case choosing $\kappa_{k_2}=\frac{4}{\lambda}, \hdots,\kappa_{k_{u_{b_2}}}=(\frac{4}{\lambda})^{u_{b_2}-1}, \kappa_{\cl{b_2}{b_2+1}}=(\frac{4}{\lambda})^{u_{b_2}-1}$ the step $b_2$ relations hold. %\footnote{index on the $c$ is mod $N$; i.e. the $c$ arrow may be $\cl{N}{0}$}

Thus inductively we consider the Step $t$ relations with $\kappa_{\cl{1}{2}},\hdots,\kappa_{\cl{t-1}{t}}$ and $k_2,\hdots,k_{V_t}$ %\footnote{Note that these both are non empty since $\alpha_{b_2}>2$}
already chosen such that all relations up to and including Step $t-1$ are satisfied.  If $\alpha_{b_2+1}=2$ choose $\kappa_{\cl{t+1}{t}}=\kappa_{\cl{t-1}{t}}$ then the Step $t$ relations hold. %\footnote{even if $t=N$.}
 Else choose $\kappa_{v_t}=\kappa_{\cl{t-1}{t}}(\frac{4}{\lambda}),\hdots,\kappa_{u_t}=\kappa_{\cl{t-1}{t}}(\frac{4}{\lambda})^{(u_t-v_t)+1},\kappa_{\cl{t}{t+1}}=\kappa_{\cl{t-1}{t}}(\frac{4}{\lambda})^{(u_t-v_t)+1}$ then the Step $t$ relations hold. %\footnote{again index on the $c$ is mod $N$; i.e. the $c$ arrow may be $\cl{N}{0}$}
This concludes the induction step, hence the result follows.

The final statement is now immediate since by inspection of the symmetric presentation the algebra obtained by everywhere swapping $+$ and $-$ is just $D_{n,q}^{\prime}(-4)$.
\end{proof}

\section{The moduli space of representations}
In this section we use quiver GIT on the reconstruction algebra $D_{n,q}$ to obtain the minimal resolution of $\C{2}/\mathbb{D}_{n,q}$ and so obtain the slightly stronger statement that the special representations not only give the dual graph, but also the whole space.  The surprising thing here is that although the spaces involved are not toric, the proofs follow almost immediately from the toric case.

As in type $A$, fix for the rest of this paper the moduli presentation of the reconstruction algebra, the dimension vector $\alpha=(1,1,\hdots,1)$ and the generic stability condition $\theta=(-(N+2),1,\hdots,1)$.   Notice that a representation $M$ of dimension vector $\alpha$ is $\theta$-stable if and only if it is generated from vertex $\star$, i.e. for every vertex in the representation $M$ there is a non-zero path in $M$ from $\star$ to that vertex.  

Throughout this section we use the moduli presentation convention that $\CL{0}{1}=\cl{0}{+}\cl{+}{1}$ whereas $\CL{0}{t}=\CL{0}{1}\cl{1}{2}\hdots\cl{t-1}{t}$ for any $2\leq t\leq N$.  For $\mathbb{D}_{n,q}$ with $Z_f$ reduced, we claim (and prove in Lemma~\ref{opencover}) that the moduli space is covered by the following $N+3$ open sets:
\[
\begin{array}{rcc}
U_0&\CL{0}{N}\neq 0, \cl{0}{-}\neq 0 & (\an{1}{-},\cl{N}{0},\an{0}{N})\\
&  \vdots &  \\
\begin{array}{c}U_t\\ _{1\leq t\leq N-1}
\end{array}&\CL{0}{N-t}\neq 0, \cl{0}{-}\neq 0, \AN{0}{N-t+1}\neq 0 &  (\an{1}{-},\cl{N-t}{N-t+1}\an{N-t+1}{N-t})\\
&  \vdots &  \\
U_{N}&\cl{0}{+}\neq 0, \cl{0}{-}\neq 0, \AN{0}{1}\neq 0 & (\an{1}{-},\an{1}{+},\cl{+}{1})\\
U_{+}&\cl{0}{+}\neq 0, \AN{0}{1}\neq 0, \an{1}{-}\neq 0 & (\cl{0}{-},\an{1}{+},\an{-}{0})\\
U_{-}&\cl{0}{-}\neq 0, \AN{0}{1}\neq 0, \an{1}{+}\neq 0 & (\cl{0}{+},\an{1}{-},\an{+}{0})
\end{array}
\]
where in the above we have stated for reference the result of Lemma~\ref{moduli_is_smooth_hypers}, which gives the position of where (if we change basis so that the specified non-zero arrows are actually the identity) the co-ordinates can be read off the quiver.  In fact in Lemma~\ref{moduli_is_smooth_hypers} we prove a little more; we show, also reading off the quiver, that these open sets are given abstractly by the following smooth hypersurfaces in $\C{3}$:
\[
\begin{array}{rcc}
{}_{(0\leq t\leq N-1)}& U_t  & a(1-4b^{\sum_{i=1}^{N}\eta_{t}^{(i)}}c^{\sum_{i=1}^{N}\theta_{t}^{(i)}})=b^{\eta_{t}^+}c^{\theta_t^+}\\
&U_{N}&  a(c-4)=bc\\
&U_{+}&  b(a^2c+4)=ac\\
&U_{-}&  b(a^2c-4)=ac
\end{array}
\]
%where the $\eta$ and $\theta$ are specified combinatorics introduced in the proof of Lemma~\ref{moduli_is_smooth_hypers}.  Also by taking an arbitrary stable module $M$ and simply changing basis appropriately the above open cover translates into ratios of polynomials as:
%\[
%\begin{array}{rl}
%U_0
%%&\C{}[\frac{\an{1}{-}\CL{0}{1}}{\cl{0}{-}},\CL{0}{N}\cl{N}{0},\frac{\an{0}{N}}{\CL{0}{N}}]
%& \C{}[\frac{(xy)^{r_3}(+)^3}{(-)},(+)^{2+2\Delta_{N+1}+\Gamma_{N+1}}(-)^{\Gamma_{N+1}},\frac{xy}{(+)^{2+2\Delta_{N}+\Gamma_{N}}(-)^{\Gamma_{N}}}]\\
%& \vdots   \\
%\begin{array}{c}U_t\\ _{1\leq t\leq N-1}
%\end{array}
%%&  \C{}[\frac{\an{1}{-}\CL{0}{1}}{\cl{0}{-}},\frac{\CL{0}{N-t+1}}{\AN{0}{N-%t+1}},\frac{\AN{0}{N-t}}{\CL{0}{N-t}}]
%&\C{}[\frac{(xy)^{r_3}(+)^3}{(-)},\frac{(+)^{2+2\Delta_{N+1-t}+\Gamma_{N+1-t}}(-)^{\Gamma_{N+1-t}}}{(xy)^{i_{N+1-t}}},\frac{(xy)^{i_{N-t}}}{(+)^{2+2\Delta_{N-t}+\Gamma_{N-t}}(-)^{\Gamma_{N-t}}}]\\
%&  \vdots   \\
%U_{N}
%%&\C{}[\frac{\an{1}{-}\AN{0}{1}}{\cl{0}{-}},\frac{\an{1}{+}\AN{0}{1}}{\cl{0}{+}},\frac%{\cl{0}{+}\cl{+}{1}}{\AN{0}{1}}]
%&\C{}[\frac{(xy)^{n-q}(+)}{(-)},\frac{(xy)^{n-q}(-)}{(+)},\frac{(+)^2}{(xy)^q}]\\
%U_{+}
%%&\C{}[\frac{\cl{0}{-}}{\AN{0}{1}\an{1}{-}},\frac{\an{1}{+}\AN{0}{1}}{\cl{0}{+}},\an{
%%-}{0}\AN{0}{1}\an{1}{-}]
%&\C{}[\frac{(-)}{(xy)^{n-q}(+)},\frac{(xy)^{n-q}(-)}{(+)},(xy)^{r_2}(+)^2]\\
%U_{-}
%%&\C{}[\frac{\cl{0}{+}}{\AN{0}{1}\an{1}{+}},\frac{\an{1}{-}\AN{0}{1}}{\cl{0}{-}},\an{%+}{0}\AN{0}{1}\an{1}{+}]
%&\C{}[\frac{(+)}{(xy)^{n-q}(-)},\frac{(xy)^{n-q}(+)}{(-)},(xy)^{r_2}(-)^2]
%\end{array}
%\]
%where $(+)=x^q+y^q$ and $(-)=x^q-y^q$.  

Note that there is a choice of coordinates in $U_N$ since for the third coordinate $c$ we could instead choose $d=\frac{\cl{0}{-}\cl{-}{1}}{\AN{0}{1}}=\frac{(-)^2}{(xy)^{q}}$ since they differ by 4.  Picking this alternative co-ordinate changes the defining equation to $ad=b(4-d)$.  Although trivial, it makes the gluing of the affine pieces slightly easier: we shall see in Theorem~\ref{minres_is_moduli} that the gluing data between the open pieces is precisely
\[
U_t\ni(a,b,c)\leftrightarrow(a,c^{-1},c^{\alpha_{N-t}}b)\in U_{t+1}
\]
for $0\leq t\leq N-2$, whereas for $t=N-1$ the gluing data is
\[
U_{N-1}\ni (a,b,c)\leftrightarrow (ca,c^{\alpha_1-1}b,c^{-1})\in U_N.
\]
The choice of coordinate in $U_N$ gives the final two glues %\footnote{$U_N=\C{}[\frac{(xy)^{s+1}(+)}{(-)},\frac{(xy)^{s+1}(-)}{(+)},\frac{(-)^2}{(xy)^s}]\ni (a,b,d)$, $U_N=\C{}[\frac{(xy)^{s+1}(+)}{(-)},\frac{(xy)^{s+1}(-)}{(+)},\frac{(+)^2}{(xy)^s}] \ni (a,b,c)$}
\[
\begin{array}{c}
U_N\ni (a,b,d) \leftrightarrow (a^{-1},b,a^2d)\in U_{+} \\
U_N \ni (a,b,c)\leftrightarrow (b^{-1},a,b^2c)\in U_{-}
\end{array}
\]
and note that these two do not change from example to example. 

We now proceed to prove these statements.  To prove that the open sets mentioned above do indeed form an open cover of the moduli space it is convenient to denote $\CLm{0}{1}=\cl{0}{-}\cl{-}{1}$,  to denote $\CLm{0}{t}=\CLm{0}{1}\cl{1}{2}\hdots\cl{t-1}{t}$  for all $2\leq t\leq N$, and further also to define the following open sets:
\[
\begin{array}{rc}
V_0:&\CLm{0}{N}\neq 0, \cl{0}{+}\neq 0 \\
&  \vdots   \\
\begin{array}{c}V_t\\ _{1\leq t\leq N-1}
\end{array}:&\CLm{0}{N-t}\neq 0, \cl{0}{+}\neq 0, \AN{0}{N-t+1}\neq 0
\end{array}
\]
\begin{lemma}\label{Vs_OK_by_Us}
For every $0\leq k\leq N-1$, the open set $V_k$ is contained
inside the union of $U_k$ and $U_N$.
\end{lemma}
\begin{proof}
Suppose $M$ belongs to $V_k$, then necessarily
\[
0\neq \cl{0}{-}\cl{-}{1}=\cl{0}{+}\cl{+}{1}-4\AN{0}{1}.
\]
If $\cl{+}{1}=0$ then $\AN{0}{1}\neq 0$ and so $M$ is in $U_N$. If
$\cl{+}{1}\neq 0$ then $M$ is in $U_k$.
\end{proof}
\begin{lemma}\label{opencover}
The open sets $U_0,\hdots,U_N,U_+,U_-$ completely cover the moduli
space.
\end{lemma}
\begin{proof}
Suppose $M$ is a stable module; we must show that $M$ belongs to
one of the open sets in the statement.  Note first that if
$\cl{0}{+}=\cl{0}{-}=0$ then the relation
$\cl{0}{+}\cl{+}{1}-\cl{0}{-}\cl{-}{1}=4\AN{0}{1}$
forces $\AN{0}{1}=0$, which is impossible since $M$ must be
generated as a module from vertex $\star$.  Hence we can assume that
either $\cl{0}{+}\neq 0$ or $\cl{0}{-}\neq 0$.\\
\emph{Case 1:} Both $\cl{0}{+}\neq 0$ and $\cl{0}{-}\neq 0$. Now
if $\an{0}{N}=0$ then the only way a nonzero path can reach vertex
$N$ is if either $\cl{+}{1}\CL{1}{N}\neq 0$ or
$\cl{-}{1}\CL{1}{N}\neq 0$.  In the first case $M$ belongs to
$U_0$, whereas in the second case $M$ belongs to $V_0$ thus by
Lemma~\ref{Vs_OK_by_Us} either $U_0$ or $U_N$.  Hence we may
assume $\an{0}{N}\neq 0$.  Now if $\an{N}{N-1}=0$ then by a
similar argument $M$ is either in $U_1$ or $V_1$, hence by
Lemma~\ref{Vs_OK_by_Us} either $U_1$ or $U_N$. Thus we can assume
that also $\an{N}{N-1}\neq 0$ and so $\AN{0}{N-1}\neq 0$.
Continuing in this manner either $M$ is in $U_{N-1}$ or we can
assume that $\AN{0}{1}\neq 0$, in which case $M$ is in $U_N$.\\
\emph{Case 2:} $\cl{0}{+}\neq 0$ but $\cl{0}{-}=0$.  Then certainly $\an{1}{-}\neq 0$ since we have to reach vertex $-$ with a non-zero path.  Now if $\AN{0}{1}=0$ then by the non-monomial relation it follows that $\cl{0}{+}\cl{+}{1}=0$ and so we cannot reach vertex $1$ by a non-zero path.  This cannot happen, hence $\AN{0}{1}\neq 0$ and so we are in $U_+$. \\
\emph{Case 3:} $\cl{0}{-}\neq 0$ but $\cl{0}{+}=0$.  In this case we are in $U_-$ by the symmetric argument to Case 2.
\end{proof}

\begin{lemma}\label{moduli_is_smooth_hypers}
Each open set $U_0,\hdots,U_N,U_+,U_-$ is a smooth
hypersurface in $\C{3}$.  More precisely the equation of these
open sets as a hypersurface in $\C{3}$ with co-ordinates $a,b,c$
are given as follows:\\
\[
\begin{array}{rcc}
{}_{(0\leq t\leq N-1)}& U_t  & a(1-4b^{\sum_{i=1}^{N}\eta_{t}^{(i)}}c^{\sum_{i=1}^{N}\theta_{t}^{(i)}})=b^{\eta_t^+}c^{\theta_t^+}\\
&U_{N}&  a(c-4)=bc\\
&U_{+}&  b(a^2c+4)=ac\\
&U_{-}&  b(a^2c-4)=ac
\end{array}
\]
\end{lemma}
\begin{proof}
(i) In $U_0$ change basis so that all the specified non-zero arrows equal the identity.  By Remark~\ref{recon_A_part} the calculation of \cite[4.2]{Wemyss_reconstruct_A} shows every arrow (except for the moment $\cl{0}{-}=1$, $\an{-}{0}$, $\cl{-}{1}$, $\an{1}{-}$) is determined by a monomial in $\cl{N}{0}$ and $\an{0}{N}$, and the algorithmic relations play no further role.  Define
\[
\begin{array}{cll}
&\eta_{0}^{+}=\t{the power of }\,\,\cl{0}{N}\,\,\t{in}\,\,\an{1}{+}   & \theta_0^+=\t{the power of }\,\,\an{N}{0}\,\,\t{in}\,\,\an{1}{+} \\
(1\leq i\leq N)&\eta_0^{(i)}=\t{the power of}\,\,\cl{0}{N}\,\,\t{in}\,\,\an{i+1}{i} & \theta_0^{(i)}=\t{the power of}\,\,\an{N}{0}\,\,\t{in}\,\,\an{i+1}{i}\\
\end{array}
\]
where by $\an{N+1}{N}$ we mean $\an{0}{N}$.  From this it is clear that $\AN{0}{1}=\cl{0}{N}^{\sum_{i=1}^{N}\eta_{0}^{(i)}}\an{N}{0}^{\sum_{i=1}^{N}\theta_{0}^{(i)}}$.
Now we are left with the variables $\cl{N}{0}$, $\an{0}{N}$, $\an{-}{0}$, $\an{1}{-}$ and $\cl{-}{1}$ subject to the four relations
\[
\begin{array}{rcl}
\an{+}{0}&=&\an{-}{0}\\
\an{-}{0}&=&\cl{-}{1}\an{1}{-}\\
\an{1}{+}&=&\an{1}{-}\cl{-}{1}\\
1-\cl{-}{1}&=&4\AN{0}{1}
\end{array}
\]
and so really there are only three variables $\cl{N}{0}$, $\an{0}{N}$ and  $\an{1}{-}$ subject to the one relation
\[
\an{1}{-}(1-4\AN{0}{1})=\an{+}{0}.
\]
Hence it suffices to show how to put $\an{+}{0}$ and $\AN{0}{1}$ in terms of $\cl{N}{0}$ and $\an{0}{N}$.  But by the above %\footnote{and using the preprojective relation to see $\an{1}{+}=\an{+}{0}$} 
this becomes
\[
\an{1}{-}(1-4\cl{N}{0}^{\sum_{i=1}^{N}\eta_{0}^{(i)}}\an{0}{N}^{\sum_{i=1}^{N}\theta_{0}^{(i)}})=\cl{N}{0}^{\eta_0^+}\an{0}{N}^{\theta_0^+}.
\]
(ii) The proof of the case $U_t$ for $1\leq t\leq N-1$ is identical to the above --- after setting the specified non-zero elements to be the identity we are down to the three variables $\cl{N-t}{N-t+1}$, $\an{N-t+1}{N-t}$ and  $\an{1}{-}$ with only one relation
\[
\an{1}{-}(1-4\AN{0}{1})=\an{+}{0}
\]
so again it suffices to show how to put $\an{+}{0}$ and $\AN{0}{1}$ in terms of $\cl{N-t}{N-t+1}$ and $\an{N-t+1}{N-t}$.  Define %\footnote{Note most of these will be zero since $\an{0}{N}=1$ etc}
\[
\begin{array}{cll}
&\eta_{t}^{+}=\t{the power of}\,\,\cl{0}{N}\,\,\t{in}\,\,\an{1}{+}   & \theta_t^+=\t{the power of}\,\,\an{N}{0}\,\,\t{in}\,\,\an{1}{+} \\
(1\leq i\leq N)&\eta_t^{(i)}=\t{the power of}\,\,\cl{0}{N}\,\,\t{in}\,\,\an{i+1}{i} & \theta_t^{(i)}=\t{the power of}\,\,\an{N}{0}\,\,\t{in}\,\,\an{i+1}{i}\\
\end{array}
\]
then it is clear that the equation is simply
\[
\an{1}{-}(1-4\cl{N-t}{N-t+1}^{\sum_{i=1}^{N}\eta_{t}^{(i)}}\an{N-t+1}{N-t}^{\sum_{i=1}^{N}\theta_{t}^{(i)}})=\cl{N-t}{N-t+1}^{\eta_t^+}\an{N-t+1}{N-t}^{\theta_t^+}.
\]
(iii) For $U_N$  set the specified non-zero arrows to be 1 then by Remark~\ref{recon_A_part} and the calculation \cite[4.2]{Wemyss_reconstruct_A} every arrow (except for the moment $\cl{0}{-}=1$, $\an{-}{0}$, $\cl{-}{1}$, $\an{1}{-}$) is determined by a monomial in $\cl{+}{1}$ and $\an{1}{+}$, and the algorithmic relations play no further role.  Thus we are down to $\cl{+}{1}$, $\an{1}{+}$, $\an{-}{0}$, $\an{1}{-}$ and $\cl{-}{1}$ subject to the four relations
\[
\begin{array}{rcl}
\an{+}{0}&=&\an{-}{0}\\
\an{-}{0}&=&\cl{-}{1}\an{1}{-}\\
\an{1}{+}\cl{+}{1}&=&\an{1}{-}\cl{-}{1}\\
\cl{+}{1}-\cl{-}{1}&=&4.
\end{array}
\]
But these reduce to the three variables $\an{1}{-}$, $\an{1}{+}$ and $\cl{+}{1}$ subject to the relation
\[
\an{1}{+}\cl{+}{1}=\an{1}{-}(\cl{+}{1}-4).
\]
Note that since $\cl{+}{1}-\cl{-}{1}=4$ we actually have a choice of coordinate between $\cl{+}{1}$ and $\cl{-}{1}$; making the other choice changes the equation in the obvious way.\\
(iv) For $U_+$: after setting the specified non-zero arrows to be 1, by Remark~\ref{recon_A_part} and the calculation \cite[4.2]{Wemyss_reconstruct_A} every arrow (except for the moment $\cl{0}{-}$, $\an{-}{0}$, $\cl{-}{1}$, $\an{1}{-}=1$) is determined by a monomial in $\cl{+}{1}$ and $\an{1}{+}$ and the algorithmic relations play no further role.  Thus we are down to $\cl{+}{1}$, $\an{1}{+}$, $\cl{0}{-}$, $\an{-}{0}$ and $\cl{-}{1}$ subject to the four relations
\[
\begin{array}{rcl}
\an{+}{0}&=&\cl{0}{-}\an{-}{0}\\
\an{-}{0}\cl{0}{-}&=&\cl{-}{1}\\
\an{1}{+}\cl{+}{1}&=&\cl{-}{1}\\
\cl{+}{1}-\cl{0}{-}\cl{-}{1}&=&4
\end{array}
\]
But now these reduce to the three variables $\cl{0}{-}$, $\an{1}{+}$ and $\an{-}{0}$ subject to the relation
\[
\an{-}{0}\cl{0}{-}=\an{1}{+}(4+\an{-}{0}\cl{0}{-}^2)
\]
(v) $U_-$ follows from $U_+$ by a symmetrical argument.
\end{proof}
%\begin{remark}
%\t{The above proof also justifies the position of the coordinates of the open sets as stated in the beginning of this section.  Labelling the arrows in the reconstruction algebra with their corresponding polynomials, for any open set if we simply change bases to set the specified non-zero arrows to be 1 we also obtain easily the ratios of polynomials as stated previously.}
%\end{remark}
\begin{thm}\label{minres_is_moduli}
Consider $\mathbb{D}_{n,q}$ with $\nu=0$.  With the given dimension vector and stability as above, the moduli of representations of the corresponding reconstruction algebra is precisely the minimal resolution of $\C{2}/\mathbb{D}_{n,q}$.
\end{thm}
\begin{proof}
The open cover is smooth and irreducible, thus we can restrict our attention to the exceptional locus.

First, it's easy to see that the only gluing data needed are the glues mentioned in the introduction to this section --- for example if a stable module $M$ belongs to $U_0$ with coordinates (obtained \emph{after} setting the specified non-zero arrows to be the identity) just $(\an{1}{-},\cl{N}{0},\an{0}{N})$ then clearly for it to belong to $U_1$ requires $\an{0}{N}\neq 0$.  Now by just looking at the conditions of the non-zero arrows determining the other open sets, it's clear that if $U_0$ glues to any other of the open sets it necessarily also has to satisfy $\an{0}{N}\neq 0$,  thus necessarily this glue is through $U_1$.  Continuing in this fashion we see that the $N+2$ stated glues are all that are necessary.

Now it is also easy to see that the open pieces glue (in coordinates) in precisely the way mentioned in the introduction to this section --- this follows directly from type $A$ and in fact we can see this from the proof of the last lemma: to see why, for $0\leq t\leq N-2$
\[
U_t\ni(a,b,c)\leftrightarrow(a,c^{-1},c^{\alpha_{N-t}}b)\in U_{t+1}
\]
is true, notice that given $(a,b,c)\in U_t$ then by the proof of the above lemma the scalar in the position of $\an{N-t}{N-t-1}$ is $bc^{\alpha_{N-t}-1}$, so changing basis at vertex $N-t$ by dividing every arrow into vertex $t$ by $c$ whilst  multiplying every arrow out by $c$ yields the result.  The remaining three glues are also done by inspection.

Thus by inspecting our open cover and gluing data we see precisely the dual graph of the minimal resolution.  In fact, one must check that in each open affine, the $\mathbb{P}^1$'s from the gluing data form precisely the locus which is contracted, but this is an easy calculation for which we suppress the details.

To check that the self-intersection numbers are correct, we use adjunction: for example for the curve $C$ in the glue between $U_N$ (co-ordinates $a,b,c$ subject to $f=(a(c-4)-bc)$) and $U_-$ (co-ordinates $\tt{a},\tt{b},\tt{c}$ subject to $\tt{f}=(\tt{b}(\tt{a}^2\tt{c}-4)-\tt{a}\tt{c}$)) given by 
\[
U_N\ni (a,b,c)\leftrightarrow (b^{-1},a,b^2c)=(\tt{a},\tt{b},\tt{c})\in U_{-}
\]
we have
\[
\frac{d\tt{a}\wedge d\tt{c}}{\partial\tt{f}/\partial\tt{b}}=\frac{d\tt{a}\wedge d\tt{c}}{\tt{a}^2\tt{c}-4}=\frac{d(b^{-1})\wedge d(b^2c)}{c-4}=\frac{db\wedge dc}{c-4}=\frac{db\wedge dc}{\partial f/\partial a}
\]
and so $K_{\w{X}}\cdot C=0$, thus by adjunction $C^2=-2$.  Continuing in this fashion we see that none of the exceptional curves are ($-1$)-curves, thus our resolution is minimal.
\end{proof}

\end{document}